\documentclass{article}
\usepackage{amssymb,latexsym,amsmath}
\usepackage{amsmath,amscd}
\usepackage{graphicx}
\usepackage[all,arc,curve,color,frame,line]{xy}
\usepackage{graphicx}
\usepackage{subfig}
\usepackage{array}
\usepackage{placeins}
\usepackage{adjustbox}
\usepackage{subcaption}

\usepackage{biblatex}
\PassOptionsToPackage{usenames,dvipsnames,svgnames}{xcolor}

\noexpand\extrafloats\RequirePackage{etex}
\usepackage{tikz}
\usetikzlibrary{arrows,shapes,trees,positioning,automata}
\usetikzlibrary{bending}
\usetikzlibrary{patterns,decorations.pathreplacing}
\usepackage{tikz-network}
\usepackage{caption}

\usepackage{wrapfig}
\usepackage[left=1in, right=1in, top=1.25in]{geometry}
\usepackage{color}
\usepackage{amsthm}

\newtheorem{theorem}{Theorem}
\newtheorem{proposition}{Proposition}[theorem]

\theoremstyle{definition}
\newtheorem{example}{Example}[section]

\theoremstyle{definition}
\newtheorem{definition}{Definition}

\theoremstyle{definition}
\newtheorem{question}{Question}[section]
\begin{document}

\title{Methods of constructing spaces with non-trivial self covers}
\author{Mathew Timm\footnote{Department of Mathematics, Bradley University, Peoria, IL 61625, U.S.A. email - mtimm@bradley.edu.  Partially supported by Fulbright U.S. Scholars Grant PS00309131 and the Department of Mathematics, Faculty of Science, University of Split.}}

\date{\today}

\maketitle
\begin{abstract}
\noindent We survey techniques for constructing spaces with non-trivial self covers.  These process include methods for building low and high dimension continua which non-trivially self.  We also discuss several related group theoretic and topological concerns. Statements of a number of open problems related to self covering phenomena are included.
\\

\noindent 2020 AMS and zbMath Mathematics Subject Classification: 54E40, 57M10, 57M15, 57M60, 57N20, 54F50.

\medskip

\end{abstract}
\large

\bigskip
The first examples of covering spaces the typical topology student encounters are those of the circle $S^1$ and it happens that all its finite sheeted covering spaces are itself.  The student then sees self covering phenomena as they develop the classification of compact surfaces. Here they learn that among the infinitely many compact surfaces the only ones which non-trivially cover themselves are the annulus $S^1\times I$, the Mobius band, the 2-torus $T^2=S^1\times S^1$, and the Klein bottle.  A bit later, students might also learn that the products $S^2\times S^1$, $S^1\times S^1 \times I$, the $n$-tori $T^n=\overset{n}{\times}S^1$, and possibly a few more manifolds, have nontrivial self covers.  However, the impression that develops through this initial introduction is that spaces which have nontrivial self covers are quite rare in the universe of all possible spaces.  As a consequence, the special nature of the self covering phenomenon and its rarity within the class of all spaces lends them a mystique which has, over the last 60 years or so, generated a large volume of work.  The purpose of this paper is to survey some of the methods for constructing spaces which non-trivially self covers and collect them in a single place.  By doing so we hope to show that such spaces are, perhaps, not so rare as the casual introduction to them suggests.  Along the way, we illustrate some of the connections between these spaces and other questions in topology and algebra.  There are, without a doubt, examples of, and methods for constructing, spaces with nontrivial self covers other than those included in this survey.  If your favorite method for doing so, or favorite example of a space with a non-trivial self cover is not included here, please let me know so that it can be included in a sequel.  Before turning to the main topic we need a few definitions.

\bigskip
A \emph{covering space} of a topological space $Y$ is a pair $(X,f)$ where $X$ is a topological space and $f:X\rightarrow Y$ is a continuous surjective map such that for each $y\in Y$ there is a neighborhood $V_y$ and a pairwise disjoint collection of neighborhoods $U_x$, $x\in f^{-1}(y)$, for which the restriction $(f|U_x):U_x\rightarrow V_y$ is a homeomorphism.  The function $f$ is said to \emph{evenly cover} $Y$ and it is called a \emph{covering space projection}. The space $X$ is the \emph{total space} of the cover.   We sometimes abbreviate and say $X$ is a covering space of $Y$ without including a reference to the covering projection.  A \emph{deck transformation} for a covering space $f:X \rightarrow Y$ is a homeomorphism $\alpha:X\rightarrow Y$ such that $f\circ \alpha = f$.  The set $Aut_Y(X,f)$ of all deck transformations for the covering space $f:X\rightarrow Y$ is a group under function composition.  It is frequently the case that the following discussion will focus on self covers with particular groups as their group of deck transformations.  Specfically, the covering space $f:X\rightarrow Y$ is \emph{regular}, if $Aut_Y(X,f)$ is ``as big as possible,'' that is, if for every $y\in Y$ and every $x_1,x_2 \in f^{-1}(y)$ there is a deck transformation $\alpha \in Aut_Y(X,f)$ such that $\alpha(x_1)=x_2$.  A regular covering $f:X\rightarrow Y$ is a \emph{$G$-regular} cover, if $Aut_Y(X,f) \cong G$.  Unless indicated otherwise, it is assumed that both domain and range spaces of a covering space are connected.  Consequently, the \emph{degree} or \emph{order} of a covering space $f:X \rightarrow Y$ is the constant $deg(f)=|f^{-1}(y)|$, $y\in Y$.

A subgroup $H$ of a group $G$ is a \emph{clone}, if $H\cong G$.  A group $G$ with a proper clone is \emph{non-cohopfian}.  For spaces, e.g., connected manifolds, with sufficiently nice local topology, if such a space has a nontrivial finite degree self cover, then its fundamental group has a proper clone and given a space $X$ such that $\pi_1(X)$ is non-cohopfian, whether or not $X$ has a self cover corresponding to a given proper clone is a problem of interest.  Note that while our focus will be on finite index clones of a given group, it is possible for a group to have infinite index clones.  For example, in the free group $F_2$ on two generators, any infinite index subgroup $H$ with a minimal generating set consisting of 2 elements is a clone of $F_2$ since all subgroups of a free group are also free.

We make frequent use of the connections between group theory and topology via the fundamental group of a topological space.  Our fundamental groups have base points but, as above, we will usually suppress it in the notation.  Recall that there is a one-to-one correspondence between covering spaces $f:X\rightarrow Y$ and subgroups of the fundamental group $\pi_1(Y,y_0)$ when $Y$ is connected, locally path connected, and semi-locally simply connected.  See Spanier \cite{Spanier1966}.  For the remainder of the paper, unless noted otherwise, a \emph{space} is a metric space.  A \emph{continuum} is a compact connected metric space.

\bigskip
A focus of this paper is methods for building low dimensional spaces, especially methods for build 1 and 2 dimensional spaces with non-trivial self covers. However, to satisfy our objective to collect together methods for building spaces with non-trivial self covers in a single place, we also provide a treatment of several processes which produce high dimensional spaces which non-trivially self cover.  The most important of these are summarized in Section \ref{sec:HighDimSelfCov}, but other methods for producing high dimensional spaces which non-trivially self cover are scattered through out the paper because they are obtained from a ``base space'' which is low dimensional.  Brief discussions of self covering phenomena in dimension 3 are in Sections \ref{sec:SeifertFibred} and \ref{sec:conclussions}.  However, this is a topic of sufficient independent interest that a more detailed discussion is warranted and will be dealt with in detail in another paper.

\bigskip
A topological space can have the property that every finite sheeted cover of the space is a self cover.  There are two ways that this can happen.  One way is for the space $Y$ to have the sort of structure that forces every local homeomorphism $f:X\rightarrow Y$ to be a global homeomorphism.  Questions about these sorts of spaces provided fodder for many conversations with my late colleague Jerry Jungck in regard to his work \cite{JungckLocalHomeo1983}. These conversations wandered through various mathematical terrains and frequently touched on questions about spaces which have the property that every self map has a fixed point (that is, spaces which satisfy the, so called, Fixed Point Property).  Because of the continuing interest in planar continua which have the Fixed Point Property, see Hagopian \cite{HAGOPIAN2007}, these conversations led to investigations of special properties of subcontinua of the plane \cite{JungckTimm1999}, \cite{JungckTimm2002}, and eventually to the development of the construction process described in Section \ref{sec:CantorSetControl}.  This process produces large numbers of planar, 1-, and 2-dimensional continua which self cover.  It is based on two observations.  First, that the (middle thirds) Cantor set has $G$-regular self covers for every finite group $G$.  Second, that the Cantor set can act as a ``control set'' for the construction of larger topological spaces which have $G$-regular covers for many groups $G$.  One of the interesting consequences of these developments is actually a question about the Hawaiian Earring $\mathcal{H}$, a space which has no nontrivial self covers.  See Section \ref{sec:HawaiianEarRing} for specifics and \cite{Timm2002Survey} for additional discussion.

\bigskip
The second process for producing low dimensional spaces with non-trivial self covers is an inverse system, shape theoretic, based approach.  This method has its origins in work by Fox \cite{Fox1972},\cite{Fox1974} on the solenoids, that is, inverse limits $\Sigma=\underset{\leftarrow}{\textrm{\textbf{lim}}}(S^1,f_n,\mathbb{N})$ of inverse sequences in which the spaces in the inverse sequene are all copies of the circle $S^1$ and all the bonding maps are finite sheeted self covers $f_n:S^1\rightarrow S^1$.  This inverse limit process produces limit spaces $\Sigma$ with complicated local topology.  Specifically, the $\Sigma$ are not locally connected.  They also have the property that their only subcontinua are copies of the real line, and consequently, the $\Sigma$'s have trivial fundamental group.  However, in spite of the triviality of their fundamental groups, they have interesting covering spaces properties.  For example, the only covering spaces of any particular $\Sigma$ are itself and the types of self coverings desired one wants in a particular $\Sigma$ can be engineered by careful selections of the degrees of the self covers $f_n:S^1\rightarrow S^1$ used as the bonding maps in the inverse sequence.  By doing so, one can produce solenoids whose only covering spaces are trivial ones (i.e., homeomorphisms) or one can select bonding maps in the sequence so that the solenoid has non-trivial self covers of almost every degree $d\in \mathbb{N}$.  In particular, one can predetermine a finite or infinite collection $\mathcal{E} \neq \emptyset$ of degrees $m \geq 2$ of self covers of a solenoid which one wants to \underline{exclude}, and produces a solenoid with self covers of every degree $d$ relatively prime to every $m\in \mathbb{N} \setminus \mathcal{E}$. When $\mathcal{E} = \mathbb{N} \setminus \{1\}$, the resulting solenoid's only covers are self homeomorphisms.  Note that it is not possible to choose $\mathcal{E}=\mathbb{N}$ since the identity map on any topological space is a self cover.  See Section \ref{sec:InvLimAndSelfCov}.

\bigskip
The third method for constructing low dimensional spaces which nontrivially self cover is modeled on the graph of groups approach to building a group.   It makes use of a directed graph $\Gamma$ and a weight function $\omega$ which assigns to each directed edge $e$ of $G$ an ordered pair $\omega(e)=(\omega^-(e),\omega^+(e))$ of non-zero integers.  The pair $(\Gamma,\omega)$ provides instructions on how to glue copies of the annulus $A=S^1\times [0,1]$ together along their boundaries.  The resulting 2-complexes, $K(\Gamma,\omega)$ then have $k$-fold self covers for every $k$ which is relatively prime to all the weights on $\Gamma$.  Their fundamental groups are members of the class of groups called \emph{generalized Bausmslag-Solitar groups} and consequently, it follows that every finitely related generalized Baumslag-Solitar group has proper finite index clones. See Section \ref{sec:GraphsOfSpaces}.

\bigskip
Section \ref{sec:DynamSystems} considers self covering phenomena in the context of complex dynamical systems and a class of groups known as a self similar groups.  Specifically, the Julia sets $J(f)$ of certain holomorphic functions $f:\widehat{\mathbb{C}} \rightarrow \widehat{\mathbb{C}}$ of the Riemann sphere $\widehat{\mathbb{C}}=\mathbb{C}\cup +\infty$ non-trivially self cover and these spaces $J(f)$ can sometimes be obtained from a self similar group acting on a rooted tree.

\bigskip
Section \ref{sec:GeoRep} contains a few comments about geometric representability of group theoretic properties within various classes of topological spaces. They build on the observation that spaces with non-trivial self covers provide for the geometric representability of the non-cohopfianism of various groups.  A few concluding remarks are in Section \ref{sec:conclussions}.   .

\bigskip
And finally, before turning to more detailed discussions of self covering phenomena, we provide a short list of a familiar spaces which have self covers.

\begin{itemize}
    \item As noted above, the circle, $S^1 = \{e^{i\cdot \theta}:0\leq \theta \leq 2 \pi\}$ has $\mathbb{Z}_m$-regular self covers for every $m\in \mathbb{N}$.  They are given by $f(e^{i\cdot \theta})=e^{i\cdot m \theta}$.
    \item The torus $S^1\times S^1$, the Klein bottle, the annulus, the Mobius band, and the products $S^1 \times B^n$, $B^n = \{(x_1,\ldots , x_n)\in \mathbb{R}^n:0\leq \sqrt{x^2_1+\cdots +x^2_n} \leq 1\}$, of the circle with the $n$-ball, $n\in \mathbb{N}$, has $\mathbb{Z}_m$-regular self covers for every $m\in \mathbb{N}$.
    \item Surprisingly, since it is a non-homogeneous continuum, David Bellamy \cites{Bellamy2005} has shown that the pseudo-circle has regular cyclic self covers of all orders.
    \item There are multiple ways to prove that the Cantor set has $G$-regular self covers for every finite group $G$.  One is given below, but the reader is invited to find as many other ways to do so as possible.
\end{itemize}

\section{High dimensional self covering}\label{sec:HighDimSelfCov}
First, consider the $n$-torus $T^n=\overset{n}{\underset{i=1}{\times}}S^1$.  All of its finite sheeted covering spaces are self covers.  Hence, for all $n\in \mathbb{N}$ and, in particular, for $n\geq 4$, we have a product with a great many inequivalent self covers.\footnote{Two covers $(X_1,f_1)$ and $(X_2,f_2)$ of a space $Y$ are \emph{equivalent} if there is homeomorphism $h:X_1\rightarrow X_2$ such that $h\circ f_2 = f_1$.}  Specifically, given any  finite abelian group $A$ on $1\leq m \leq n$ generators with $rank(A)=m$, $T^n$ has an $A$-regular self cover.

Next, consider $T^n$ and take any connected $m$-dimensional space $M$.  Then the product $T^n\times M$ has, for each self cover $f:T^n\rightarrow T^n$, a self cover $\widetilde{f}:T^n\times M \rightarrow T^n\times M$ defined by $\widetilde{f}(x,y)=(f(x),y)$.  In addition, $T^n\times M$, like $T^n$ itself, has $A$-regular self covers for every finite abelian group $A$ on at most $n$ generators.  If $M$ also has self covers, there are other self covers of the product induced by those of $M$.  When either $n$ or $m$, or both, are sufficiently large, these spaces are high dimensional.

As a slight generalization of products, assume that $M$ has a periodic homeomorphism $h:M\rightarrow M$ of order $k$.  Then the locally trivial $M$ bundle over $S^1$, $M\times I / h$ obtained by gluing $M\times 0$ to $M\times 1$ via $h$, has a degree $d = (k+1)$ self cover. Thus, when $m$ is large, this bundle construction process produces high dimensional examples of spaces which nontrivially cover themselves.  Of course, the bundle construction can also be used to produce low dimensional spaces which self cover.  In dimension 2, the Klein bottle is constructed by taking $M=S^1$ and $h$ as the antipodal map.  For a 3-dimensional example, let $M$ be the genus 2 surface and $h:M\rightarrow M$ a period 2 homeomorphism with exactly two fixed points, e.g., one that interchanges the holes.  Both of these bundles have $\mathbb{Z}_{2m+1}$-regular self covers for all $m\in \mathbb{Z}$.

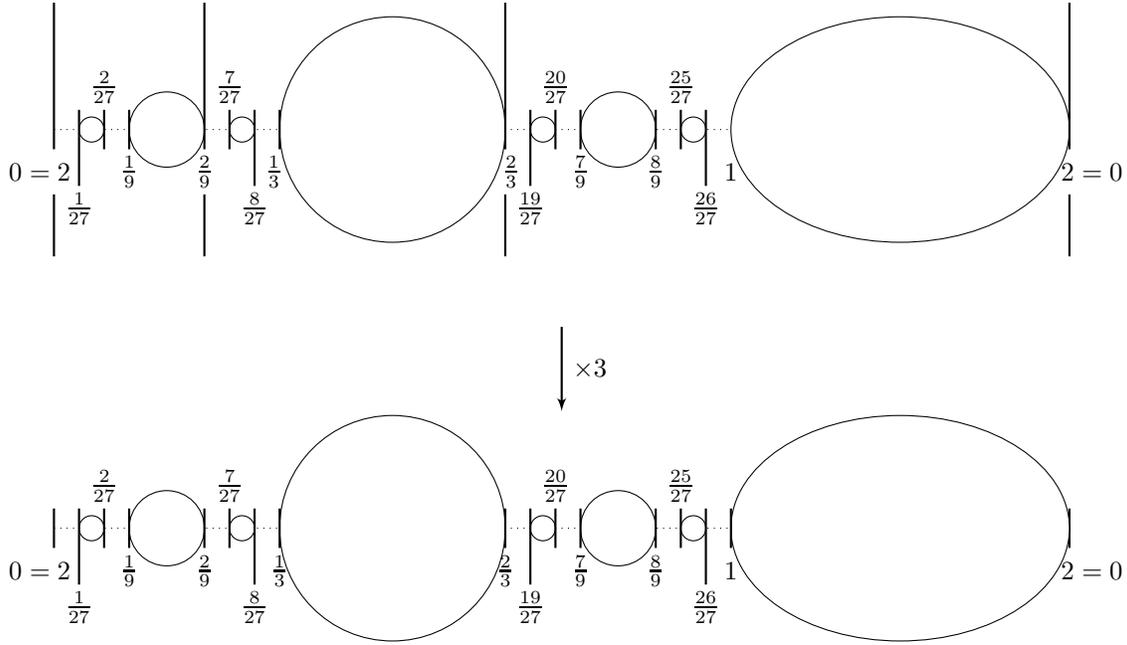
\begin{figure}[t]
\centering
\begin{tikzpicture}[scale=.75]
\node[fill=white] (2) at (12.4,-.75) {$2 = 0$};
\draw (9,0) ellipse (3cm and 2cm);
\draw[dotted] (-6,0) -- (6,0);
\draw[thick] (12,-.35) -- (12,2.25); \draw[thick] (12,-2.25)--(12,-1.15); \node (1) at (6,-.75) {$1$};
\draw[thick] (-6,-.35) -- (-6,2.25); \node[fill=white] (0) at (-6.25,-.75) {$0=2$}; \draw[thick] (-6,-2.25)--(-6,-1.15);
\draw[fill=white] (0,0) circle (2cm);
\draw[thick] (2,-.35) -- (2,2.25); \node (2/3) at (2.1,-.75) {$\frac{2}{3}$}; \draw[thick] (2,-2.25) -- (2,-1.15);
\draw[thick] (-2,-.35) -- (-2,.35); \node(1/3) at (-2.1,-.75) {$\frac{1}{3}$};
\draw[fill=white] (4,0) circle (.667cm);
\draw[thick] (4.6667,-.35) -- (4.6667,.35);
\node (8/9) at (4.6667,-.75) {$\frac{8}{9}$};
\draw[thick] (3.3333,-.35) -- (3.3333,.35); \node (7/9) at (3.3333,-.75) {$\frac{7}{9}$};
\draw[fill=white] (-4,0) circle (.667cm);
\draw[thick] (-4.6667,-.35) -- (-4.6667,.35); \node (1/9) at (-4.6667,-.75) {$\frac{1}{9}$};
\draw[thick] (-3.3333,-.35) -- (-3.3333,2.25); \node (2/9) at (-3.3333,-.75) {$\frac{2}{9}$}; \draw[thick] (-3.3333,-2.25) -- (-3.3333,-1.15);
\draw[fill=white] (2.66665,0) circle (.2222);
\draw[thick] (2.44443,-1) -- (2.44443,.35); \node (19/27) at (2.44443,-1.4) {$\frac{19}{27}$};
\draw[thick] (2.88887,-.35) -- (2.88887,.35); \node (20/27) at (2.88887,.75) {$\frac{20}{27}$};
\draw[fill=white] (5.33335,0) circle (.2222);
\draw[thick] (5.11113,-.35) -- (5.11113,.35); \node (25/27) at (5.11113,.75) {$\frac{25}{27}$};
\draw[thick] (5.55557,-1) -- (5.55557,.35); \node (26/27) at (5.55557,-1.4) {$\frac{26}{27}$};
\draw[fill=white] (-2.66665,0) circle (.2222);
\draw[thick] (-2.44443,-1) -- (-2.44443,.35); \node (8/27) at (-2.44443,-1.4) {$\frac{8}{27}$};
\draw[thick] (-2.88887,-.35) -- (-2.88887,.35); \node (7/27) at (-2.88887,.75) {$\frac{7}{27}$};
\draw[fill=white] (-5.33335,0) circle (.2222);
\draw[thick] (-5.11113,-.35) -- (-5.11113,.35); \node (2/27) at (-5.11113,.75) {$\frac{2}{27}$};
\draw[thick] (-5.55557,-1) -- (-5.55557,.35); \node (1/27) at (-5.55557,-1.4) {$\frac{1}{27}$};
\tikzset{edge/.style = {->,> = latex'}}
\draw[edge,thick] (3,-3.5) to (3,-5);
\node (x3) at (3.5,-4.25) {$\times 3$};
\end{tikzpicture}

\begin{tikzpicture}[scale=.75]
\node[fill=white] (2) at (12.4,-.75) {$2 = 0$}; \draw[thick] (12,-.35)--(12,.35);
\draw[fill=white] (9,0) ellipse (3cm and 2cm);
\draw[dotted] (-6,0) -- (6,0);
\draw[thick] (6,-.35) -- (6,.35); \node (1) at (6,-.75) {$1$};
\draw[thick] (-6,-.35) -- (-6,.35); \node[fill=white] (0) at (-6.25,-.75) {$0=2$};
\draw[fill=white] (0,0) circle (2cm);
\draw[thick] (2,-.35) -- (2,.35); \node (2/3) at (2,-.75) {$\frac{2}{3}$};
\draw[thick] (-2,-.35) -- (-2,.35); \node(1/3) at (-2,-.75) {$\frac{1}{3}$};
\draw[fill=white] (4,0) circle (.667cm);
\draw[thick] (4.6667,-.35) -- (4.6667,.35);
\node (8/9) at (4.6667,-.75) {$\frac{8}{9}$};
\draw[thick] (3.3333,-.35) -- (3.3333,.35); \node (7/9) at (3.3333,-.75) {$\frac{7}{9}$};
\draw[fill=white] (-4,0) circle (.667cm);
\draw[thick] (-4.6667,-.35) -- (-4.6667,.35); \node (1/9) at (-4.6667,-.75) {$\frac{1}{9}$};
\draw[thick] (-3.3333,-.35) -- (-3.3333,.35); \node (2/9) at (-3.3333,-.75) {$\frac{2}{9}$};
\draw[fill=white] (2.66665,0) circle (.2222);
\draw[thick] (2.44443,-1) -- (2.44443,.35); \node (19/27) at (2.44443,-1.4) {$\frac{19}{27}$};
\draw[thick] (2.88887,-.35) -- (2.88887,.35); \node (20/27) at (2.88887,.75) {$\frac{20}{27}$};
\draw[fill=white] (5.33335,0) circle (.2222);
\draw[thick] (5.11113,-.35) -- (5.11113,.35); \node (25/27) at (5.11113,.75) {$\frac{25}{27}$};
\draw[thick] (5.55557,-1) -- (5.55557,.35); \node (26/27) at (5.55557,-1.4) {$\frac{26}{27}$};
\draw[fill=white] (-2.66665,0) circle (.2222);
\draw[thick] (-2.44443,-1) -- (-2.44443,.35); \node (8/27) at (-2.44443,-1.4) {$\frac{8}{27}$};
\draw[thick] (-2.88887,-.35) -- (-2.88887,.35); \node (7/27) at (-2.88887,.75) {$\frac{7}{27}$};
\draw[fill=white] (-5.33335,0) circle (.2222);
\draw[thick] (-5.11113,-.35) -- (-5.11113,.35); \node (2/27) at (-5.11113,.75) {$\frac{2}{27}$};
\draw[thick] (-5.55557,-1) -- (-5.55557,.35); \node (1/27) at (-5.55557,-1.4) {$\frac{1}{27}$};
\end{tikzpicture}
\caption{The necklace of circles $\mathcal{N}^{(1)}_{C,S^1}$ and its $\mathbb{Z}_3$-regular cover $f_3:\mathcal{N}^{(1)}_{C,S^1}\rightarrow \mathcal{N}^{(1)}_{C,S^1}$. The vertical line segments in the picture corresponding to the rational points of the Cantor set are not part of the necklace.  They are included only for reference.  The lifts $f^{-1}_3(0=2)$ in the domain copy are indicated by the longest vertical line segments in the domain copy of $\mathcal{N}^{(1)}_{C,S^1}$.  These line segments delineate closures of fundamental domains of the covering map.  Note that on each of these fundamental domains $f_3$ is non-contractive.  In fact, it is expansive on the part of the domain copy of $\mathcal{N}^{(1)}_{C,S^1}$ between $0=2$ and $1$ in the picture and maps the ellipse sitting between $1$ and $2=0$ homeomorphically onto itself.} \label{fig:necklace_of_circles}
\end{figure}

Infinite dimensional products which clearly have nontrivial self covers include the infinite dimensional torus $T^{\infty}=\underset{i=1}{\overset{\infty}{\times}}S^1$ with the product topology and $T^{\infty}\times M$, in the product topology, for any topological space $M$.  In addition, a ``categorical'' process \cite{DelgadoTimm2017} allows one to start with any finite group $G$ and use it to build infinite dimensional product spaces with $H$-regular self covers for each proper subgroup $H$ of $G$.  We illustrate this categorical process with an example.

\bigskip
Consider the finite cyclic group $\mathbb{Z}_6$ (written additively).  It has 4 subgroups, $\mathbb{Z}_6$, $\langle 2 \rangle$, $\langle 3 \rangle$, and $0$.  Build the standard 2-complex, the so called \emph{presentation complex} $K$ for $\mathbb{Z}_6$ with $\pi_1(K) \cong \mathbb{Z}_6$.  It has one 0-cell; one 1-cell, with both endpoints attached to the 0-cell; and one 2-cell glued to the circle $S^1$, formed by the 0- and 1-cell, via the degree 6 cover $f:\partial D\rightarrow S^1$. For each subgroup $H$ of $\mathbb{Z}_6$, let $f_H:K_H\rightarrow K$ be the $H$-regular covering space of $K$ and let $S$ be the set of all subgroups of $\mathbb{Z}_6$.  Note that $K_{0}$ is a copy of $K$ and $K_{\mathbb{Z}_6}$ is the universal cover of $K$.  Take the product $\underset{H\in S}{\times}K_H$, then form the infinite product $X=\overset{\infty}{\underset{i=1}{\times}}\left(\underset{H\in S}{\times}K_H\right)$.  Now consider the map $\widetilde{f}:K_{\langle 2 \rangle}\times X \rightarrow K \times X$ defined by $\widetilde{f}(t,x) = (f_{\langle 2 \rangle}(t),x)$.  Both $K_{\langle 2 \rangle} \times X$ and $K\times X$ are homeomorphic to $X$ -- the homeomorphisms are just coordinate changes -- and the composition of these coordinate changes with $\widetilde{f}$ is a $\mathbb{Z}_2$-regular self cover of $X$.  Similarly $X$ has $H$-regular self covers for each $H < \mathbb{Z}_6$.  See \cite{DelgadoTimm2017} for a complete description of this process.

\bigskip
Other high dimensional examples of the self covering phenomena can be produced via a second categorical process called the Tinsley Construction.  See Daverman \cite{Daverman1993HyperhopfianGA} for more details.  Start with an aspherical group $G$ which has a proper finite index subgroup $H$ isomorphic to $G$ and and which satisfies additional technical conditions.  Assume also, that $K$ is a finite aspherical complex with $\pi_1(K)\cong G$.  Let $f:K_H \rightarrow K$ be the finite sheeted cover of $K$ corresponding to $H$.  Since $K$ is aspherical, so is $K_H$ and, consequently, $K$ and $K_H$ are homotopy equivalent.  The technical conditions then imply that there is a complex $L$ which collapses to both $K$ and $K_H$ and, therefore $\pi_1(L)\cong G$. Embed $L$ in a copy of Euclidean space of dimension so large that a regular neighborhood $N$ of $L$ has $dim(\partial N) \geq dim(L)+2$.  Then $\pi_1(\partial N)\cong G$, $\partial N$ bounds a regular neighborhood of both $K$ and $K^{\prime}$ and, in general, $\partial N$, nontrivially covers itself.  For example, the (soluble) Baumslag-Solitar groups, $G\cong \langle x,t:t^{-1}x^{-1}t=x^{m-1} \rangle$, $m\geq 2$ satisfies the needed conditions and, consequently, applying the Tinsley construction to each of them provide high dimensional examples of manifolds which nontrivially cover themselves.  Since $m\geq 2$, these manifolds are non-products.

\section{Spaces whose only finite sheeted covers are themselves}\label{sec:OnlyFinCovSelf}
A connected space $Y$ can have the property that every covering space map $f:X \rightarrow Y$ from a connected space $X$ onto $Y$ be a self cover by virtue of having no non-trivial covers at all.  That is, the connected space $X$ satisfies the condition
$$
\begin{array}{c}
  (\star) \textrm{ if $X$ is connected and $f:X\rightarrow Y$ is covering space map, then $f$ is a homeomorphism.}
\end{array}
$$
\noindent
Connected, simply connected manifolds $Y$ satisfy $(\star)$.  The $(\star)$-condition is a special case of the more general problem of determining conditions on the space $Y$ such that every local homeomorphism $f:X \overset{\textrm{onto}}{\longrightarrow} Y$ is a global homeomorphism.  See Jungck \cite{JungckLocalHomeo1983} and \cite{Timm2002Survey} for more detailed discussions.  This also relates to the problem of determining which functions on which spaces have fixed points.

\bigskip
Recall that a \emph{fixed point} for a function $f:X\rightarrow X$ is a point $x_0\in X$ such that $f(x_0)=x_0$.  A space $X$ which is such that every continuous self map $f:X\rightarrow X$ has a fixed point $f(x_0)=x_0$ is said to \emph{have the Fixed Point Property}.  Questions about which functions on which spaces have fixed points and questions about which spaces have the Fixed Point Property have been of interest for a very long time.  For example, Brower proved his Fixed Point Theorem in 1911 \cite{Brouwer1911}, \cite{wikiBrowerFpp}.  Two continua of interest which have the fixed point property are the pseudo arc and the Topologist's Sine Curve $X = \{(x,y):sin(\frac{1}{x}: 0<x\frac{1}{\pi} \} \cup \{(0,y): -1 \leq y \leq 1\}$. See Deyer \cite{Hamilton} and Hagopian \cite{hagopian1976fixed}.   Additional history of the investigations into which planar continua have the Fixed Point Property are in Hagopian's other papers \cite{HagopianMarsh}, \cite{HAGOPIAN2007}, \cite{Hagopian1996TheFP}.

When a space has the Fixed Point Property, it provides an obstruction to regular self covering as follows.  Suppose $X$ has the Fixed Point Property.  In addition, suppose $f:X\rightarrow X$ is a self cover, and $\alpha:X\rightarrow X$ is a deck transformation for the self cover.  The fact that $X$ has the Fixed Point Property implies that $\alpha$ is the identity on $X$.  Consequently, if $X$ has sufficiently nice local topology, e.g., $X$ is a manifold or locally finite cell complex, then $f:X\rightarrow X$ is not a regular cover unless it is the universal cover.  In which case, it follows that $\pi_1(X)$ is the trivial group.  This observation suggests several lines of investigation.  One can ask whether every simply connected space has the Fixed Point Property (the antipodal map on the 2-sphere shows this is not the case) or, more generally, about the topology of spaces with fundamental groups which are simple groups.

\bigskip
The second way a space can have all its finite sheeted covering spaces homeomorphic to itself is to have at least one degree $d\geq 2$ covering space and also satisfy the condition that all its finite sheeted covers are self covers.  The $n$-tori, for example, satisfy this condition.  When $X$ is a finite dimensional manifold or a locally finite cell complexes which has a nontrivial cover and it has the property that all of its finite sheeted covers are self covers, it follows that its fundamental group $\pi_1(M)$ has a proper finite clone and that every finite index subgroup is a clone.  This generates the obvious problem of determining the structure of non-cohopfian groups which have the property that all of their finite index subgroups are clones.  A structure theorem for these sorts of finitely generated, highly non-cohopfian groups is in \cite{RobinsonTimm1998}.  Whether or not any of these groups satisfy the conditions needed to apply the Tinsley construction is, as of this writing, unknown.  If any one of them does, the construction can be applied to produce a high dimensional manifold which has a self cover corresponding to a particular finite index clone.  However, which of these groups can be fundamental groups of spaces for which \underline{every} finite degree cover is a self cover is an interesting question.  An example of how this sort of group theoretic condition determines the structure of finite 2- and 3-dimensional cell complex which non-trivially self cover can be found in \cite{DelgadoTimm2003}.
\bigskip

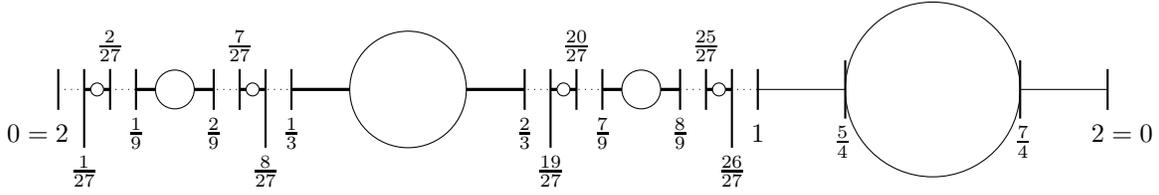
\begin{figure}
\centering
\begin{tikzpicture}[scale=.775]
\draw[thick] (12,-.35)--(12,.35); \node[fill=white] (2) at (12.25,-.75) {$2=0$};
\draw[dotted] (-6,0) -- (6,0);
\draw[thick] (6,-.35) -- (6,.35); \node (1) at (6,-.75) {$1$};
\draw[thick] (-6,-.35) -- (-6,.35); \node[fill=white] (0) at (-6.35,-.75) {$0=2$};
\draw (6,0)--(12,0); \draw[fill=white] (9,0) circle (1.5cm);
\draw[thick] (7.5,-.5)--(7.5,.5);\draw[thick] (10.5,-.5)--(10.5,.5);
\node (5/4) at (7.45,-.9) {$\frac{5}{4}$}; \node (7/4) at (10.55,-.9) {$\frac{7}{4}$};
\draw[fill=white] (0,0) circle (1cm);
\draw[very thick] (1,0) -- (2,0); \draw[thick] (2,-.35) -- (2,.35); \node (2/3) at (2,-.75) {$\frac{2}{3}$};
\draw[very thick] (-1,0) -- (-2,0); \draw[thick] (-2,-.35) -- (-2,.35); \node(1/3) at (-2,-.75) {$\frac{1}{3}$};
\draw[fill=white] (4,0) circle (.333cm);
\draw[very thick] (4.333,0) -- (4.6667,0); \draw[thick] (4.6667,-.35) -- (4.6667,.35); \node (8/9) at (4.6667,-.75) {$\frac{8}{9}$};
\draw[very thick] (3.6667,0) -- (3.3333,0); \draw[thick] (3.3333,-.35) -- (3.3333,.35); \node (7/9) at (3.3333,-.75) {$\frac{7}{9}$};
\draw[fill=white] (-4,0) circle (.333cm);
\draw[very thick] (-4.333,0) -- (-4.6667,0); \draw[thick] (-4.6667,-.35) -- (-4.6667,.35); \node (1/9) at (-4.6667,-.75) {$\frac{1}{9}$};
\draw[very thick] (-3.6667,0) -- (-3.3333,0); \draw[thick] (-3.3333,-.35) -- (-3.3333,.35); \node (2/9) at (-3.3333,-.75) {$\frac{2}{9}$};
\draw[fill=white] (2.66665,0) circle (.1111); \draw[very thick] (2.44443,0) -- (2.55554,0); \draw[thick] (2.44443,-1) -- (2.44443,.35); \node (19/27) at (2.44443,-1.4) {$\frac{19}{27}$};
\draw[very thick] (2.77776,0) -- (2.88887,0); \draw[thick] (2.88887,-.35) -- (2.88887,.35); \node (20/27) at (2.88887,.75) {$\frac{20}{27}$};
\draw[fill=white] (5.33335,0) circle (.1111); \draw[very thick] (5.11113,0) -- (5.22224,0); \draw[thick] (5.11113,-.35) -- (5.11113,.35); \node (25/27) at (5.11113,.75) {$\frac{25}{27}$};
\draw[very thick] (5.44446,0) -- (5.55557,0); \draw[thick] (5.55557,-1) -- (5.55557,.35); \node (26/27) at (5.55557,-1.4) {$\frac{26}{27}$};
\draw[fill=white] (-2.66665,0) circle (.1111); \draw[very thick] (-2.44443,0) -- (-2.55554,0); \draw[thick] (-2.44443,-1) -- (-2.44443,.35); \node (8/27) at (-2.44443,-1.4) {$\frac{8}{27}$};
\draw[very thick] (-2.77776,0) -- (-2.88887,0); \draw[thick] (-2.88887,-.35) -- (-2.88887,.35); \node (7/27) at (-2.88887,.75) {$\frac{7}{27}$};
\draw[fill=white] (-5.33335,0) circle (.1111); \draw[very thick] (-5.11113,0) -- (-5.22224,0); \draw[thick] (-5.11113,-.35) -- (-5.11113,.35); \node (2/27) at (-5.11113,.75) {$\frac{2}{27}$};
\draw[very thick] (-5.44446,0) -- (-5.55557,0); \draw[thick] (-5.55557,-1) -- (-5.55557,.35); \node (1/27) at (-5.55557,-1.4) {$\frac{1}{27}$};
\end{tikzpicture}
\caption{The 1-dimensional continuum $\mathcal{N}^{(2)}_{C,S^1}$ with $\mathbb{Z}_m$-regular self covers for all $m\in \mathbb{Z}.$  The small vertical line segments corresponding to the rational points of the Cantor set, the one at $\frac{5}{4}$, and the one at $\frac{7}{4}$ are not part of $\mathcal{N}^{(2)}_{C,S^1}$.  They are included only to indicate the location of the associated rational points.  Should the reader be inclined to include them, doing so produces another continuum with $\mathbb{Z}_m$-regular self covers for every $m \in \mathbb{N}$.  In this case the control set is $C\times I$.}
\label{fig:2nd_necklace_of_circles}
\end{figure}

\section{Cantor set controlled constuctions} \label{sec:CantorSetControl}
Let $f:X \rightarrow X$ be a self cover of the space $X$ and $x_0 \in X$.  We have $f^{-1}(x_0)= \{x\in X: f(x)=x_0\}$, and for $n\in \mathbb{N}$, set $f^{-(n+1)}(x_0)=\{x\in X:f(x)\in f^{-n}(x_0)\}$.  If $deg(f)=d \geq 2$, then $|f^{-n}(x_0)|= d^n$ and consequently $|\overset{\infty}{\underset{n=1}{\bigcup}}f^{-n}(x_0)| = \infty$.  Now, suppose that an observer of this situation notices that $X$ is a continuum and that the local topology about $x_0 \in X$ appears to be  ``exceptional'' or ``singular.''  Then, because each $n$-fold compositions $f^n=f\circ f\circ \cdots \circ f:X \rightarrow X$ is a degree $d^n$ self cover, the above observation implies that there are infinitely many points about which the local topology is the same as is that about $x_0$.  But via the Bolzano-Weierstrauss Theorem, compactness of $X$ imposes additional structure on $S=\overset{\infty}{\underset{n=1}{\bigcup}}f^{-n}(x_0)$.  In particular, Bolzano-Weierstrauss implies that this set $S$ of exceptional points has a limit point $x^{\prime}$. That is, there is a second type of exceptional points in $X$.  Considering $S^{\prime}=\overset{\infty}{\underset{n=1}{\bigcup}}f^{-n}(x^{\prime})$, Bolzano-Weierstrauss implies that the set $S^{\prime}$ of exceptional points has a limit point and that there is a then second class $S^{\prime \prime}$ of exceptional points in $X$.  This analysis can be repeated, inductively, to get sets of points which are limit points of limit points of limit points of $\ldots \ldots$ and therefore exceptional in some sense.  One space this observation brings to mind is the Cantor set.\footnote{There are other possible structures for the set of exceptional points if there are any.  It could be dense in $X$, and therefore not exceptional or see, for example, Section \ref{sec:GraphsOfSpaces}.  The sets of exceptional points for the spaces in these examples are unions of circles.}  This suggests that building continua $X$ which make the points of a Cantor set $C \subset X$ ``exceptional'' points in $X$ can produce spaces with non-trivial self covers.

Also note that the Cantor set has $G$-regular self covers for every finite group $G$.  To see this, suppose that $G$ has the presentation $\langle x_1,x_2,\ldots,x_m:r_1=1,r_2=1,\ldots,r_n=1\rangle$ with $m$ generators and $n$ relations.  Build the presentation complex $K$ with fundamental group $G = \pi_1(K,y_0)$. Also, let $C$ be the middle thirds Cantor set sitting in a copy $I=[0,1]$ of the unit interval.

Now $K$ has one 0-cell, $y_0$; $m$ 1-cells, $J_1=[0,1],J_2=[0,1],\ldots,J_m=[0,1]$ attached to $y_0$ via their endpoints; and $n$ 2-cells, attached to the bouquet of $m$ loops formed from $y_0$ and the $J_k$, $k=1,\ldots,m$, via gluing maps determined by the relations $r_1=1,r_2=1,\ldots,r_n=1$.   Identify a copy of $I$ with $[\frac{1}{4},\frac{3}{4}] \subset J_1 \subset K$.  This embeds a copy $C_1$ of the Cantor Set $C$ in $J_1$.  The universal cover $f:U\rightarrow K$ is a $G=\pi_1(K,y_0)$-regular cover of $K$.  Because $C_1$ is contained in the open interval $J_1\setminus \{0,1\} \subset K$, and open intervals are simply connected, $C_1$ is contained in a fundamental domain of the covering projection $f:U\rightarrow K$. Consequently, $f^{-1}(C_1)$ is a union of $|G|$ Cantor sets and, since a finite union of Cantor sets is a Cantor set, $f^{-1}(C_1)$ is a Cantor set.  Therefore, the restriction, $f:f^{-1}(C_1)\rightarrow C_l$, is the desired $G$-regular self cover.

\begin{figure}
\centering
\begin{tikzpicture}[scale=.775]
\draw[thick,->, > = latex'] (12,-1.5)--(12,1.5);  \node[fill=white] (2) at (12.25,-.75) {$2=0$};
\draw[dotted] (-6,0) -- (6,0);
\draw[thick] (6,-.35) -- (6,1.5); \draw[thick] (6,-1.5)--(6,-.75); \node[fill=white] (1) at (6,-.75) {$1$};
\draw[thick,->, > = latex'] (-6,-1.5) -- (-6,1.5); \node[fill=white] (0) at (-6.35,-.75) {$0=2$};
\draw (6,0)--(12,0); \draw[fill=white] (9,0) circle (1.5cm);
\node (5/4) at (7.45,-.9) {$\frac{5}{4}$}; \node (7/4) at (10.55,-.9) {$\frac{7}{4}$};
\draw[fill=white] (0,0) ellipse (1cm and 1.5cm);
\draw[very thick] (1,0) -- (2,0); \draw[thick] (2,-1.5) -- (2,1.5); \node[fill=white] (2/3) at (2,-.75) {$\frac{2}{3}$};
\draw[very thick] (-1,0) -- (-2,0); \draw[thick] (-2,-1.5) -- (-2,1.5); \node[fill=white] (1/3) at (-2,-.75) {$\frac{1}{3}$};
\draw[fill=white] (4,0) ellipse (.333cm and 1.5cm);
\draw[very thick] (4.333,0) -- (4.6667,0); \draw[thick] (4.6667,-1.5) -- (4.6667,1.5); \node[fill=white] (8/9) at (4.6667,-.75) {$\frac{8}{9}$};
\draw[very thick] (3.6667,0) -- (3.3333,0); \draw[thick] (3.3333,-1.5) -- (3.3333,1.5); \node[fill=white] (7/9) at (3.3333,-.75) {$\frac{7}{9}$};
\draw[fill=white] (-4,0) ellipse (.333cm and 1.5cm);
\draw[very thick] (-4.333,0) -- (-4.6667,0); \draw[thick] (-4.6667,-1.5) -- (-4.6667,1.5); \node[fill=white] (1/9) at (-4.6667,-.75) {$\frac{1}{9}$};
\draw[very thick] (-3.6667,0) -- (-3.3333,0); \draw[thick] (-3.3333,-1.5) -- (-3.3333,1.5); \node[fill=white] (2/9) at (-3.3333,-.75) {$\frac{2}{9}$};
\draw[fill=white] (2.66665,0) ellipse (.1111cm and 1.5cm);
\draw[very thick] (2.44443,0) -- (2.55554,0); \draw[thick] (2.44443,-1.5) -- (2.44443,1.5); 
\draw[very thick] (2.77776,0) -- (2.88887,0); \draw[thick] (2.88887,-1.5) -- (2.88887,1.5); 
\draw[fill=white] (5.33335,0) ellipse (.1111cm and 1.5cm);
\draw[very thick] (5.11113,0) -- (5.22224,0); \draw[thick] (5.11113,-1.5) -- (5.11113,1.5); 
\draw[very thick] (5.44446,0) -- (5.55557,0); \draw[thick] (5.55557,-1.5) -- (5.55557,1.5); 
\draw[fill=white] (-2.66665,0) ellipse (.1111cm and 1.5cm);
\draw[very thick] (-2.44443,0) -- (-2.55554,0); \draw[thick] (-2.44443,-1.5) -- (-2.44443,1.5); 
\draw[very thick] (-2.77776,0) -- (-2.88887,0); \draw[thick] (-2.88887,-1.5) -- (-2.88887,1.5); 
\draw[fill=white] (-5.3335,0) ellipse (.1111cm and 1.5cm);
\draw[very thick] (-5.11113,0) -- (-5.22224,0); \draw[thick] (-5.11113,-1.5) -- (-5.11113,1.5); 
\draw[very thick] (-5.44446,0) -- (-5.55557,0); \draw[thick] (-5.55557,-1.5) -- (-5.55557,1.5); 
\end{tikzpicture}
\caption{ The continuum $\mathcal{N}^{(3)}_{C\times [-1,1],S^1}.$ It has a $\mathbb{Z}_m$-regular self covers for each $m\in\mathbb{Z}$.  The $x$-coordinates at the points $\frac{i}{27}$ have been removed to reduce clutter in the picture.  The control set in this case is the product $C\times [-1,1]$. The direction of the gluing on $0\times [-1,1]$ and $2\times [-1,1]$ are indicated by the arrows in the picture.} \label{fig:3rd_necklace_of_circles}
\end{figure}
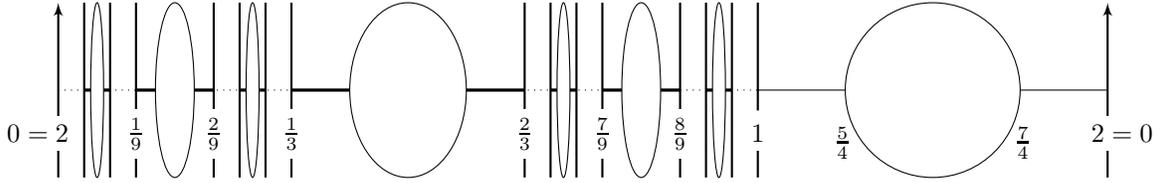

\bigskip
In this section we look at ways in spaces with nontrivial self covers can be obtained by using the self covers of the Cantor set to exercise control over the construction. There is a technical definition of ``control set,'' see \cite{Timm2020} for details, but the idea is that one specifies a space $F$ in advance which has some desired property, then one uses it as a frame on which to construct a space $X$ which also has the desired property.  In these first constructions the middle thirds Cantor set $C$ is the control set and we exploit the fact that it has $\mathbb{Z}_m$-regular self covers for all $m\in \mathbb{Z}$ in order to produce 1-dimensional continua with the same property.  These continua are not manifolds, in fact, they are not even locally finite cell complexes, since the Cantor set is also a set of exceptional points in each.  They are, however, metric continua.

\bigskip
From this continua theoretic perspective, covering space maps can be characterized via the following: if $f:X\rightarrow Y$ is a proper local homeomorphism of a metric continuum $X$ onto the metric continuum $Y$, then it is a covering space map.  Recall that $f$ is \emph{proper}, if $f^{-1}(K)$ is compact, whenever $K\subset Y$ is compact.  See Jungck \cite{JungckLocalHomeo1983} for a more general discussion of properties of local homeomorphisms.  Bellamy's result \cite{Bellamy2005}, that the pseudo circle has non-trivial self covers of all orders, is a continua theoretic result of considerable interest.  His proof is similar to the Cantor set controlled constructions in the sense that he takes some care with how the pseudo circle sits in the plane.  He exploits properties of certain analytic maps $f:\mathbb{C}\rightarrow \mathbb{C}$ to complete the proof.

\begin{example}\label{ex:NecklaceCircles}\emph{A necklace of circles.}
Start with the interval $[0,2]$ sitting in the plane, say on the $x$-axis.  Build the middle-thirds Cantor set $C$ in $[0,1]\subset[0,2]$.  For each pair of rational points $\frac{i}{3^j}, \ \frac{i+1}{3^j}$ in $C$, attach a copy of $S^1$ of radius $\frac{1}{2\cdot 3^j}$, with center at the midpoint $\frac{2i+1}{2\cdot 3^j}$ of $[\frac{i}{3^j},\frac{i+1}{3^j}]$. Next, delete the interior of $[1,2]$ and attach a circle of radius $\frac{1}{2}$ and center $\frac{3}{2}$.  Finally, identify $0$ and $2$ to form a necklace $\mathcal{N}^{(1)}_{C,S^1}$ of circles.  The necklace $\mathcal{N}^{(1)}_{C,S^1}$ has $\mathbb{Z}_m$-regular self-covers $f_m$ for every $m \in \mathbb{N}$.  A picture of $\mathcal{N}^{(1)}_{C,S^1}$ and the $\mathbb{Z}_3$-regular cover are in Figure \ref{fig:necklace_of_circles}.  The interested reader can write out formulae for the $f_m$.

In the construction note that the circles attached to the rational points of the Cantor set $C$ have diameters which decrease to $0$ \emph{in levels}.
At the $0^{\textrm{th}}$ level there is one (topological) circle of diameter $1$, at the $1^{\textrm{st}}$ level there is one circle of diameter $\frac{1}{3}$, at the $2^{\textrm{nd}}$ level there are two circles of diameter $\frac{1}{9}$, at the $3^{\textrm{rd}}$ level there are four circles of diameter $\frac{1}{27}$ and, in general, at the $j^{\textrm{th}}$ level there are $2^{j-1}$ circles of diameter $\frac{1}{3^j}$.

\bigskip
This necklace $\mathcal{N}^{(1)}_{C,S^1}$ of circles can be modified in various ways to produce other continua with $\mathbb{Z}_m$-regular self covers for many $m \in \mathbb{Z}$ while continuing to use the Cantor set as a control set. The first, $\mathcal{N}^{(2)}_{C,S^1}$ is another necklace of circles with the diameters of the circles decreasing to 0 in levels.  It is illustrated in Figure \ref{fig:2nd_necklace_of_circles}.   The reader is invited to draw the picture for the $\mathbb{Z}_3$-regular cover of $\mathcal{N}^{(2)}_{C,S^1}$ or other $\mathbb{Z}_m$-regular covers $f_m:\mathcal{N}^{(2)}_{C,S^1}\rightarrow \mathcal{N}^{(2)}_{C,S^1}$.

The construction of the $\mathcal{N}^{(2)}_{C,S^1}$ can also be described as a replacement of each of the circles with diameters decreasing in levels to 0 by copies of the continuum $K=\{(x,y):x^2+y^2=1\}\cup [1,2]\times 0\cup [-2,-1]\times 0$, the union of $S^1$ and two intervals, which have diameters decreasing to 0 in levels.

\begin{figure}
\centering
\begin{tikzpicture}[scale=.8]
\node[fill=white] (2) at (12.4,-.75) {$2 = 0$}; \draw[thick] (12,-.35)--(12,.35);
\draw[fill=white] (9,0) ellipse (3cm and 2cm);
\draw[dashed] (12,0)  arc[start angle=0, end angle=180, x radius=3cm, y radius=.75cm];
\draw (6,0)  arc[start angle=180, end angle=360, x radius=3cm, y radius=.75cm];
\draw[dotted] (-6,0) -- (6,0);
\draw[thick] (6,-.35) -- (6,.35); \node (1) at (6,-.75) {$1$};
\draw[thick] (-6,-.35) -- (-6,.35); \node[fill=white] (0) at (-6.25,-.75) {$0=2$};
\draw[fill=white] (0,0) circle (2cm);
\draw[thick] (2,-.35) -- (2,.35); \node (2/3) at (2,-.75) {$\frac{2}{3}$};
\draw[thick] (-2,-.35) -- (-2,.35); \node(1/3) at (-2,-.75) {$\frac{1}{3}$};
\draw[dashed] (2,0)  arc[start angle=0, end angle=180, x radius=2cm, y radius=.5cm];
\draw (-2,0)  arc[start angle=180, end angle=360, x radius=2cm, y radius=.5cm];
\draw[fill=white] (4,0) circle (.667cm);
\draw[dashed] (4.667,0)  arc[start angle=0, end angle=180, x radius=.667cm, y radius=.16675cm];
\draw (3.333,0)  arc[start angle=180, end angle=360, x radius=.667cm, y radius=.16675cm];
\draw[thick] (4.6667,-.35) -- (4.6667,.35);
\node (8/9) at (4.6667,-.75) {$\frac{8}{9}$};
\draw[thick] (3.3333,-.35) -- (3.3333,.35); \node (7/9) at (3.3333,-.75) {$\frac{7}{9}$};
\draw[fill=white] (-4,0) circle (.667cm);
\draw[dashed] (-3.333,0)  arc[start angle=0, end angle=180, x radius=.667cm, y radius=.16675cm];
\draw (-4.667,0)  arc[start angle=180, end angle=360, x radius=.667cm, y radius=.16675cm];
\draw[thick] (-4.6667,-.35) -- (-4.6667,.35); \node (1/9) at (-4.6667,-.75) {$\frac{1}{9}$};
\draw[thick] (-3.3333,-.35) -- (-3.3333,.35); \node (2/9) at (-3.3333,-.75) {$\frac{2}{9}$};
\draw[fill=white] (2.66665,0) circle (.2222);
\draw[dashed] (2.88885,0)  arc[start angle=0, end angle=180, x radius=.2222cm, y radius=.05555cm];
\draw (2.44445,0)  arc[start angle=180, end angle=360, x radius=.2222cm, y radius=.05555cm];
\draw[thick] (2.44443,-1) -- (2.44443,.35); \node (19/27) at (2.44443,-1.4) {$\frac{19}{27}$};
\draw[thick] (2.88887,-.35) -- (2.88887,.35); \node (20/27) at (2.88887,.75) {$\frac{20}{27}$};
\draw[fill=white] (5.33335,0) circle (.2222);
\draw[dashed] (5.55555,0)  arc[start angle=0, end angle=180, x radius=.2222cm, y radius=.05555cm];
\draw (5.11115,0)  arc[start angle=180, end angle=360, x radius=.2222cm, y radius=.05555cm];
\draw[thick] (5.11113,-.35) -- (5.11113,.35); \node (25/27) at (5.11113,.75) {$\frac{25}{27}$};
\draw[thick] (5.55557,-1) -- (5.55557,.35); \node (26/27) at (5.55557,-1.4) {$\frac{26}{27}$};
\draw[fill=white] (-2.66665,0) circle (.2222);
\draw[thick] (-2.44443,-1) -- (-2.44443,.35); \node (8/27) at (-2.44443,-1.4) {$\frac{8}{27}$};
\draw[dashed] (-2.44445,0)  arc[start angle=0, end angle=180, x radius=.2222cm, y radius=.05555cm];
\draw (-2.88885,0)  arc[start angle=180, end angle=360, x radius=.2222cm, y radius=.05555cm];
\draw[thick] (-2.88887,-.35) -- (-2.88887,.35); \node (7/27) at (-2.88887,.75) {$\frac{7}{27}$};
\draw[fill=white] (-5.33335,0) circle (.2222);
\draw[dashed] (-5.11115,0)  arc[start angle=0, end angle=180, x radius=.2222cm, y radius=.05555cm];
\draw (-5.55555,0)  arc[start angle=180, end angle=360, x radius=.2222cm, y radius=.05555cm];
\draw[thick] (-5.11113,-.35) -- (-5.11113,.35); \node (2/27) at (-5.11113,.75) {$\frac{2}{27}$};
\draw[thick] (-5.55557,-1) -- (-5.55557,.35); \node (1/27) at (-5.55557,-1.4) {$\frac{1}{27}$};
\end{tikzpicture}
\caption{Cantor's Pearl Necklace $\mathcal{N}^{(2)}_{C,S^2}$.  The vertical line segments in the picture corresponding to the rational points of the Cantor set are included only for reference and are not part of the Pearl Necklace.  However, should the reader wish to include them, this does produces another necklace with $\mathbb{Z}_m$-regular self covers for every $m \in \mathbb{Z}$.} \label{fig:CantorsPearlNecklace}
\end{figure}
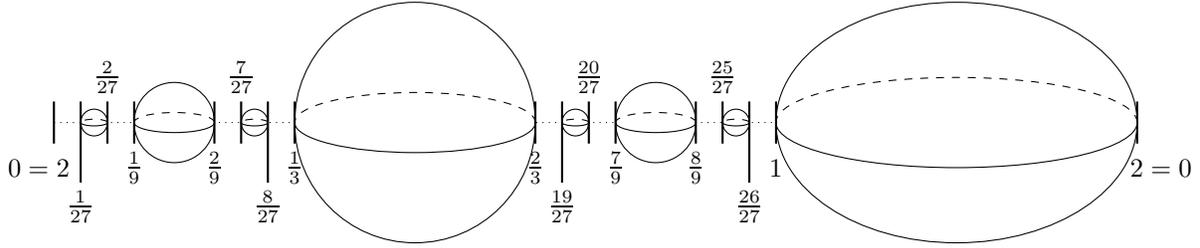

\bigskip
In fact, this can be done in general.  Fix a triple $(K,x_l,x_r)$ where $K$ is a metric continuum and $x_l \neq x_r$ are distinct points in $K$ which determine its diameter. Let $\mathcal{K}$ be a countably infinite collection of copies of $(K,x_l,x_r)$ with diameters which decrease to 0 in levels.  Specifically, in $\mathcal{K}$, the $0^{\textrm{th}}$ level contains one copy of $K$ of diameter $1$, the $1^{\textrm{st}}$ level contains one copy of $K$ of diameter $\frac{1}{3}$ and, in general, for $j \geq 1$, the $j^{\textrm{th}}$ level contains $2^{j-1}$ copies of $K$ of diameter $\frac{1}{3^{j}}$. Attach each copy of $K$ in $\mathcal{K}$ of the appropriate diameter to the rational points $\frac{i}{3^p}$ and $\frac{i+1}{3^p}$ by identifying $x_l$ with $\frac{i}{3^p}$ and $x_r$ with $\frac{i+1}{3^p}$.  Call the resulting necklace $\mathcal{N}^{(2)}_{C,K}$.  Then for each $K$, there are $\mathbb{Z}_m$-regular self covers $f_m:\mathcal{N}^{(2)}_{C,K}\rightarrow \mathcal{N}^{(2)}_{C,K}$.  Note that when $K$ is non-planar or high dimensional, it follows that $\mathcal{N}^{(2)}_{C,K}$ is also, respectively, non-planar or high dimensional.

\begin{example}\emph{Cantor's Pearl Necklace.}  Of particular interest is the necklace $\mathcal{N}^{(2)}_{C,S^2}$ where $\mathcal{K}$ is a collection of countably infinitely many copies of $(S^2,x_l,x_r)$, with diameters decreasing to 0 in levels where $(x_l,x_r)$ is a pair of antipodal points in the 2-sphere $S^2$.  See Figure \ref{fig:CantorsPearlNecklace} and \cite{BedDelTimm2006} Example 5.1 and Figure 9.  This necklace clearly has uncountably generated $\pi_2$.

\begin{question}
What is the structure of  $\pi_2(\mathcal{N}^{(2)}_{C,S^2})$?
\end{question}
\end{example}

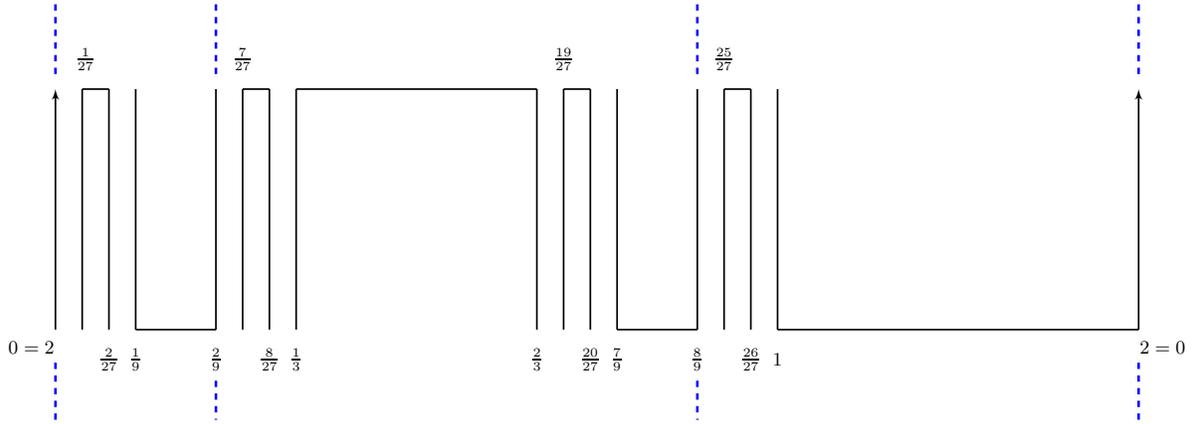
\begin{figure}[t]
\centering
\scalebox{.8}{\begin{tikzpicture}\tikzstyle{every node}=[font=\small][scale=1]
\node[fill=white] (2) at (12.4,-2.3) {$2 = 0$}; \draw[thick,->, > = latex'](12,-2)--(12,2);
\draw[thick] (6,-2)-- (12,-2);
\draw[thick] (6,-2) -- (6,2); \node (1) at (6,-2.5) {$1$};
\draw[thick,->, > = latex'](-6,-2) -- (-6,2); \node[fill=white] (0) at (-6.4,-2.3) {$0=2$};
\draw[thick] (2,-2) -- (2,2); \node (2/3) at (2,-2.5) {$\frac{2}{3}$};
\draw[thick] (-2,-2) -- (-2,2); \node(1/3) at (-2,-2.5) {$\frac{1}{3}$};
\draw[thick] (-2,2)--(2,2);
\draw[thick] (4.6667,-2) -- (4.6667,2);
\node (8/9) at (4.6667,-2.5) {$\frac{8}{9}$};
\draw[thick] (3.3333,-2) -- (3.3333,2); \node (7/9) at (3.3333,-2.5) {$\frac{7}{9}$};
\draw[thick] (3.3333,-2)--(4.6667,-2);
\draw[thick] (-4.6667,-2) -- (-4.6667,2); \node (1/9) at (-4.6667,-2.5) {$\frac{1}{9}$};
\draw[thick] (-3.3333,-2) -- (-3.3333,2); \node (2/9) at (-3.3333,-2.5) {$\frac{2}{9}$};
\draw[thick] (-4.6667,-2)--(-3.3333,-2);
\draw[thick] (2.44443,-2) -- (2.44443,2); \node (19/27) at (2.44443,2.5) {$\frac{19}{27}$};
\draw[thick] (2.88887,-2) -- (2.88887,2); \node (20/27) at (2.88887,-2.5) {$\frac{20}{27}$};
\draw[thick] (2.44443,2)--(2.88887,2);
\draw[thick] (5.11113,-2) -- (5.11113,2); \node (25/27) at (5.11113,2.5) {$\frac{25}{27}$};
\draw[thick] (5.55557,-2) -- (5.55557,2); \node (26/27) at (5.55557,-2.5) {$\frac{26}{27}$};
\draw[thick] (5.11113,2)--(5.55557,2);
\draw[thick] (-2.44443,-2) -- (-2.44443,2); \node (8/27) at (-2.44443,-2.5) {$\frac{8}{27}$};
\draw[thick] (-2.88887,-2) -- (-2.88887,2); \node (7/27) at (-2.88887,2.5) {$\frac{7}{27}$};
\draw[thick] (-2.88887,2)--(-2.44443,2);
\draw[thick] (-5.11113,-2) -- (-5.11113,2); \node (2/27) at (-5.11113,-2.5) {$\frac{2}{27}$};
\draw[thick] (-5.55557,-2) -- (-5.55557,2); \node (1/27) at (-5.5,2.5) {$\frac{1}{27}$};
\draw[thick] (-5.11113,2)--(-5.55557,2);
\draw[very thick, dashed, color=blue] (-6,2.25)--(-6,3.5); \draw[very thick, dashed, color=blue] (-6,-2.55)--(-6,-3.5);
\draw[very thick, dashed, color=blue] (-3.3333,-2.85)--(-3.3333,-3.5);\draw[very thick, dashed, color=blue] (-3.3333,2.25)--(-3.3333,3.5);
\draw[very thick, dashed, color=blue] (4.6667,2.25)--(4.6667,3.5); \draw[very thick, dashed, color=blue] (4.6667,-2.85)--(4.6667,-3.5);
\draw[very thick, dashed, color=blue] (12,2.25)--(12,3.5); \draw[very thick, dashed, color=blue] (12,-2.55)--(12,-3.5);
\end{tikzpicture}}
\caption{The circle-like Knaster Cup-cap continuum. The dashed blue in the picture indicate the closures of fundamental domains of a $\mathbb{Z}_3$-regular self cover.  The arrows on the segments $0\times[-1,1]$ and $2\times [-1,1]$ show the gluing map $(0,t)\mapsto (2,t)$ used to form the circle-like Cup-cap continuum.} \label{fig:CircleCup-cap}
\end{figure}

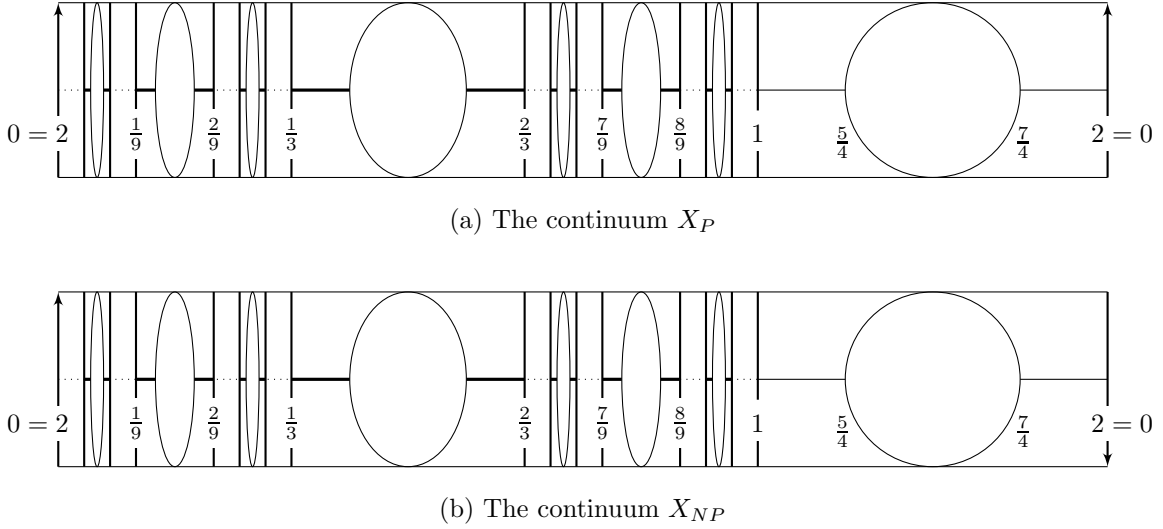
\begin{figure}[t]
\centering
\begin{tikzpicture}[scale=.775]
\draw[thick,->, > = latex'] (12,-1.5)--(12,1.5);  \node[fill=white] (2) at (12.25,-.75) {$2=0$};
\draw[dotted] (-6,0) -- (6,0);
\draw[thick] (6,-.35) -- (6,1.5); \draw[thick] (6,-1.5)--(6,-.75); \node[fill=white] (1) at (6,-.75) {$1$};
\draw[thick,->, > = latex'] (-6,-1.5) -- (-6,1.5); \node[fill=white] (0) at (-6.35,-.75) {$0=2$};
\draw (6,0)--(12,0); \draw[fill=white] (9,0) circle (1.5cm);
\node (5/4) at (7.45,-.9) {$\frac{5}{4}$}; \node (7/4) at (10.55,-.9) {$\frac{7}{4}$};
\draw[fill=white] (0,0) ellipse (1cm and 1.5cm);
\draw[very thick] (1,0) -- (2,0); \draw[thick] (2,-1.5) -- (2,1.5); \node[fill=white] (2/3) at (2,-.75) {$\frac{2}{3}$};
\draw[very thick] (-1,0) -- (-2,0); \draw[thick] (-2,-1.5) -- (-2,1.5); \node[fill=white] (1/3) at (-2,-.75) {$\frac{1}{3}$};
\draw[fill=white] (4,0) ellipse (.333cm and 1.5cm);
\draw[very thick] (4.333,0) -- (4.6667,0); \draw[thick] (4.6667,-1.5) -- (4.6667,1.5); \node[fill=white] (8/9) at (4.6667,-.75) {$\frac{8}{9}$};
\draw[very thick] (3.6667,0) -- (3.3333,0); \draw[thick] (3.3333,-1.5) -- (3.3333,1.5); \node[fill=white] (7/9) at (3.3333,-.75) {$\frac{7}{9}$};
\draw[fill=white] (-4,0) ellipse (.333cm and 1.5cm);
\draw[very thick] (-4.333,0) -- (-4.6667,0); \draw[thick] (-4.6667,-1.5) -- (-4.6667,1.5); \node[fill=white] (1/9) at (-4.6667,-.75) {$\frac{1}{9}$};
\draw[very thick] (-3.6667,0) -- (-3.3333,0); \draw[thick] (-3.3333,-1.5) -- (-3.3333,1.5); \node[fill=white] (2/9) at (-3.3333,-.75) {$\frac{2}{9}$};
\draw[fill=white] (2.66665,0) ellipse (.1111cm and 1.5cm);
\draw[very thick] (2.44443,0) -- (2.55554,0); \draw[thick] (2.44443,-1.5) -- (2.44443,1.5); 
\draw[very thick] (2.77776,0) -- (2.88887,0); \draw[thick] (2.88887,-1.5) -- (2.88887,1.5); 
\draw[fill=white] (5.33335,0) ellipse (.1111cm and 1.5cm);
\draw[very thick] (5.11113,0) -- (5.22224,0); \draw[thick] (5.11113,-1.5) -- (5.11113,1.5); 
\draw[very thick] (5.44446,0) -- (5.55557,0); \draw[thick] (5.55557,-1.5) -- (5.55557,1.5); 
\draw[fill=white] (-2.66665,0) ellipse (.1111cm and 1.5cm);
\draw[very thick] (-2.44443,0) -- (-2.55554,0); \draw[thick] (-2.44443,-1.5) -- (-2.44443,1.5); 
\draw[very thick] (-2.77776,0) -- (-2.88887,0); \draw[thick] (-2.88887,-1.5) -- (-2.88887,1.5); 
\draw[fill=white] (-5.3335,0) ellipse (.1111cm and 1.5cm);
\draw[very thick] (-5.11113,0) -- (-5.22224,0); \draw[thick] (-5.11113,-1.5) -- (-5.11113,1.5); 
\draw[very thick] (-5.44446,0) -- (-5.55557,0); \draw[thick] (-5.55557,-1.5) -- (-5.55557,1.5); 
\draw (-6,1.5) -- (12,1.5); \draw (-6,-1.5) -- (12,-1.5);
\node (Xp) at (3,-2.25) {\textrm{(a) The continuum }$X_P$};
\node (spacer) at (3,-3.25) {$ \ $};
\end{tikzpicture}

\begin{tikzpicture}[scale=.775]
\draw[thick,->, > = latex'] (12,1.5)--(12,-1.5);  \node[fill=white] (2) at (12.25,-.75) {$2=0$};
\draw[dotted] (-6,0) -- (6,0);
\draw[thick] (6,-.35) -- (6,1.5); \draw[thick] (6,-1.5)--(6,-.75); \node[fill=white] (1) at (6,-.75) {$1$};
\draw[thick,->, > = latex'] (-6,-1.5) -- (-6,1.5); \node[fill=white] (0) at (-6.35,-.75) {$0=2$};
\draw (6,0)--(12,0); \draw[fill=white] (9,0) circle (1.5cm);
\node (5/4) at (7.45,-.9) {$\frac{5}{4}$}; \node (7/4) at (10.55,-.9) {$\frac{7}{4}$};
\draw[fill=white] (0,0) ellipse (1cm and 1.5cm);
\draw[very thick] (1,0) -- (2,0); \draw[thick] (2,-1.5) -- (2,1.5); \node[fill=white] (2/3) at (2,-.75) {$\frac{2}{3}$};
\draw[very thick] (-1,0) -- (-2,0); \draw[thick] (-2,-1.5) -- (-2,1.5); \node[fill=white] (1/3) at (-2,-.75) {$\frac{1}{3}$};
\draw[fill=white] (4,0) ellipse (.333cm and 1.5cm);
\draw[very thick] (4.333,0) -- (4.6667,0); \draw[thick] (4.6667,-1.5) -- (4.6667,1.5); \node[fill=white] (8/9) at (4.6667,-.75) {$\frac{8}{9}$};
\draw[very thick] (3.6667,0) -- (3.3333,0); \draw[thick] (3.3333,-1.5) -- (3.3333,1.5); \node[fill=white] (7/9) at (3.3333,-.75) {$\frac{7}{9}$};
\draw[fill=white] (-4,0) ellipse (.333cm and 1.5cm);
\draw[very thick] (-4.333,0) -- (-4.6667,0); \draw[thick] (-4.6667,-1.5) -- (-4.6667,1.5); \node[fill=white] (1/9) at (-4.6667,-.75) {$\frac{1}{9}$};
\draw[very thick] (-3.6667,0) -- (-3.3333,0); \draw[thick] (-3.3333,-1.5) -- (-3.3333,1.5); \node[fill=white] (2/9) at (-3.3333,-.75) {$\frac{2}{9}$};
\draw[fill=white] (2.66665,0) ellipse (.1111cm and 1.5cm);
\draw[very thick] (2.44443,0) -- (2.55554,0); \draw[thick] (2.44443,-1.5) -- (2.44443,1.5); 
\draw[very thick] (2.77776,0) -- (2.88887,0); \draw[thick] (2.88887,-1.5) -- (2.88887,1.5); 
\draw[fill=white] (5.33335,0) ellipse (.1111cm and 1.5cm);
\draw[very thick] (5.11113,0) -- (5.22224,0); \draw[thick] (5.11113,-1.5) -- (5.11113,1.5); 
\draw[very thick] (5.44446,0) -- (5.55557,0); \draw[thick] (5.55557,-1.5) -- (5.55557,1.5); 
\draw[fill=white] (-2.66665,0) ellipse (.1111cm and 1.5cm);
\draw[very thick] (-2.44443,0) -- (-2.55554,0); \draw[thick] (-2.44443,-1.5) -- (-2.44443,1.5); 
\draw[very thick] (-2.77776,0) -- (-2.88887,0); \draw[thick] (-2.88887,-1.5) -- (-2.88887,1.5); 
\draw[fill=white] (-5.3335,0) ellipse (.1111cm and 1.5cm);
\draw[very thick] (-5.11113,0) -- (-5.22224,0); \draw[thick] (-5.11113,-1.5) -- (-5.11113,1.5); 
\draw[very thick] (-5.44446,0) -- (-5.55557,0); \draw[thick] (-5.55557,-1.5) -- (-5.55557,1.5); 
\draw (-6,1.5) -- (12,1.5); \draw (-6,-1.5) -- (12,-1.5);
\node (Xp) at (3,-2.25) {\textrm{(b) The continuum }$X_{NP}$};
\end{tikzpicture}
\caption{The two continua $X_P$ with $\mathbb{Z}_m$-regular self covers for all $m$ and $X_{NP}$ with $\mathbb{Z}_{2m}$-regular self covers for all $m$ and control set the product $C\times [-1,1]$. The tops and bottoms of the rectangles in both pictures are parts of the continua. The direction of the gluing map on $[-1,1]$ at $0=2$ and $2=0$ are indicated by the arrows in the picture.}
\label{fig:Xp_and_Xnp}
\end{figure}

\noindent  Note that one could also produce another Cantor's Necklaces by taking $\mathcal{K}$ to be a collection of copies of $n$-spheres $(S^n,x_l,x_r)$, with diameters which decrease to 0 in levels, $(x_l,x_r)$ antipodal points in $S^n$, or a collection of copies of closed $n$-balls $(B^n,x_l,x_r)$, with diameters which decrease to 0 in levels, $(x_l,x_r)$ antipodal points in $\partial B^n$.  Each of these spaces can be visually represented via pictures similar to those in Figures \ref{fig:necklace_of_circles} and \ref{fig:CantorsPearlNecklace} and all have $\mathbb{Z}_m$-regular self covers for every $m \in \mathbb{Z}$.

\begin{figure}[t]
\centering
\begin{tikzpicture}[scale=.7]
\draw[very thick] (-6,-3) rectangle (12,3);
\node[fill=white] (1onY) at (-6.15,3) {1}; \node[fill=white] (-1onY) at (-6.25,-3) {-1};
\node[fill=white] (2) at (12.4,-.75) {$2 = 0$};
\draw[very thick,->, > = latex'] (-6,.35) -> (-6,2.9);
\draw[very thick,->, > = latex'] (12,.35) -> (12,2.9);
\draw[fill=blue!30] (9,0) ellipse (3cm and 2cm);
\draw[dotted] (-6,0) -- (6,0);
\draw[thick] (6,-.35) -- (6,.35); \node (1) at (6,-.75) {$1$};
\draw[thick] (-6,-.35) -- (-6,.35); \node[fill=white] (0) at (-6.25,-.75) {$0=2$};
\draw[fill=blue!30] (0,0) circle (2cm);
\draw[thick] (2,-.35) -- (2,3); \node (2/3) at (2,-.75) {$\frac{2}{3}$}; \draw[thick] (2,-3) -- (2,-1.15);
\draw[thick] (-2,-.35) -- (-2,.35); \node(1/3) at (-2,-.75) {$\frac{1}{3}$};
\draw[fill=blue!30] (4,0) circle (.667cm);
\draw[thick] (4.6667,-.35) -- (4.6667,.35);
\node (8/9) at (4.6667,-.75) {$\frac{8}{9}$};
\draw[thick] (3.3333,-.35) -- (3.3333,.35); \node (7/9) at (3.3333,-.75) {$\frac{7}{9}$};
\draw[fill=blue!30] (-4,0) circle (.667cm);
\draw[thick] (-4.6667,-.35) -- (-4.6667,.35); \node (1/9) at (-4.6667,-.75) {$\frac{1}{9}$};
\draw[thick] (-3.3333,-.35) -- (-3.3333,3); \node (2/9) at (-3.3333,-.75) {$\frac{2}{9}$}; \draw[thick] (-3.3333,-3) -- (-3.3333,-1.15);
\draw[fill=blue!30] (2.66665,0) circle (.2222);
\draw[thick] (2.44443,-1) -- (2.44443,.35); \node (19/27) at (2.44443,-1.4) {$\frac{19}{27}$};
\draw[thick] (2.88887,-.35) -- (2.88887,.35); \node (20/27) at (2.88887,.75) {$\frac{20}{27}$};
\draw[fill=blue!30] (5.33335,0) circle (.2222);
\draw[thick] (5.11113,-.35) -- (5.11113,.35); \node (25/27) at (5.11113,.75) {$\frac{25}{27}$};
\draw[thick] (5.55557,-1) -- (5.55557,.35); \node (26/27) at (5.55557,-1.4) {$\frac{26}{27}$};
\draw[fill=blue!30] (-2.66665,0) circle (.2222);
\draw[thick] (-2.44443,-1) -- (-2.44443,.35); \node (8/27) at (-2.44443,-1.4) {$\frac{8}{27}$};
\draw[thick] (-2.88887,-.35) -- (-2.88887,.35); \node (7/27) at (-2.88887,.75) {$\frac{7}{27}$};
\draw[fill=blue!30] (-5.33335,0) circle (.2222);
\draw[thick] (-5.11113,-.35) -- (-5.11113,.35); \node (2/27) at (-5.11113,.75) {$\frac{2}{27}$};
\draw[thick] (-5.55557,-1) -- (-5.55557,.35); \node (1/27) at (-5.55557,-1.4) {$\frac{1}{27}$};
\tikzset{edge/.style = {->,> = latex'}}
\draw[edge,thick] (3,-3.5) to (3,-5);
\node (x3) at (3.5,-4.25) {$\times 3$};
\end{tikzpicture}

\begin{tikzpicture}[scale=.7]
\draw[very thick] (-6,-3) rectangle (12,3);
\node[fill=white] (1onY) at (-6.15,3) {1}; \node[fill=white] (-1onY) at (-6.25,-3) {-1};
\node[fill=white] (2) at (12.4,-.75) {$2 = 0$};
\draw[very thick,->, > = latex'] (-6,.35) -> (-6,2.9);
\draw[very thick,->, > = latex'] (12,.35) -> (12,2.9);
\draw[fill=blue!30] (9,0) ellipse (3cm and 2cm);
\draw[dotted] (-6,0) -- (6,0);
\draw[thick] (6,-.35) -- (6,.35); \node (1) at (6,-.75) {$1$};
\draw[thick] (-6,-.35) -- (-6,.35); \node[fill=white] (0) at (-6.25,-.75) {$0=2$};
\draw[fill=blue!30] (0,0) circle (2cm);
\draw[thick] (2,-.35) -- (2,.35); \node (2/3) at (2,-.75) {$\frac{2}{3}$};
\draw[thick] (-2,-.35) -- (-2,.35); \node(1/3) at (-2,-.75) {$\frac{1}{3}$};
\draw[fill=blue!30] (4,0) circle (.667cm);
\draw[thick] (4.6667,-.35) -- (4.6667,.35);
\node (8/9) at (4.6667,-.75) {$\frac{8}{9}$};
\draw[thick] (3.3333,-.35) -- (3.3333,.35); \node (7/9) at (3.3333,-.75) {$\frac{7}{9}$};
\draw[fill=blue!30] (-4,0) circle (.667cm);
\draw[thick] (-4.6667,-.35) -- (-4.6667,.35); \node (1/9) at (-4.6667,-.75) {$\frac{1}{9}$};
\draw[thick] (-3.3333,-.35) -- (-3.3333,.35); \node (2/9) at (-3.3333,-.75) {$\frac{2}{9}$};
\draw[fill=blue!30] (2.66665,0) circle (.2222);
\draw[thick] (2.44443,-1) -- (2.44443,.35); \node (19/27) at (2.44443,-1.4) {$\frac{19}{27}$};
\draw[thick] (2.88887,-.35) -- (2.88887,.35); \node (20/27) at (2.88887,.75) {$\frac{20}{27}$};
\draw[fill=blue!30] (5.33335,0) circle (.2222);
\draw[thick] (5.11113,-.35) -- (5.11113,.35); \node (25/27) at (5.11113,.75) {$\frac{25}{27}$};
\draw[thick] (5.55557,-1) -- (5.55557,.35); \node (26/27) at (5.55557,-1.4) {$\frac{26}{27}$};
\draw[fill=blue!30] (-2.66665,0) circle (.2222);
\draw[thick] (-2.44443,-1) -- (-2.44443,.35); \node (8/27) at (-2.44443,-1.4) {$\frac{8}{27}$};
\draw[thick] (-2.88887,-.35) -- (-2.88887,.35); \node (7/27) at (-2.88887,.75) {$\frac{7}{27}$};
\draw[fill=blue!30] (-5.33335,0) circle (.2222);
\draw[thick] (-5.11113,-.35) -- (-5.11113,.35); \node (2/27) at (-5.11113,.75) {$\frac{2}{27}$};
\draw[thick] (-5.55557,-1) -- (-5.55557,.35); \node (1/27) at (-5.55557,-1.4) {$\frac{1}{27}$};
\end{tikzpicture}
\caption{The $\mathbb{Z}_3$-regular cover $f_3:A^{(1)}_C\rightarrow A^{(1)}_C$ of Cantor's Annulus $A^{(1)}_C$ with countably infinitely many holes and with the boundaries of the holes tangent to the Cantor set at rational points of the Cantor set. The lifts $f^{-1}_3(0=2 \times [-1,1])$ of the line segment $0=2 \times [-1,1]$ in the range are shown in the domain copy of $A^{(1)}_C$.  They delineate closures of fundamental domains of the covering map. Note that $f_3$ is an expansive map on the lefthand and middle fundamental domains and on the portion of the righthand fundamental domain between $\frac{2}{3}$ and $1$.  It is the identity on the portion of righthand fundamental domain between $1$ and $2=0$.}\label{fig:CountPuncAnnTangent}
\end{figure}
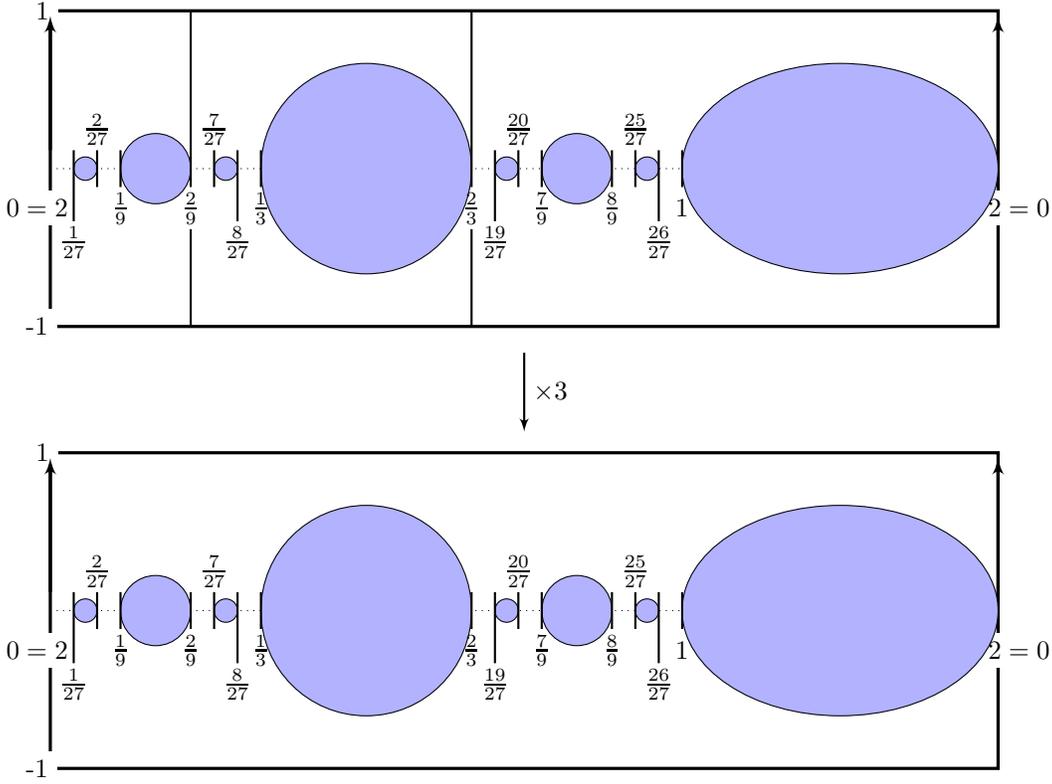

Another example of interest occurs by taking $\mathcal{K}$ to be a countably infinite collections of copies of a knotted arc $\alpha$.  Let $x_l$ be one endpoint of $\alpha$, $x_r$ is the other, and assume that the diameters of the copies of $(\alpha,x_l,x_r)$ decrease to 0 in levels.  The resulting continuum, $\mathcal{N}^{(2)}_{C,\alpha}$ is a continuum with $\mathbb{Z}_m$-regular self covers for all $m\in \mathbb{N}$.  In fact, doing this construction in $S^3$ for different $\alpha$ produces continua  $\mathcal{N}^{(2)}_{C,\alpha}$, which are all different wild knotted $S^1$'s in $S^3$ with a Cantor's set of wild points.  Provided they sit in $S^3$ in the correct way, see \cite{Nanyes_Timm_2004}, they are periodic in the sense made precise by Erica Flapan \cite{Flapan_1985}.  In this work she shows that tame knots can have at most finitely many periods.  By \cite{Nanyes_Timm_2004}, the wild $\mathcal{N}^{(2)}_{C,\alpha}$'s all have infinitely many periods, showing that the tameness hypothesis in \cite{Flapan_1985} is necessary.

\bigskip
A third class of examples of continua $X$ with $\mathbb{Z}_m$-regular self covers for all $m\in \mathbb{N}$ can be obtained by attaching copies of a planar continuum $K$ to the Cantor set whose ``widths'' decrease to 0 in levels, but for which the diameters do not.  We illustrate with a particular example.  To develop a more general process requires addressing complications arising from the copies of $K$ sitting between the pairs of rational points in $C$.  One must attach a copy of a continuum $Y$ at each of the rational points of the Cantor set in such a way that the continua attached between the pairs of rational points will converge to the copies of $Y$ in the appropriate way.

Consider the interval $[0,2] \times 0$ contained in the rectangle $[0,2]\times [-1,1] \subset \mathbb{R}^2$.  Choose the middle thirds Cantor set $C \times 0 \subset [0,1]\times 0 \subset [0,2]\times 0$, and the Cantor ladder $C\times [-1,1]$ sitting in the rectangle.  Note that the Cantor ladder contains $\{(r,0)\in C \times 0: r\textrm{ is a rational point in }C \}\times 0\cup C\times [-1,1] \cup 2\times [-1,1]$.  Glue the left side $0\times [-1,1]$ of this continuum to its right side via $(0,y)\sim (2,y)$ to form a continuum $X$.  Now, take a collection $\mathcal{K}$ of copies of $K = \{(x,y):x^2+y^2=1\} \cup [-2,-1]\times 0 \cup [1,2]\times 0$ with ``widths'' which decrease to 0 in levels, but whose ``height'' does not.  Between each pair of $\frac{i}{3^j},\frac{i+1}{3^j}$ of rational points in $C\times 0$, delete the open interval $(\frac{i}{3^j},\frac{i+1}{3^j})\times 0$, and attach a copy of $K$ of height 2 and of the appropriate width.  Also delete $(1,2) \times 0$ and attach a copy of $K$ of width 1 and height equal to 2.  The control set in this case is $C \times [-1,1]$.   A picture of the resulting necklace $\mathcal{N}^{(3)}_{C\times I,S^1}$ is in Figure \ref{fig:3rd_necklace_of_circles}.  Note its similarity to the continuum in \ref{fig:Xp_and_Xnp}.  The control set for this continuum is also $C \times [-1,1]$.

\bigskip
The necklace $\mathcal{N}^{(3)}_{C \times I,S^1}$ can be further modified to produce additional 1-dimensional continua which have $\mathbb{Z}_m$-regular self covers for infinitely many $m\in \mathbb{N}$.  We do so by using the boundary rectangle $[0,2]\times (\{-1,1\})) \cup (\{0,2\}\times [-1,1])$ of $[0,2]\times I$ as a frame.  We illustrate two such.  The first is the continuum $X_P=\mathcal{N}^{(3)}_{C\times I,S^1}\cup \left(\{-1,1\}\times [0,2]\right)$, which clearly embeds in an annulus and is, therefore, planar.  It has $\mathbb{Z}_m$-regular self covers for all $m \in \mathbb{N}$.  For the second, start with $X_P$ but change the gluing map $0\times [-1,1] \overset{0\times id}{\longrightarrow} 2\times[-1,1]$ from the identity on the second factor to the orientation reversing homeomorphism $(0,y) \mapsto (2, -y)$ on the second factor.  Call the resulting space $X_{NP}$.  Pictures of $X_P$ and $X_{NP}$ are in Figure \ref{fig:Xp_and_Xnp}.  It is easy to see that $X_{NP}$, while 1-dimensional, is non-planar since it contains copies of the complete bipartite graph $K_{3,3}$.  Inspecting its construction, it is also clear that it is a subcontinuum of the Mobius band.  Like the Mobius band, $X_{NP}$ has $\mathbb{Z}_{2m+1}$-regular self covers of all odd degrees $2m+1$, $m\in \mathbb{N}\cup 0$.  As is the case for the annulus and the Mobious band, there are also even order $\mathbb{Z}_{2m}$-regular covers $f_{2m}:X_P\rightarrow X_{NP}$ for all $m\in \mathbb{N}$.

Finally, note that each of these continua can be thought of as being obtained from $C\times I \cup \left([0,2]\times \{-1,1\}\right)$ by carefully attaching copies of the continuum $([-2,-1]\times 0)\cup \{(x,y):x^2+y^2=1\}\cup([1,2]\times 0$ to it.  It is clear that one can apply this sort of process for any planar continuum $K$, provided one takes care with the attaching processes and provides sufficient control over the ``widths'' of the copies of $K$ which are being attached.  For each such $K$ then, these processes produce a planar continuum $X_{P,K}$ with $\mathbb{Z}_m$-regular self covers and a non-planar continuum $X_{NP,K}$ with $\mathbb{Z}_{2m-1}$-regular self covers for all $m \in \mathbb{N}$.  It is also the case that there are even order covers $f_{2m}:X_{P,K}\rightarrow X_{NP,K}$ for all $m\in \mathbb{N}$.

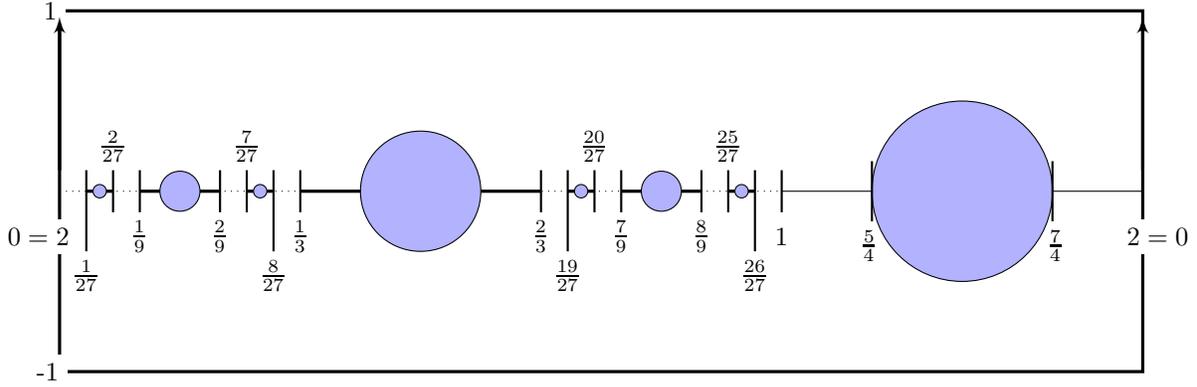
\begin{figure}[t]
\centering
\begin{tikzpicture}[scale=.8]
\draw[very thick] (-6,-3) rectangle (12,3);
\node[fill=white] (1onY) at (-6.15,3) {1}; \node[fill=white] (-1onY) at (-6.2,-3) {-1}; \node[fill=white] (2) at (12.25,-.75) {$2=0$};
\draw[very thick,->, > = latex'] (-6,.35) -> (-6,2.9);
\draw[very thick,->, > = latex'] (12,.35) -> (12,2.9);
\draw[dotted] (-6,0) -- (6,0);
\draw[thick] (6,-.35) -- (6,.35); \node (1) at (6,-.75) {$1$};
\draw[thick] (-6,-.35) -- (-6,.35); \node[fill=white] (0) at (-6.35,-.75) {$0=2$};
\draw (6,0)--(12,0); \draw[fill=blue!30] (9,0) circle (1.5cm);
\draw[thick] (7.5,-.5)--(7.5,.5);\draw[thick] (10.5,-.5)--(10.5,.5);
\node (5/4) at (7.45,-.9) {$\frac{5}{4}$}; \node (7/4) at (10.55,-.9) {$\frac{7}{4}$};
\draw[fill=blue!30] (0,0) circle (1cm);
\draw[very thick] (1,0) -- (2,0); \draw[thick] (2,-.35) -- (2,.35); \node (2/3) at (2,-.75) {$\frac{2}{3}$};
\draw[very thick] (-1,0) -- (-2,0); \draw[thick] (-2,-.35) -- (-2,.35); \node(1/3) at (-2,-.75) {$\frac{1}{3}$};
\draw[fill=blue!30] (4,0) circle (.333cm);
\draw[very thick] (4.333,0) -- (4.6667,0); \draw[thick] (4.6667,-.35) -- (4.6667,.35); \node (8/9) at (4.6667,-.75) {$\frac{8}{9}$};
\draw[very thick] (3.6667,0) -- (3.3333,0); \draw[thick] (3.3333,-.35) -- (3.3333,.35); \node (7/9) at (3.3333,-.75) {$\frac{7}{9}$};
\draw[fill=blue!30] (-4,0) circle (.333cm);
\draw[very thick] (-4.333,0) -- (-4.6667,0); \draw[thick] (-4.6667,-.35) -- (-4.6667,.35); \node (1/9) at (-4.6667,-.75) {$\frac{1}{9}$};
\draw[very thick] (-3.6667,0) -- (-3.3333,0); \draw[thick] (-3.3333,-.35) -- (-3.3333,.35); \node (2/9) at (-3.3333,-.75) {$\frac{2}{9}$};
\draw[fill=blue!30] (2.66665,0) circle (.1111); \draw[very thick] (2.44443,0) -- (2.55554,0); \draw[thick] (2.44443,-1) -- (2.44443,.35); \node (19/27) at (2.44443,-1.4) {$\frac{19}{27}$};
\draw[very thick] (2.77776,0) -- (2.88887,0); \draw[thick] (2.88887,-.35) -- (2.88887,.35); \node (20/27) at (2.88887,.75) {$\frac{20}{27}$};
\draw[fill=blue!30] (5.33335,0) circle (.1111); \draw[very thick] (5.11113,0) -- (5.22224,0); \draw[thick] (5.11113,-.35) -- (5.11113,.35); \node (25/27) at (5.11113,.75) {$\frac{25}{27}$};
\draw[very thick] (5.44446,0) -- (5.55557,0); \draw[thick] (5.55557,-1) -- (5.55557,.35); \node (26/27) at (5.55557,-1.4) {$\frac{26}{27}$};
\draw[fill=blue!30] (-2.66665,0) circle (.1111); \draw[very thick] (-2.44443,0) -- (-2.55554,0); \draw[thick] (-2.44443,-1) -- (-2.44443,.35); \node (8/27) at (-2.44443,-1.4) {$\frac{8}{27}$};
\draw[very thick] (-2.77776,0) -- (-2.88887,0); \draw[thick] (-2.88887,-.35) -- (-2.88887,.35); \node (7/27) at (-2.88887,.75) {$\frac{7}{27}$};
\draw[fill=blue!30] (-5.33335,0) circle (.1111); \draw[very thick] (-5.11113,0) -- (-5.22224,0); \draw[thick] (-5.11113,-.35) -- (-5.11113,.35); \node (2/27) at (-5.11113,.75) {$\frac{2}{27}$};
\draw[very thick] (-5.44446,0) -- (-5.55557,0); \draw[thick] (-5.55557,-1) -- (-5.55557,.35); \node (1/27) at (-5.55557,-1.4) {$\frac{1}{27}$};
\end{tikzpicture}
\caption{The annulus $A^{(2)}_C$ with countably infinitely many holes which accumulate on a Cantor set.} \label{fig:ctbl_punc_disk_cantor_set_accum}
\end{figure}
\end{example}

One does not always have to replace the deleted intervals used to build the Cantor set with copies of the same space at every level in order to obtain a continuum which non-trivially self covers.  One can choose a collection $\{K_0,K_j,\ldots, K_n\}$ of non-homeomorphic continua and an infinite collection $\{\{K_{0j},K_{1j},\ldots, K_{nj}\}:j\in \mathbb{N}\}$, where $K_{ij}$ is homeomorphic to $K_j$.  One then attaches copies of these of the appropriate diameter between pairs of rational points in the Cantor set via a shuffling process: attach one copy $K_{01}$ of $K_0$ of diameter 1 to $1$ and $2=0$; attach one copy $K_{11}$ of $K_1$ of diameter $\frac{1}{3}$ to $\frac{1}{3}$ and $\frac{2}{3}$; attach two copies $K_{21}$,$K_{22}$ of $K_2$ of diameter $\frac{1}{9}$ to, respectively, the pairs $\{\frac{1}{9},\frac{2}{9}\}$ and $\{\frac{7}{9},\frac{8}{9}\}$; attach 4 copies $K_{31},K_{32},K_{33},K_{34}$ of $K_3$ of diameter $\frac{1}{27}$ to the pairs $\frac{i}{27}$ and $\frac{i+1}{27}$ of rational points in the Cantor set; $\ldots \ldots$ ; attach $2^{n-1}$ copies $K_{n1},K_{n2},\ldots,K_{n,2^{n-1}}$ of $K_n$ of diameter $\frac{1}{3^n}$ to the pair of rational points $\frac{i}{3^n}$ and $\frac{i+1}{3^n}$ in the Cantor set; then at the $n+$first level start over by sewing in $2^{n+1}$ copies of $K_0$; at the $n+$second level sew in $2^{n+2}$ copies of $K_1$; etc.

Two simple examples of this sort of process were brought to my attention by Logan Hoehn in the Dynamics and Continuum Theory session of the 2025 Summer Conference on Topology.  The first is a circle-like Knaster Cup-cap continuum.  See Figure \ref{fig:CircleCup-cap}.  It is constructed as follows.

Consider the middle thirds Cantor set $C\times 0 \subset [0,1]\times 0$ in the circle $[0,2]\times 0/((2,0)=(0,))$ but sitting in the annulus $[0,2] \times [-1,1]/~$ formed by gluing $0\times [-1,1]$ to $2\times [-1,1]$ via the identity on the second factor.  Form Cantor's Fence $C \times [-1,1]$ inside this annulus.  At the $0^{\textrm{th}}$ level, glue in the segment $[1,2]\times -1$.  This forms the \emph{cup} $1\times [-1,1] \cup [1,2]\times -1 \times 2 \cup [-1,1]$.  At the $1^{\textrm{st}}$ level, glue in the segment $[\frac{1}{3},\frac{2}{3}] \times 1$ to form the \emph{cap} $\frac{1}{3} \times [-1,1] \cup [\frac{1}{3},\frac{2}{3} \times 1 \cup \frac{2}{3}\times [-1,1]$.  At the $2^{\textrm{nd}}$ level, in the two intervals from $[\frac{i}{9}\times, \frac{i+1}{9}]\to\times -1$ determined by the rational points of the Cantor set with denominator 9, form 2 cups by gluing in segments of length  $\frac{1}{9}$.  At the $3^{\textrm{rd}}$ level, in the four intervals $[\frac{i}{27},\frac{i+1}{27}] \times 1$ determined by rational points of the Cantor set with denominator 27, form four caps by gluing in segments of length $\frac{1}{27}$.  Continue this process, forming the appropriate number of cups and caps in alternate levels.  Note that the widths of the cups and caps decrease to 0 in levels, but their diameters do not, and that the control set for this construction is $C\times [-1,1]$.  This Cup-cap space has $\mathbb{Z}_m$-regular self covers for every $m \in \mathbb{N}$.  Both the Cup-cap space and its $\mathbb{Z}_3$-regular cover are illustrated in Figure \ref{fig:CircleCup-cap}.  Note that the classical Knaster Cup-cap continuum omits the line segment forming the $0^{\textrm{th}}$-level cup on the interval $[1,2]$.  A quotient space of the Knaster Cup-cap space is the Knaster V-$\Lambda$ continuum formed by shrinking the line segments forming the top of each cap and bottom of each cup to a point.  Applying the same quotienting process to the circle-like cup-cap space produces a circle-like V-$\Lambda$ continuum.  The quotienting process is compatible with the $\mathbb{Z}_m$-regular self covers which the circle-like Cup-cap space has and, consequently, the circle-like V-$\Lambda$ continuum has $\mathbb{Z}_m$-regular self covers for all $m\in \mathbb{N}$ as well.  See Kuratowski, \cite{Kuratowski1968}, page 191, for descriptions of the original Knaster Cup-cap and V-$\Lambda$ continua.  These continua are also called the Cantor Organ and Cantor Accordion, respectively, in the literature.

\begin{example}\label{ex:CountPuncAnnTangent}\emph{Cantor's Annulus.}
Begin with the annulus formed by gluing the ends $0\times [-1,1]$ and $2\times [-1,1]$ of the rectangle $[0,2]\times [-1,1]$ together via $(0,y)\mapsto (2,y)$ and consider the middle thirds Cantor set $C\subset [0,1]\times 0$. Delete the interior of a disk centered at $(\frac{3}{2},0)$, tangent to the points $(1,0)$ and $(2,0)=(0,0)$, and contained in the interior of the annulus.  Then, between each pair $(\frac{i}{3^j},0)$ and $(\frac{i+1}{3^j},0)$ of points in $C$, choose the disk of radius $\frac{1}{2\cdot 3^j}$ and centered at the midpoint $(\frac{2i+1}{3^j},0)$ of the interval $[\frac{i}{3^j},\frac{i+1}{3^j}]\times 0$. Delete the interiors of all these disks. Call the resulting countable punctured annulus $A^{(1)}_C$.  Note that all the rational points in $C$ are contained in the boundaries of the disks and so are in $A^{(1)}_C$.  Also note that diameters of these removed disks decrease in layers to 0:  there is 1 disk of diameter 1 which is removed, one disk of diameter $\frac{1}{3}$ which is removed, 2 disks of diameter $\frac{1}{9}$ which are removed, and in the $j^{\textrm{th}}$ level, $j \geq 1$, there are $2^{j-1}$ disks of diameter $\frac{1}{3^j}$ which are removed.  The Cantor set $C$ is again functioning as a control set in this example: it controls how the holes in the countably infinitely punctured annulus $A^{(1)}_C$ are accumulating.  See Figure \ref{fig:CountPuncAnnTangent}.

This punctured annulus $A^{(1)}_C$ has $\mathbb{Z}_m$-regular covers $f_m:A^{(1)}_C\rightarrow A^{(1)}_C$ for every $m\in \mathbb{N}$.  The $\mathbb{Z}_3$ self cover is illustrated in Figure \ref{fig:CountPuncAnnTangent}.  The reader is invited to write out the formula for $f_3$ suggested by the picture and the general formula for $f_m$.  There are recursive approaches to defining $\mathbb{Z}_m$-regular self covers of two 2-complexes related to $A^{(1)}_C$ in \cite{BedDelTimm2006}, Examples 5.1 and 5.2.
\end{example}

\begin{example}{\emph{(Cantor's Annulus with Pearls)}.}\label{ex:CantorAnnulusPearls}
Example \ref{ex:CountPuncAnnTangent} can be built upon in several ways.  First, take a countable collection $\mathcal{S}$ of pairs $(S^2,S^1)$ of 2-spheres, and its equatorial $S^1$, with diameters which decrease to 0 in levels.  Attach these $S^2$'s to the boundaries of the holes of the appropriate diameter in $A^{(1)}_C$ via a homeomorphism from the equatorial $S^1$ to the boundary of the hole.   This produces a continuum $A^{(1)}_{C,S^2}$ containing countably infinitely many $S^2$'s with $C$ serving as the control set for how the $S^2$'s accumulate. $A^{(1)}_{C,S^2}$ has $\mathbb{Z}_m$-regular self covers for all $m\in \mathbb{N}$.  Again note that this sort of construction can be carried out more generally.  One can attach copies of any preferred continuum $K$ to the boundaries of the holes in Cantor's Annulus provided the diameters of the copies of $K$ decrease to 0 in levels and the attaching process is done in a consistent manner from one copy of $K$ to the next.   The continuum $A^{(1)}_{C,K}$ has $\mathbb{Z}_m$-regular self covers for all $m\in \mathbb{N}$.  More details about this sort of generalized process using the Cantor set as a control set are in \cite{Timm2020}. Also note that there is a strong deformation retract of $A^{(1)}_{C,S^2}$ onto Cantor's Pearl Necklace $\mathcal{N}^{(2)}_{C,S^2}$.  Consequently, they have the same homotopy type.  We again ask what the structure of $\pi_2(\mathcal{N}^{(2)}_{C,S^2})$ is.
\end{example}

What happens if, instead of deleting the interiors of a collection of disks whose boundaries contain the rational points of the Cantor set $C$ as in $A^{(1)}_C$, one deletes the interiors of slightly smaller disks which lying between the pair of points $(\frac{i}{3^j},0)$ and $(\frac{i+1}{3^j},0)$ of $C$?  For example, for each pair $(\frac{i}{3^j},0)$ and $(\frac{i+1}{3^j},0)$ of $C$, delete the interior of a disk with center at the midpoint of each $[\frac{i}{3^j},\frac{i+1}{3^j}]$ and with diameter equal to one half its length.  Call the resulting space $A^{(2)}_{C}$.  Observe that $A^{(2)}_C$ can also be obtained by attaching a collection of copies of small annuli $K=S^1\times I$ with diameters which decrease to 0 in levels to the holes in $A^{(1)}_C$.  We can, therefore write, $A^{(2)}_C = A^{(1)}_{C,S^1\times I}$.

Consequently, it follows from Example \ref{ex:CantorAnnulusPearls} that the space $A^{(2)}_C = A^{(1)}_{C,S^1\times I}$ has $\mathbb{Z}_m$-regular self covers for all $m\in \mathbb{N}$.  The Cantor set is again the control set. It controls how the holes in this infinitely punctured annulus accumulate.

\begin{figure}
\centering
\includegraphics[width=5in]{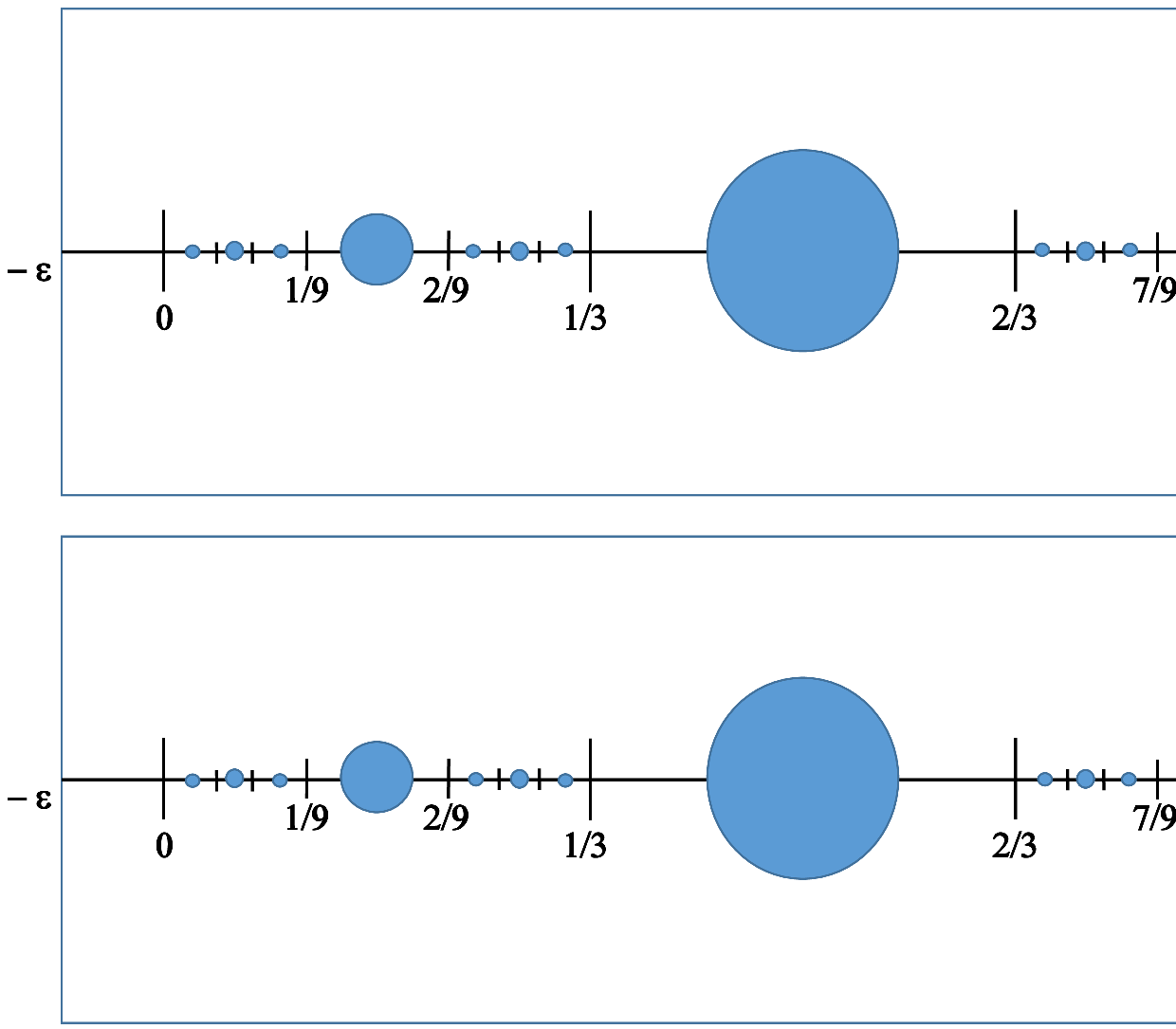}\\
\caption{Two disks with countably many holes and control set $C$.  The top disk with holes is $\mathcal{D}^{(1)}_C$, the bottom one is $\mathcal{D}^{(2)}_C$. Illustrations from \cite{Timm2020}.}\label{fig:2CountPuncDisks}
\end{figure}

\bigskip
We now look at two situations which, on first consideration, appear to be different from each other and appear to be fundamentally different from the annuli $A^{(i)}_C$, $i \in \{1,2\}$.  These are two disks $\mathcal{D}^{(i)}_C$, $i\in\{1,2\}$ with countably many holes where the boundaries of the holes don't contain rational points of $C$ but, how the holes accumulate, is still controlled by $C$ as in $A^{(2)}_C$. The $\mathcal{D}^{(i)}_C$ are illustrated in Figure \ref{fig:2CountPuncDisks}.

\begin{figure}
  \centering
  \includegraphics[width=4.5in]{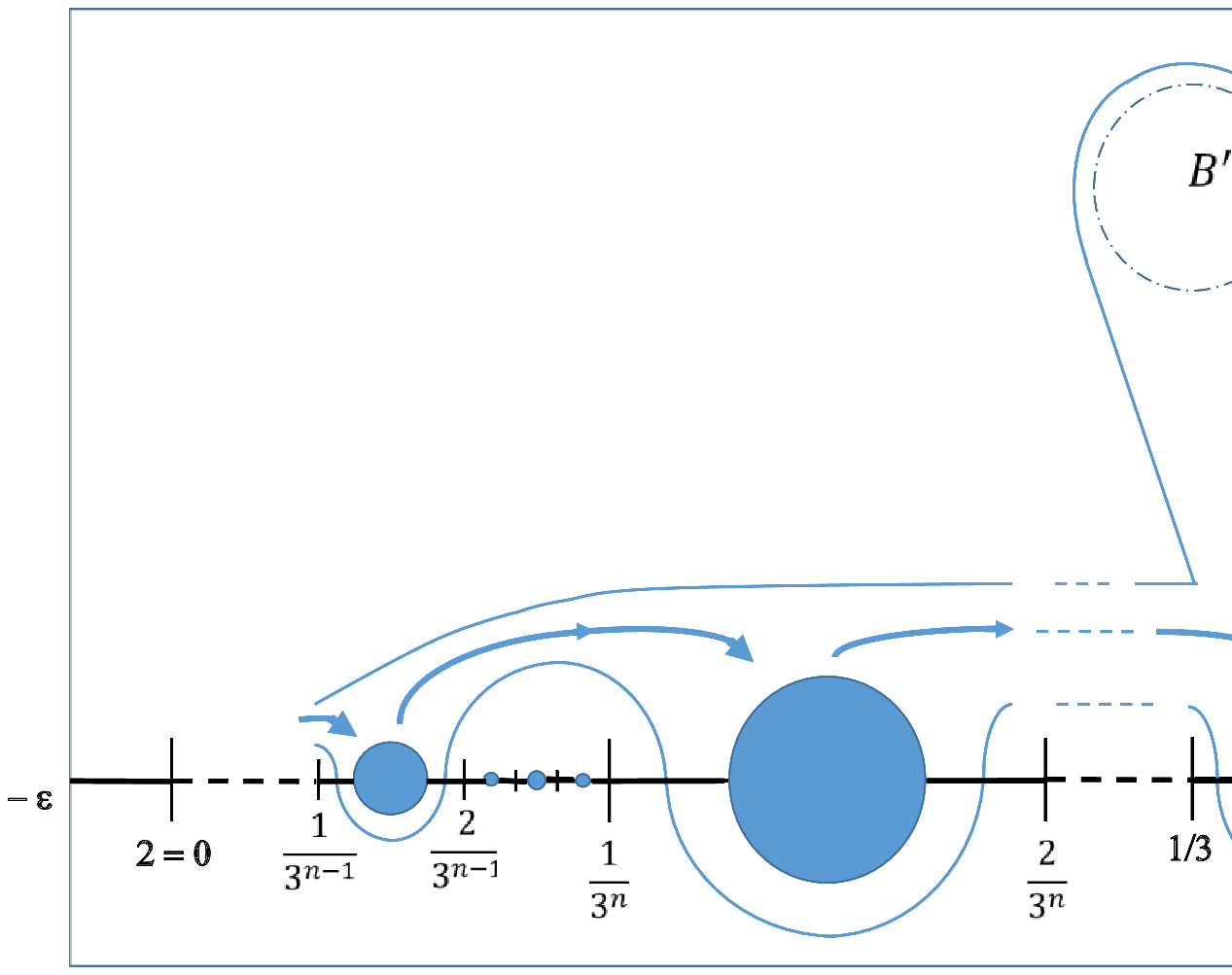}\\
  \caption{The homotopies $f_t$ and $h_t$.  Illustration from \cite{Timm2020}.}\label{fig:CountPuncDiskHtpys}
\end{figure}

It turns out that these two disks with countably many holes are homeomorphic, in fact, they are ambiently isotopic.  Perhaps surprisingly, they are just different representations of the countably puncture annulus $A^{(2)}_C$.  The ambient isotopies $f_t$ and $h_t$ whose end maps $f_1$ and $h_1$ are homeomorphisms of $\mathcal{D}^{(1)}_C$ onto, respectively $\mathcal{D}^{(2)}_C$ and $A^{(2)}_C$ are shown in Figure \ref{fig:CountPuncDiskHtpys}.  Some explanation is needed.

Let $S^1_n=\{(x,y):(x-\frac{1}{2 \cdot 3^{n-1}})^2+y^2=(\frac{1}{2\cdot 3^n})^2\}.$  That is, $S^1_n$ is the boundary of the hole which has its center at the midpoint of the interval $[\frac{1}{3^n},\frac{2}{3^n}]$ used in the construction of the the Cantor set $C$.  Choose a closed neighborhood $N$ of $\overset{\infty}{\underset{n=1}{\bigcup}}S^1_n$ such that

\begin{enumerate}
  \item $S^1_n \in N$ for all $n\in \mathbb{N}$,
  \item $(0,0) \in \partial N$, and
  \item all the points of $C\setminus \{(0,0)\}$ are in $D\setminus{N}$.
\end{enumerate}

\begin{figure}
  \centering
  \includegraphics[width=4.5in]{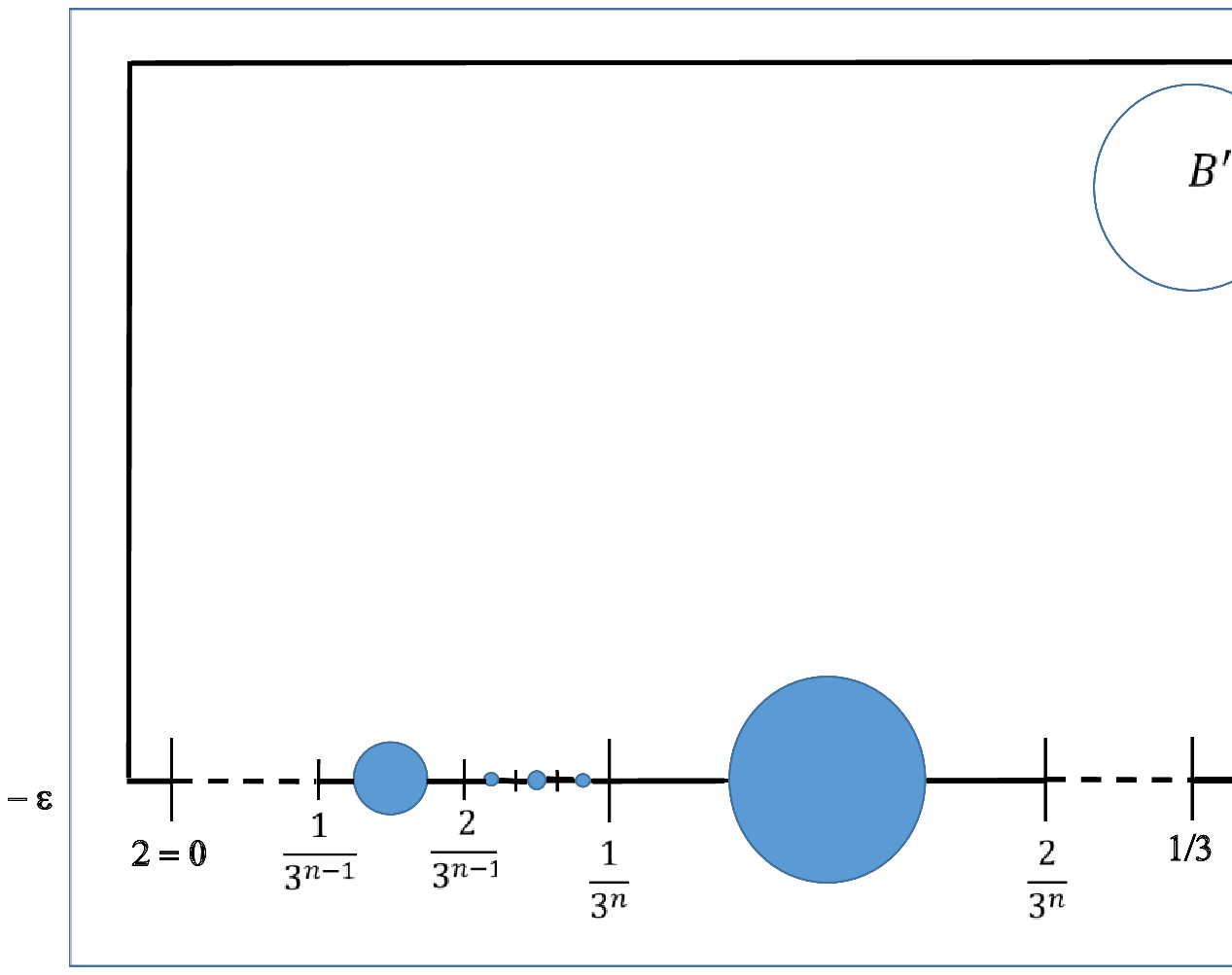}\\
  \caption{$D^{(1)}_C \cong D^{(2)}_C$ as $A^{(1)}_C$ with $B^{\prime}_{00}$ as the central hole in $A^{(1)}_C$.  The inner
  rectangle in bold represents the circle $[0,2] \times 0 / \big( (0,0)\sim (2,0) \big)$ of $A^{(1)}_C$.  The illustration is from \cite{Timm2020}.}\label{DiskAsAnnulus}
\end{figure}

\noindent The disks $B_{00}$ and $B^{\prime}$ shown with dotted outlines in Figure \ref{fig:CountPuncDiskHtpys} and unfilled circles in Figure \ref{fig:CountPuncDiskHtpys} are the images of the hole surrounded by $S^1_1$ under the end maps of the isotopies $f_t$ and $h_t$, respectively.  The isotopy $f_t$ shifts $S^1_1$ onto the boundary $\partial B_{00}$ and, for $n\geq 2$, shifts $S^1_{n}$ onto $S^1_{n-1}$ via trajectories in $int(N)$.  That is, $f_t$ is the ambient isotopy which changes $\mathcal{D}^{(1)}_C$ into $\mathcal{D}^{(2)}_C$.  Follow $f_t$ by the second isotopy $h_t$ which shifts $\partial B_{00}$ to $\partial B^{\prime}_{00}$, again shifts $S^1_1$ onto $\partial B_{00}$, and again, for $n \geq 2$, shifts $S^1_n$ onto $S^1_{n-1}$ all along trajectories contained in $int(N)$.  On $[-\epsilon,2+\epsilon]\times[-1,1]\setminus int(N)$, both $f_t$ and $h_t$ are the identity for all $t$.  The map $h_t$ is the isotopy showing that disk $\mathcal{D}^{(2)}_C$ with countably many holes and the annulus $A^{(2)}_C$ with countably many holes are the same space.  The representation of $\mathcal{D}^{(2)}_C$ as $A^{(2)}_C$ is in Figure \ref{DiskAsAnnulus}.  The rectangle outlined in thick black in the interior of the larger rectangle $[-\epsilon,2+\epsilon]\times [-1,1]$ represents the center circle of the annulus $A^{(2)}_C$.

Since $\mathcal{D}^{(1)}_C$ is homeomorphic to $A^{(1)}_C$, is is clear that $\mathcal{D}^{(1)}_C$ has $\mathbb{Z}_m$-regular self covers for all $m \in \mathbb{N}$.  In fact, $\mathcal{D}^{(1)}_C \cong A^{(2)}_C$ has many more regular self covers.

\begin{theorem}\cite{DelgadoTimm2017},\cite{Timm2020}
$\mathcal{D}^{(1)}_C \cong A^{(2)}_C$ has $G$-regular self covers for every finite group $G$.
\end{theorem}

The proof depends on several observations.  The first is that a finite group is the quotient of a free group on some minimum number, say $q$, generators by a finite index normal subgroup of index $|G|$, which is necessarily another finitely generated free group.  Second, this quotient can be realized as the group of deck transformations for a $G$-regular, degree $|G|$, covering space $f:D_p\rightarrow D_q$ of a disk $D_q$ with $q$ holes by a disk $D_p$ with $p\geq q$ holes.  Third, the homotopy $h_t$ can be applied repeatedly to obtain disks with countably many holes controlled by $C$ with any finite number of holes which are ``far'' from the $x$-axis in $[-\epsilon,2+\epsilon]\times [-1,1]$.

Now, let $V$ be a small disk contained in the interior of a fundamental domain of $f$.  Then $U=f(V)$ is a small disk in $D_q$ which is evenly covered by $f$. In $int(U)$ construct a copy of $D^{(1)}_C$ by specifying a nicely embedded copy $J\subset int(U)$ of the interval $[0,2]$.  Next,  construct a copy $C$ of the middle thirds Cantor set in $[0,1]\subset J$.  Finally, use $C$ as the controls set to form the holes in this small copy of $D^{(1)}$ contained in $int(U)$. Let $\mathcal{D}$ denote the collection of the disks with diameters decreasing to 0 in levels whose interiors are deleted to form the holes. The third observation in the preceding paragraph $D_q\setminus (\cup\{int(D)\in \mathcal{D}\})$ is also a copy of the disk $D^{(1)}_C$ with countably many holes controlled by a Cantor set.  Because $U$ is evenly covered by $f$, so is $U\setminus \{int(D):D\in \mathcal{D}\}$.  Consequently, $f^{-1}(U\setminus (\cup\{int(D):D\in \mathcal{D}\})$ is a union of $|G|$ copies of $D^{(1)}_C$ all of which are ``far'' from the original $p$ holes in $D_p$. Since the union of the $|G|$ copies of the Cantor set $C$ contained in $f^{-1}\left(U^{}\setminus (\cup\{int(D):D\in \mathcal{D}\})\right)$ is also a Cantor set, apply the third observation in the preceding paragraph a second time and it follows that $D_p\setminus f^{-1}\left(D^{}_q\setminus (\cup\{int(D):D\in \mathcal{D}\})\right)$ is a copy of $\mathcal{D}^{(2)}_C$.  Consequently the restriction
$$
f:D_p\setminus f^{-1}\left(D^{}_q\setminus (\cup\{int(D):D\in \mathcal{D}\})\right)\longrightarrow D_q\setminus \left(\cup\{int(D):D\in \mathcal{D}\}\right) $$

\noindent is a $G$-regular self cover. In particular $\mathcal{D}^{(2)}_C$ has $G$-regular self covers for every finite nonabelian group $G$.  More details of the above are in \cite{Timm2020}.

In the same paper it is observed that by fixing a continuum $K$ and carefully attaching copies of it with diameters which decrease to 0 in levels to the boundaries of the holes in $\mathcal{D}^{(1)}_C$, one can obtain additional examples of continua $\mathcal{D}^{(1)}_{C,K}$ with $G$-regular covers for certain choices of $K$ and finite group $G$.  A complication is that for a particular self cover $f:\mathcal{D}^{(1)}_C\rightarrow \mathcal{D}^{(1)}_C$ the degrees of the restrictions of $f$ to different circles in $\partial \mathcal{D}^{(1)}_C$ can vary.  For example, if one wants to construct a $\mathcal{D}^{(1)}_{C,K}$ by starting with one of the standard double covers $f:D_3\rightarrow D_2$, of the disk $D_2$ with two holes by the disk $D_3$ with three holes, one cannot then choose $K$ to be the Mobius band glued to the boundaries of the holes via homeomorphisms from $\partial K$ to the boundaries of holes.  The problem here is that the boundary $S^1$ of one of the holes in $D_3$ must map to the boundary of one of the holes in $D_2$ via a double cover and the Mobius band has no double self cover.  On the other hand, $f:D_3\rightarrow D_2$ can be used to construct a $D^{(1)}_{C,K}$ which will have a degree 2 self cover if one chooses $K$ to be an annulus and attaches copies of it to the boundaries of the holes vie a homeomorphism from their center circles to the boundary of the hole to which each is attached.

Finally, the disk with holes $\mathcal{D}^{(1)}_C$ is a 2-dimensional space with $G$-regular self covers for every finite group $G$. The above process used to show this can be thought of as the first step needed to show that an important 1-dimensional continuum has $G$-regular self covers for every finite group $G$.

Let $G$ be a finite group and $f:\mathcal{D}^{(1)}_C\rightarrow \mathcal{D}^{(1)}_C$ a $G$-regular self cover.  Cover the range copy of $\mathcal{D}^{(1)}_C$ by a collection $\mathcal{D}$ of closed disks with pairwise disjoint interiors such that $\cup\{D:D\in \mathcal{D}\}$ is dense in $D^{(1)}_C$.  Delete the union of the interiors of the disks.  By Whyburn \cite{Whyburn1958}, the resulting continuum $\mathcal{S} = D^{(1)}_C\setminus (\cup\{D:D\in \mathcal{D}\})$ is a copy of the Sierpinski Carpet.  It is also the case that $f^{-1}(\mathcal{S})$ is a copy of the Sierpinki Carpet.  Consequently, it follows that

\begin{theorem}
The Sierpinski Carpet has $G$-regular self covers for every finite group $G$.
\end{theorem}

\noindent See \cite{DelgadoTimm2017} for more detail.

\begin{question}
Are all finite sheeted covers of the disk $\mathcal{D}^{(1)}_C$ self covers? What about for the Sierpinski Carpet?
\end{question}

\section{Inverse limits and self covering}\label{sec:InvLimAndSelfCov}
Fox based some of his work on the solenoids on Borsuk's work \cite{Borsuk1968} on shape theory.  Their initial contributions to shape theory have been refined, extended, and applied in various ways by large numbers of mathematicians.  Among those most relevant to the current discussion are those by various combinations of Eda, Kocieć Bilan, Mandić, Mardešić, Matijević, and Uglešić, \cite{EdaEtAl2005}, \cite{EdaMatijevic2006}, \cite{EdaMatijevic2008}, \cite{KoceicBilan2014}, \cite{KoceicUglesic2007}, \cite{MardesicMatijevic2001}, \cite{Matijevic2003}, \cite{Matijevic2007}.  This body of work investigates the topology and group theory of various generalized solenoids, including their covering space properties.  For example, in Eda et al \cite{EdaEtAl2005}, it is shown that covering spaces of limits of inverse sequences $\underset{\leftarrow}{\textrm{\textbf{lim}}}(T^2,f_n,\mathbb{N})$, $T^2=S^1\times S^1$, with bonding maps $f_n:T^2\rightarrow T^2$ self covers for all $n\in\mathbb{N}$, are also inverse limits of the same sorts of inverse sequences, but that these covering spaces are not necessarily self covers.  The structure of the course shape groups of these sorts of inverse limits, in particular, of the solenoids, is investigated in Kocieć Bilan \cite{KoceicBilan2014}.  It is worth noting that the Cantor set makes an appearance here.  The classical solenoids are Cantor set bundles over $S^1$.  We need a few more definitions to continue this discussion.

\begin{figure}[t] 
\centering
    \begin{tikzpicture}
    \node (C1) at (0,-2) {$S^1$};
    \node (C2) at (2,-2) {$S^1$};
    \node (C3) at (4,-2) {$S^1$};
    \node (C4) at (6,-2) {$S^1$};

    \begin{scope}[every path/.style={-triangle 60}]
       \draw (1.6,-2) -- (.4,-2);
       \draw (3.6,-2) -- (2.4,-2);
       \draw (5.6,-2) -- (4.4,-2);
       \draw (7.2,-2) -- (6.4,-2);
    \end{scope}
       \draw [dashed] (7.2,-2) -- (9,-2);
    \node (m1) at (1.1,-1.75) {$\times m_1$};
    \node (m2) at (3.1, -1.75) {$\times m_2$};
    \node (m3) at (5.1, -1.75) {$\times m_3$};
    \node (m4) at (7.1,-1.75) {$\times m_4$};
    \node (b) at (4.5,-2.75) { };
    \node (solenoid) at (9.5,-2) {$ : \Sigma_{\textrm{\textbf{m}}}$};
\end{tikzpicture}
\caption{The classical solenoid. The notation is further simplified in the picture.  The $\times m_k$ notation over the arrows denotes the degree of the self cover.}\label{fig:ClassicSolenoid}
\end{figure}

An \emph{inverse system} consists of the following. (1)A partially ordered set $I$ with partial order $\leq$.  (2)For each $i\in I$ a topological space $X_i$.  (3) For each pair $i\leq j$ in $I$ a function $g^j_i:X_j\rightarrow X_i$ such that $g^i_i = id_{X_i}$ and if $i\leq j \leq k$ then $g^k_i=g^k_j\circ g^j_i$.  The inverse system is denoted by the triple $(X_i,g^j_i,I)$ or when the index set is understood, by the pair $(X_i,g^j_i)$.  The functions $g^j_i$ are the \emph{bonding maps} of the system.  When the index set for the inverse system is the natural numbers $\mathbb{N}$ with its usual ordering, the system is called an \emph{inverse sequence} and the notation is simplified.  The bonding maps between consecutive spaces in an inverse sequence are denoted by $g_n:X_{n+1}\rightarrow X_n$.  The bonding maps from $X_{n+k}$ onto $X_n$ are the compositions $g_{n+k-1}\circ \cdots \circ g_{n+1} \circ g_n$.  They are typically not listed in the notation for the inverse system.

The \emph{inverse limit} of the inverse system is denoted $\underset{\leftarrow}{\textrm{\textbf{lim}}}(X_i,g^j_i,I)$ and is the subspace of the product $\underset{i\in I}{\times}X_i$, with the product topology, defined by
$$
\underset{\leftarrow}{\textrm{\textbf{lim}}}(X_i,g^j_i,I)=\{(x_i)_{i\in I}\in \underset{i\in I}{\times}X_i: x_i=g^j_i(x_j) \textrm{ for all } j\in I \textrm{ with } i\leq j\}
$$

\noindent together with a collection $\Pi$ of projection maps $\pi_i:\underset{\leftarrow}{\textrm{\textbf{lim}}}(X_i,g^j_i,I)\rightarrow X_i$ which, for $i\leq j$ satisfy $\pi_i=g^j_i\circ \pi_j.$  The points $(x_i)_{i\in I} \in \underset{\leftarrow}{\textrm{\textbf{lim}}}(X_i,g^j_i,I)$ are called \emph{threads}.  For our purposes\footnote{A more general definition is possible.  See Engelking \cite{Engelking1989}} a \emph{map} of inverse systems, denoted $\Phi:(X_i,g^i_j,I)\rightarrow (Y_i,h^i_j,I)$ is a collection $\Phi$ of functions $\varphi_i:X_i\rightarrow Y_i$ such that, if $i\leq j$, then $\phi_i\circ h^j_i = g^j_i\circ \phi_j$.  A map of inverse systems $\Phi:(X_i,g^i_j,I)\rightarrow (Y_i,h^i_j,I)$ induces a map $\widetilde{\Phi}:\underset{\leftarrow}{\textrm{\textbf{lim}}}(X_i,g^j_i,I)\rightarrow \underset{\leftarrow}{\textrm{\textbf{lim}}}(Y_i,h^j_i,I)$.

Note that when all $X_i$ are Hausdorff, the inverse limit is a closed subset of the Cartesian product.  The proof of this and other basic properties of inverse systems can be found in Engelking \cite{Engelking1989}, primarily in Section 2.5, but also scattered throughout the rest of the text.

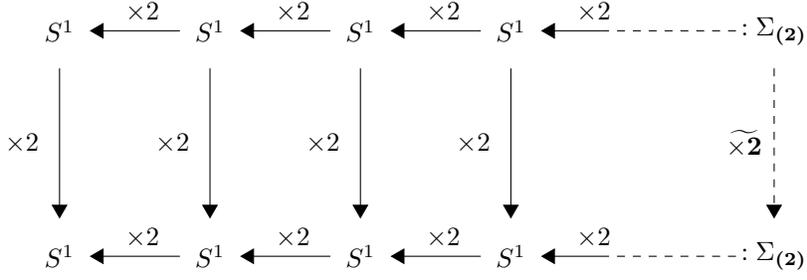
\begin{figure}[t]
\centering
\noindent \begin{tikzpicture}
    \node (C1) at (0,1) {$S^1$};
    \node (C2) at (2,1) {$S^1$};
    \node (C3) at (4,1) {$S^1$};
    \node (C4) at (6,1) {$S^1$};

    \begin{scope}[every path/.style={-triangle 60}]
       \draw (1.6,1) -- (.4,1);
       \draw (3.6,1) -- (2.4,1);
       \draw (5.6,1) -- (4.4,1);
       \draw (7.2,1) -- (6.4,1);
    \end{scope}
       \draw [dashed] (7.2,1) -- (9,1);
    \node (m1) at (1.1,1.25) {$\times 2$};
    \node (m2) at (3.1,1.25) {$\times 2$};
    \node (m3) at (5.1,1.25) {$\times 2$};
    \node (m4) at (7.1,1.25) {$\times 2$};
    \node (b) at (4.5,-2.75) { };
    \node (solenoid) at (9.5,1) {$ : \Sigma_{\textrm{\textbf{(2)}}}$};

    \node (C1) at (0,-2) {$S^1$};
    \node (C2) at (2,-2) {$S^1$};
    \node (C3) at (4,-2) {$S^1$};
    \node (C4) at (6,-2) {$S^1$};
    \begin{scope}[every path/.style={-triangle 60}]
       \draw (1.6,-2) -- (.4,-2);
       \draw (3.6,-2) -- (2.4,-2);
       \draw (5.6,-2) -- (4.4,-2);
       \draw (7.2,-2) -- (6.4,-2);
    \end{scope}
       \draw [dashed] (7.2,-2) -- (9,-2);
    \node (m1) at (1.1,-1.75) {$\times 2$};
    \node (m2) at (3.1, -1.75) {$\times 2$};
    \node (m3) at (5.1, -1.75) {$\times 2$};
    \node (m4) at (7.1,-1.75) {$\times 2$};
    \node (lim) at (9.5,-2) {$ : \Sigma_{\textrm{\textbf{(2)}}}$};
  \begin{scope}[every path/.style={-triangle 60}]
    \draw (0,.5) -- (0,-1.5);
    \draw (2,.5) -- (2,-1.5);
    \draw (4,.5) -- (4,-1.5);
    \draw (6,.5) -- (6,-1.5);
    \draw[dashed] (9.5,.5) -- (9.5,-1.5);
  \end{scope}
    \node (d1) at (-.5,-.5) {$\times 2$};
    \node (d2) at (1.5,-.5) {$\times 2$};
    \node (d3) at (3.5,-.5) {$\times 2$};
    \node (d4) at (5.5,-.5) {$\times 2$};
    \node (dm) at (9.1,-.5) {$\widetilde{\times \textbf{2}}$};
\end{tikzpicture}
\caption{The inverse sequence $(S^1,\times 2, \mathbb{N})$, its inverse limit $\Sigma_{\textrm{\textbf{(2)}}}$, and the $\times \textbf{2}:\left(S^1,\times 2, \mathbb{N}\right)\rightarrow \left(S^1,\times 2, \mathbb{N}\right)$ self map of the inverse sequence. Note that $\widetilde{\times \textbf{2}}:\Sigma_{\textrm{\textbf{(2)}}}\rightarrow \Sigma_{\textrm{\textbf{(2)}}}$ is a homeomorphism.} \label{fig:The2SolenoidSelfMap}
\end{figure}

\begin{example}\label{TheClassicalSolenoids} \emph{The classical solenoids.} Let $S^1=\{e^{i\cdot\theta}\in \mathbb{C}:0\leq \theta\leq 2\pi\}$.  For $k\in \mathbb{N}$, let $m_k \in \mathbb{N}$ and $\textrm{\textbf{m}}=(m_1,m_2,m_3,\ldots)$.  Define $g_k:S^1 \rightarrow S^1$ by $g_k(e^{i\cdot \theta})=e^{i\cdot m_k \theta}$.  Let $\Sigma_{\textrm{\textbf{m}}}=\underset{\leftarrow}{\textrm{\textbf{lim}}} (S^1,g_k,\mathbb{N})$.  The inverse limit
$\Sigma_{\textrm{\textbf{m}}}$ is a classical \emph{solenoid}.  A picture of this situation is in Figure \ref{fig:ClassicSolenoid}.

\bigskip
Let $p$ be a prime.  Define a self map of the solenoid $\widetilde{\Phi}: \Sigma_{\textrm{\textbf{m}}} \rightarrow \Sigma_{\textrm{\textbf{m}}}$ by letting the map of the inverse sequence $\Phi:(S^1,f_m,\mathbb{N})\rightarrow (S^1,g_k,\mathbb{N})$ be given by $\varphi_k(e^{i\cdot \theta})=e^{i\cdot p\theta}$, the standard degree $p$ self cover of the circle.  It is well known that the induced map $\widetilde{\Phi}:\Sigma_{\textrm{\textbf{m}}}\rightarrow \Sigma_{\textrm{\textbf{m}}}$ is a covering space map. When $p$ divides at most finitely many of the degrees $m_i$ of the bonding maps in the sequence, the induced map $\widetilde{\Phi}:\Sigma_{\textrm{\textbf{m}}} \rightarrow \Sigma_{\textrm{\textbf{m}}}$ is, in fact, a $\mathbb{Z}_p$-regular self cover.  Otherwise, it is a homeomorphism.  That is, all covers of the classical solenoids are self covers.

In certain situations it is easy to see why the induced map $\widetilde{\Phi}$ on the solenoid is a homeomorphism and not a degree $p$ cover.  Let $\textrm{\textbf{m}}=(2,2,2,\ldots)=\textrm{\textbf{(2)}}$.  Consider the the dyadic solenoid $\Sigma_{\textrm{\textbf{(2)}}}$ where all bonding maps between consecutive pairs of $S^1$'s in the sequence $(S^1,\times 2, \mathbb{N})$ are the standard degree 2 self cover of $S^1$.  Define the map of inverse sequences $\times \textbf{2}:\left(S^1,\times 2, \mathbb{N}\right)\rightarrow \left(S^1,\times 2, \mathbb{N}\right)$ where the map on the each circle in the sequence is also the standard 2-fold self cover.  See Figure \ref{fig:The2SolenoidSelfMap}.  In the figure, trace the diagram along the top row from right to left, starting at $\Sigma_{\textrm{\textbf{(2)}}}$, all the way to the left side of the diagram, then down the left edge of the diagram.  Doing so, it should be clear that the induced map $\widetilde{\times \textbf{2}}$ is a coordinate change on $\Sigma_{\textrm{\textbf{(2)}}}$.
\end{example}

There are various ways these classical results can be generalized.  The first is to ask about covering space properties of the inverse limit of an inverse system of circles when the index set is arbitrary.  Using the shape theoretic techniques developed in Mardsić and Matijeveć \cite{MardesicMatijevic2001} to investigate properties of a special type of covering space called an \emph{overlay}, see e.g., Fox \cite{Fox1972},\cite{Fox1974} and Moore \cite{Moore1978}, Matijevi\'c \cite{Matijevic2003} proves that all finite sheeted covers of solenoids $\underset{\leftarrow}{\textrm{\textbf{lim}}}(S^1,f^j_i,I)$ in which $I$ is an arbitrary infinite index set are self covers.

One may also consider generalizations of the classical solenoid construction by considering limits of inverse systems $(X_i, f^j_i, I)$ where the $X_i$ are copies of a space $X$, other than the circle, which has non-trivial self covers.  The natural place to begin this analysis is when the $X_i$ are copies of the $n$-torus $T^n=\overset{n}{\underset{k=1}{\times}}S^1$.  In \cite{EdaEtAl2005} Eda, Mandi\'c, and Matijevi\'c show that even for inverse limits of inverse sequences of the 2-torus, there are selections of self covers $f_n:T^2\rightarrow T^2$ as bonding maps for which $\underset{\leftarrow}{\textrm{\textbf{lim}}}(T^2,f_n,\mathbb{N})$ has finite index covers which are not self covers.  They use interactions between the topology of inverse limits and the number theory of the $p$-adic integers to produce their examples.  The study of these sorts of inverse systems is of particular interest, and is continued by the first and third authors in \cite{EdaMatijevic2006}, because it leads to the classification of subgroups of $\mathbb{Q}\times \mathbb{Q}$ which are finite index super groups and subgroups of torsion free abelian groups $A$ of rank $2$ in \cite{EdaMatijevic2008}. (Every torsion free abelian group of rank 2 is a subgroup of $\mathbb{Q}\times \mathbb{Q}$.)

Finally, one can generalize the classical solenoid construction by using as spaces in the inverse system, spaces for which only some of the finite sheeted covers are self covers. The analysis of these sorts of spaces is begun by Matijevi\'c in \cite{Matijevic2007} where she classifies the finite sheeted covering spaces of the inverse limit of inverse sequences in which the spaces are all copies of the Klein bottle and the bonding maps are self covers.

\begin{question}
What do the covering spaces of inverse limits of inverse systems in which the spaces are copies of the disk $D^{(1)}_C$ or the Sierpinski Gasket and the bonding maps are finite sheeted self covers of, respectively, $D^{(1)}_C$ or the Gasket, look like?  Can they be classified?
\end{question}

The shape groups and course shape groups, see Kocieć Bilan \cite{KoceicBilan2014} and Kocieć Bilan and Uglešić \cite{KoceicUglesic2007}, of the solenoids are also of interest.  Questions about how the solenoids embed in $S^3$, e.g., the knot theory of the solenoids are addressed by Connor et al in \cite{ConnorMeilRepo2015} and Jian et al in \cite{JiangEtAl2006}.  The topology of a class of 1-dimensional spaces called \emph{pseudo-solenoids} is investigated in Sturm \cite{STURM2014}.

\section{The Hawaiian Earing}\label{sec:HawaiianEarRing}
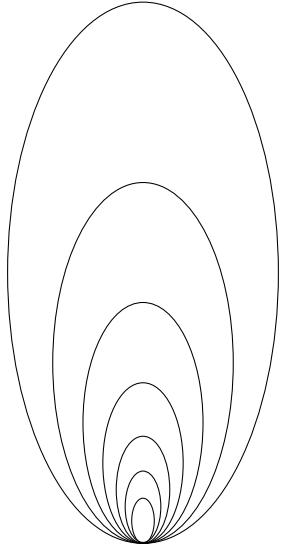
\begin{figure}
\centering
\begin{tikzpicture}[scale=.6]
\draw (0,0) ellipse (3cm and 6cm);
\draw (0,-2) ellipse (2cm and 4cm);
\draw (0,-3.33) ellipse (1.33cm and 2.67cm);
\draw (0,-4.22) ellipse (.89cm and 1.78cm);
\draw (0,-4.813) ellipse (.593cm and 1.187cm);
\draw (0,-5.189) ellipse (.395cm and .791cm);
\draw (0,-5.5) ellipse (.25cm and .5cm);
\end{tikzpicture}
\caption{ A Hawaiian Earring.}\label{HawaiianEarring}
\end{figure}

The Hawaiian Earring $\mathcal{H}$ is the one point join of countably infinitely many circles whose diameters decrease to 0 with the topology  induced by thinking of it as a subspace of Euclidean 2-space.  See Figure \ref{HawaiianEarring}.  It is of interest because it shows that in order for there to be a strict correspondence between covering spaces (with base points) and subgroups of the fundamental group of a connected, locally path connected space $Y$, $Y$ must also be \emph{semi-locally simply connected}.  The Hawaiian Earring provides the classical example of a space where this correspondence is broken: $\mathcal{H}$ is not semi-locally simply connected at the join point and it has no covering space corresponding to the trivial subgroup $\langle 1 \rangle$.  It has no such covering space because each lift under a covering projection, of every neighborhood of the join point $x_0$ in the Hawaiian Earring, contains a a copy of $\mathcal{H}$.  In addition, because the join point $x_0$ of the circles in the Hawaiian Earring must lift to $d\geq 2$ points in the domain of a degree $d\geq 2$ covering space map, it is clear that the Hawaiian Earing $\mathcal{H}$ has no non-trivial self covers.  However, $\mathcal{H}$ is of interest in the context of this paper because it does have finite degree homotopy-self covers.  That is, the total space of each finite degree cover of the Hawaiian Ear Ring has the homotopy type of $\mathcal{H}$ and consequently, its fundamental group has proper finite index clones.  To see this is not difficult.

Let $x_0$ denote the join point for the loops of $\mathcal{H}$. Let $f: X\rightarrow \mathcal{H}$ be a degree $d < \infty$ covering of the Hawaiian Earring.  $X$ has the structure of a finite graph $G$ on $d$ vertices (the $d$ lifts of $x_0$ under $f^{-1}$) with a copy of $\mathcal{H}$ at each vertex.  Choose a maximal (spanning) subtree $T$ of $G$ and contract it to a point $v_0\in f^{-1}(x_0)$.  The resulting space $X^{\prime}$ has a single 0-cell, $v_0$, with $d$ copies $H_1,H_2,\ldots, H_d$ of $\mathcal{H}$ attached to it plus a collection of finitely many, say $k$ additional loops attached to it.  These $k$ loops are formed from the non-$T$ edges of $G$ when $T$ is contracted to $v_0$.  To show $X^{\prime}$ is homeomorphic to $\mathcal{H}$, map $v_0$ to $x_0$.  Then, map the $k$ additional loops to the $k$ largest loops of $\mathcal{H}$.  Finally, map $H_1\cup H_2 \cup \ldots \cup H_d$ to $\mathcal{H}$ around $v_0$ into $\mathcal{H}$ via a shuffling processes.  That is, map the largest loop in $H_1$ to the $k+1^{\textrm{st}}$ largest loop of $\mathcal{H}$, the largest loop in $H_2$ to the $k+2^{\textrm{nd}}$ largest loop in $\mathcal{H}$, $\ldots \ldots$, the largest loop in $H_d$ to the $k+d^{\textrm{th}}$ largest loop in $\mathcal{H}$; then repeat for the collection of second largest loops from the $H_i$; and so on.

While $\pi_1(\mathcal{H},x_0)$ has no universal covering space, i.e., no covering space corresponding to certain very high index subgroups of the fundamental group, the above argument shows that the fundamental group of $\mathcal{H}$ has covering spaces corresponding to a great many (perhaps all?) of its low index subgroups.  In particular, since the above can be used to show that $\mathcal{H}$ has, up to homotopy type, $G$-regular self covers for every finite group $G$, it follows that for each finite group $G$, $\pi_1(\mathcal{H},x_0)$ has a normal finite index subgroups $N$ such the quotient group $\pi_1(\mathcal{H},x_0)/N$ is isomorphic to $G$.

\bigskip
Questions about the structure of $\mathcal{H}$ and its fundamental group have been addressed by a number of people.  First, it is an uncountable group, and accordingly, uncountable generated.  It is also not a free group. See de Smit \cite{deSmit1992}.  Cannon and Connor have a sequence of three papers \cite{CannonConnor2006}, \cite{CannonConnor2000-2},  \cite{CannonConnor2000}, which explore the combinatorial structure of the Hawaiian Earring, its fundamental group, and other uncountably generated groups they call \emph{big groups}.  In \cite{RobinsonTimm1998}, a structure theorem for finitely generated groups which have the property that all their finite index subgroups are isomorphic to the whole group is proved.  In 1999 Bob Edwards gave a series of talks titled \emph{Introduction to the Essentiality of $p$-adic Cantor Group Classifying Spaces, and the Nonexistenceof Free Cantor Group Actions on ENR's}\footnote{A \emph{Cantor group}, in this context, is any topological group whose underlying topological space is the Cantor set.} at the June Workshop in Geometric Topology \cite{Edwards1999} at the University of Wisconsin, Milwaukee.  The motivation for the talks was provided by Edwards' interest in the Hilbert-Smith conjecture.  One of the constructions to which participants in the Workshop were introduced by Bob were spaces he called the \emph{Hawaiian solenoids}.  These spaces are inverse limits of inverse systems in which the spaces are copies of the Hawaiian Earring and the bonding maps in the system are covering space maps.  McGinnis \cite{McGinnis2019} shows that when the bonding maps are regular covers, these inverse limits satisfy unique path lifting and gives a condition which assures the limits are path connected.  Fischer and Zastrow \cite{FisherZastrow2013} investigate topological and group theoretic properties of semi-coverings of $\mathcal{H}$.  Their results generalize some of our classical theory of covering spaces.

\bigskip
This collection of results suggest two questions.  Given how complicated the topology of $\mathcal{H}$ is, answers to them may be surprising.  For example, while $\mathcal{H}$ seems to not be too different from the one point join $\mathcal{S}=\underset{y_0=(0,0)}{\cup}\{S^1_i: i\in \mathbb{N}\}$, $S^1_i=\{(x,y):x^2+(y-\frac{i}{2})^2=\left(\frac{i}{2}\right)^2\}$ of a countably infinite collection of circles of increasing diameter, $\pi_1(\mathcal{H},x_0)$ is an uncountably generated non-free group, while $\pi_1(\mathcal{S},y_0)$ is a countably generated free group.

\begin{question}
Can a structure theorem similar to that of \cite{RobinsonTimm1998} be proven for big groups which have all their finite index subgroups isomorphic to the whole group?  And, of continuing interest, what more can be said about the topology and fundamental group of the Hawaiian solenoids?
\end{question}

\begin{figure}[t]
    \centering
\begin{tikzpicture}[every node/.style={transform shape}]
    \draw (0,0) circle (.25cm) (3,0) circle (.25cm) (6,0) circle (.25);
    \node (u) at (0,0) {$u$};
    \node (u) at (3,0) {$u$};
    \node (v) at (6,0) {$v$};
    \begin{scope}[every path/.style={-triangle 60}]
       \draw (3.25,0) -- (5.75,0);
       \draw (-.5,.15) ++(310:.5) arc (335:25:.5);
    \end{scope}
       \node (r) at (-1.3,-.25) {$t$};
       \node (v00w1) at (-.45,.7) {$n$};
       \node (v00w2) at (-.2,-.65) {$m$};
       \node (u30w) at (3.5,.25) {$m$};
       \node (v60w) at (5.25,.25) {$n$};
        \node (a) at (-.5,-1.25) {$(\Gamma,\omega)=[u^m,u^n]$};
        \node (b) at (4.5,-1.25) {$(\Gamma,\omega)=[u^m,v^n]$};
\end{tikzpicture}
\caption{Simple GBS graphs.} \label{fig:SingleEdge}
\end{figure}
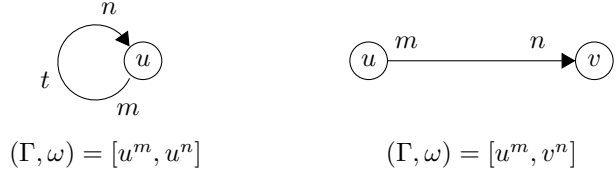

\section{Seifert fibred spaces}\label{sec:SeifertFibred}
 According to Hempel \cite{Hempel2004}, a 3 dimensional manifold $X$ is \emph{Seifert fibred} if the following conditions on $X$ are satisfied.

\begin{enumerate}
    \item There is a pairwise disjoint collection $\{S^1_{\alpha} : \alpha \in A\}$ of circles, called \emph{fibres}, such that $X = \underset{\alpha \in A}{\cup}S^1_{\alpha}$.
    \item Each fibre $S^1_{\alpha}$ has a closed neighborhood $N_{\alpha}$ which which is homeomorphic to the solid torus $S^1 \times D^2$ and is a union of fibres.
    \item There is a covering map $p_{\alpha}:S^1 \times D^2 \rightarrow N_{\alpha}$ such that $p^{-1}_{\alpha}(S^1_{\alpha})$ is connected and
        \begin{enumerate}
            \item For each $x\in D^2$, $p$ maps each $S^1 \times x$ onto some $S^1_{\alpha^{\prime}}\subset N_{\alpha}$.
            \item $(S^1\times D^2, p)$ is a $\mathbb{Z}_m$-regular cover of $N_{\alpha}$ for some $m \in \mathbb{N}$.
        \end{enumerate}
\end{enumerate}

The space $S^1\times D^2$, fibred by the collection $\{S^1\times x:x\in D^2\}$, is the \emph{model} Seifert fibred space.  Also associated to the Seifert fibred 3-manifold $X$ is the \emph{base space} $B$.  $B$ is the quotient space of $X$ obtained by collapsing each fibre $S^1_{\alpha}$ to a point.  Conditions 2, 3(a), and 3(b) are the \emph{regularity conditions} of the fibration. A space which satisfies only the first condition in the definition is said to be \emph{foliated} by circles.  In this case, the circles are called \emph{leaves} of the foliation, not fibres.

We typically write the Seifert fibred 3-manifold $X$ and its Seifert fibration with the base space as
$$
S^1 \longrightarrow X \overset{\eta}{\longrightarrow} B
$$

\noindent where $B$ is the base space and $\eta$ is the quotient map.  For a Seifert fibred 3-manifold, $B$ is a 2-manifold.

Not every Seifert fibred 3-manifold non-trivially self covers, but some important subclasses of them do.  The most familiar of these are the torus knot complements and, of course, the products $S^1 \times F$, $F$ a 2-dimensional surface, compact or not, with or without boundary.

Seifert fibred phenomena occur in other dimensions as well, provided one changes the space used as the model fibred spaces and adapts the regularity conditions in an appropriate manner.  For example, it is clear that the products $X = S^1 \times Y$, for $Y$ any $n$ dimensional manifold produce $(n+1)$-dimensional manifolds which are Seifert fibred (the model fibred space space here is $S^1\times B^n$ where $B^n$ is the $n$-ball) and these $(n+1)$ dimensional spaces also non-trivially self cover.  And, in dimension 2, the annulus $S^1\times [0,1]$, the 2-torus $S^1\times S^1$, the Mobius band, and the Klein bottle are Seifert fibred spaces which non-trivially self cover.  Here the model fibred space is the annulus $S^1\times [0,1]$.

In this context, the Klein bottle is an interesting example.  It has two inequivalent Seifert fibrings\footnote{Two Seifert fibrations of a Seifert fibred space $X$ are equivalent if there is a Seifert fibration preserving homeomorphism $h:X \rightarrow X$ of $X$ with itself.} each of which have Seifert fibred respecting non-trivial self covers.  One of the Seifert fibrations of the Klein bottle is obtained by thinking of it as a pair of Mobius bands, each fibred by its center circle and circles which wrap around the center circle twice, glued together along their boundaries.  The second is that induced by first fibring the product $S^1\times [0,1]$ by the circles $S^1 \times t$, then gluing $S^1\times 0$ to $S^1\times 1$ via the antipodal map.  The first Seifert fibred Klein bottle has odd degree Seifert fibre respecting self covers for every odd integer and the second has degree $2^k$ Seifert fibre respecting self covers for every $k\in \mathbb{N}$.

The annulus, 2-torus, Mobius band, and Klein bottle are members of the class of generalized Baumslag-Solitar complexes, all of which are Seifert fibred 2-complexes.  However, for these complexes there is more than one model fibered space which can cover neighborhoods of the fibres in the complexes.  They are the spaces $S^1 \times K_n$ where $K_n =\{\rho \cdot \textrm{exp}(i\cdot \frac{2k\pi}{n}): 0\leq \rho \leq 1, k = 0,1,\ldots n-1 \}$, is the \emph{$n$-od}, the union of $n$ copies of $[0,1]$ joined at their origins.  Allowing this sort of variation for the regularity condition(s) occurs in the 3-manifold case as well.  See Gonz\'alez-Acu\~na et al \cite{GonzalezAcunaEtAl1994} where the model fibred spaces are either the solid torus $S^1\times D^2$ fibred by $S^1\times x$, as above, or fibred solid Klein bottles.  We look at the self covering properties of the generalized Baumslag-Solitar complexes in the next section.

Finally, observe that some spaces which are foliated by circles, but which are not Seifert fibred spaces, also non-trivially self cover.  For example, consider Cantor's Pearl Necklace and glue the entire circle $[0,2]/(0=2)$ back in.  Each $\mathbb{Z}_m$-regular self cover of Cantor's Pearl Necklace mentioned above clearly extents to the resulting space $X=\mathcal{N}^{(2)}_{C,S^2} \cup [0,2]/(0=2)$.   Foliate $X$ as follows.  From each $S^2$ in the necklace form an open annulus by deleting the rational points $\{\frac{i}{3^j},\frac{i+1}{3^j}$ in the Cantor set $ C \subset [0,1] \subset [0,2]/(0=2)$ from the sphere in the necklace which contains them. Then foliate the open annulus by a collection of pairwise disjoint embedded circles in the open annulus, each of which generates the annulus' fundamental group.  Then, include the circle $\mathcal{C} \cup [0,2]/(2=0)$ as a leaf in the foliation.  Observe that each leaf in this foliation of $\mathcal{N}^{(2)}_{C,S^2} \cup [0,2]/(0=2)$ has a closed neighborhood which is a union of leaves of the foliation and that each such neighborhood has a $\mathbb{Z}_m$-regular self cover for any $m$ of your choosing.  However the structure of the foliated neighborhoods containing the leaves in the $S^2 \setminus \{\frac{i}{3^j},\frac{i+1}{3^j}\}$'s versus the foliated neighborhoods of the leaf $[0,2]/(0=2)$ so very different that they lack any sort of unifying structure.

\section{Graphs of spaces}\label{sec:GraphsOfSpaces}
Let $\Gamma$ be a directed graph with vertex set $V(\Gamma)$ and directed edge set $E(\Gamma)$, loops and multiple edges between pairs of edges allowed.  A directed edge $e \in E(\Gamma)$ is denoted by $e = [e^-,e^+]=[u,v]$ where $e^- = u$ is the \emph{initial} vertex and $e^+ = v$ is the \emph{terminal} vertex of $e$.  A \emph{graph of topological spaces associated to $\Gamma$} is a topological space $X(\Gamma)$ obtained as follows:
\begin{itemize}
  \item For each $v\in V(\Gamma)$ there is a space $X_v$.
  \item For each directed edge $e=[e^-,e^+]\in E(\Gamma)$, there is a space $X_e$, a pair of disjoint subset $S_{e^{\pm}}$, and a pair of maps $f_{e^{\pm}}:S_{e^{\pm}}\rightarrow X_{e^{\pm}}$.
  \item $X(\Gamma)$ is formed by using $\mathcal{F} = \{f_{e^{\pm}}:e\in E(\Gamma)\}$ as gluing maps.
\end{itemize}

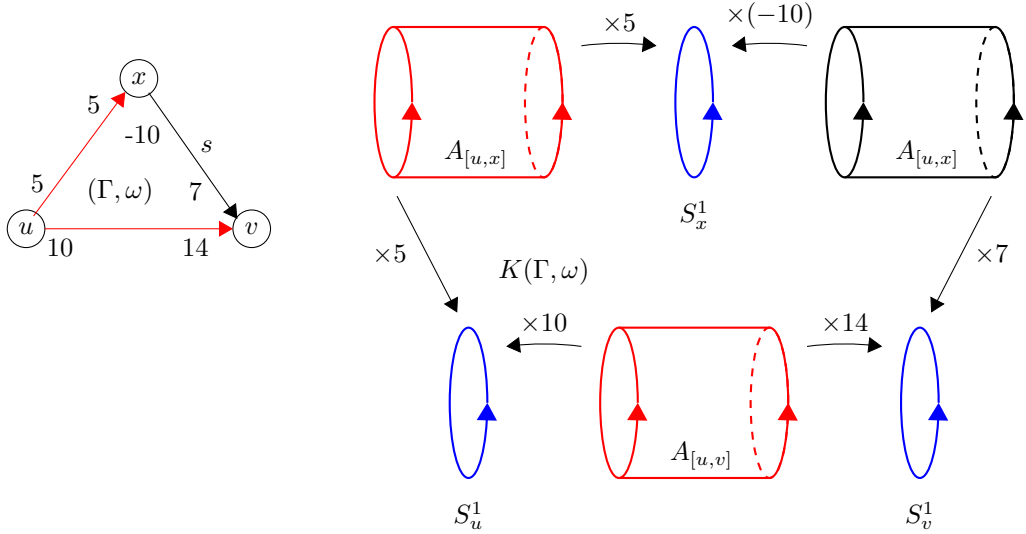
\begin{figure}
\centering
\begin{tikzpicture}
    \draw (0,0) circle (.25cm) (3,0) circle (.25cm) (1.5,2) circle (.25); 
    \node (u) at (0,0) {$u$};
    \node (v) at (3,0) {$v$};
    \node (x) at (1.5,2) {$x$};
    \begin{scope}[every path/.style={-triangle 60}]
       \draw[red] (.25,0) -- (2.75,0);
       \draw[red] (.1,.17) -- (1.32,1.82);
       \draw (1.63,1.81) -- (2.81,.13);
    \end{scope}
       \node (s) at (2.4,1.1) {$s$};
       \node (m00w1) at (.15,.6) {5};
       \node (m00w2) at (.45,-.25) {10};
       \node (m30w1) at (2.25,-.25) {14};
       \node (m30w3) at (2.25,.5) {7};
       \node (m1.5,2w1) at (.9,1.65) {5};
       \node (m1.5,2w2) at (1.55,1.25) {-10};
  \node (name) at (1.25,.5) {$(\Gamma,\omega)$};
  \node (a) at (0,-4) { };
\end{tikzpicture}
\hspace{1cm}
\begin{tikzpicture}
  \begin{scope}[every path/.style={-triangle 60}]
    \draw[red, thick] (0,0) arc(0:360:.25cm and 1cm);
    \draw[red, dashed, thick] (2,0) arc(0:360:.25cm and 1cm);
    \draw[blue,thick] (-2,0) arc(0:360:.25cm and 1cm);
    \draw[blue,thick] (4,0) arc(0:360:.25cm and 1cm);
    \node (S1u) at (-2.25,-1.5) {$S^1_u$};
    \node (S1v) at ( 3.75,-1.5) {$S^1_v$};
    \draw (-.75,.75) arc(60:120:1cm and .25cm);
    \draw (2.25,.75) arc(120:60:1cm and .25cm);
  \end{scope}
    \draw[red, thick] (1.75,-1) arc(-90:90:.25cm and 1cm);
    \draw[red,thick] (-.25,-1) -- (1.75,-1);
    \draw[red,thick] (-.25,1) -- (1.75,1);
    \node (Auv) at (.85,-.7) {$A_{[u,v]}$};
    \node (p0) at (-1.25,1.05) {$\times 10$};
    \node (p1) at (2.75,1.05) {$\times 14$};
  \begin{scope}[every path/.style={-triangle 60}]
    \draw[red, thick] (-3,4) arc(0:360:.25cm and 1cm);
    \draw[red, dashed, thick] (-1,4) arc(0:360:.25cm and 1cm);
    \draw[blue,thick] (1,4) arc(0:360:.25cm and 1cm);
    \node (S1x) at ( .75,2.5) {$S^1_x$};
    \draw (-.75,4.75) arc(120:60:1cm and .25cm);
    \draw (-3.2,2.75) -- (-2.4,1.2);
  \end{scope}
    \draw[red, thick] (-1.25,3) arc(-90:90:.25cm and 1cm);
    \draw[red,thick] (-3.25,3) -- (-1.25,3);
    \draw[red,thick] (-3.25,5) -- (-1.25,5);
    \node (Aux) at (-2.15,3.3) {$A_{[u,x]}$};
    \node (p0) at (-3.3,1.975) {$\times 5$};
    \node (p1) at (-.25,5.05) {$\times 5$};

  \begin{scope}[every path/.style={-triangle 60}]
    \draw[thick] (3,4) arc(0:360:.25cm and 1cm);
    \draw[dashed, thick] (5,4) arc(0:360:.25cm and 1cm);
    \draw (2.25,4.75) arc(60:120:1cm and .25cm);
    \draw (4.7,2.75) -- (3.9,1.2);
  \end{scope}
    \draw[thick] (4.75,3) arc(-90:90:.25cm and 1cm);
    \draw[thick] (2.75,3) -- (4.75,3);
    \draw[thick] (2.75,5) -- (4.75,5);
    \node (Aux) at (3.85,3.3) {$A_{[u,x]}$};
    \node (p0) at (4.7,1.975) {$\times 7$};
    \node (p1) at (1.75,5.15) {$\times (-10)$};
\node (2complex) at (-1.25,1.75) {$K(\Gamma,\omega)$};
\end{tikzpicture}
\caption{A GBS graph $(\Gamma,\omega)$ and its associated $K(\Gamma,\omega)$. The maximal subtree in $\Gamma$ is $T=[u,v]\cup[u,x]$.}\label{fig:GBSgraphAndComplex}
\end{figure}

In what follows, more structure is imposed on the graphs.  Those of primary interest are weighted, or labeled, directed graphs.  Let $\mathbb{Z}^*=\mathbb{Z}\setminus 0$.  Associate to the directed graph $\Gamma$ a \emph{weight function} $\omega : E(\Gamma) \rightarrow \mathbb{Z}^*\times\mathbb{Z}^*$ given by $\omega(e)= \omega(e^-,e^+)= (\omega^-(e),\omega^+(e))$.  The weighted directed graph is denoted $(\Gamma, \omega)$.  It is convenient to have a notation for weighted directed edges $e=[u,v], \ \omega(e)=(m,n)$.  In this notation, a weighted directed loop is denoted $(e,\omega)=[u^m,u^n]$ and a weighted directed straight edge is denoted $(e,\omega)=[u^m,v^n]$.  Pictures of these simple weighted directed graphs are in Figure \ref{fig:SingleEdge}.  A slightly more complicated one is in Figure \ref{fig:GBSgraphAndComplex}.  We call these graphs \emph{GBS graphs} for reasons explained below.

For the remainder of this section we let $S^1=\{e^{i\cdot \theta}\in \mathbb{C}:0\leq \theta \leq 2\pi\}$ be an oriented circle with the orientation given by the indicated parametrization.  Build a 2-dimensional graph complex, $K(\Gamma,\omega)$ as follows.  For each $v\in V(\Gamma)$, let $X_v=S^1_v$ be copy of the oriented circle $S^1$.  For each directed edge $e\in E(\Gamma)$, let $A_e=I \times S^1$, $I=[0,1]$ with the orientations on $\partial A_e = \{0,1\}\times S^1_e$ given by the above parametrization of $S^1$.  Define gluing maps $f_{e^{\pm}}$ on $\partial A_e$ as follows.

\begin{itemize}
  \item $f_{e^-}:\partial_{e^-} A=0\times S^1_e \rightarrow S^1_{e^-}$ is $f(0,e^{i\cdot \theta})=\textrm{exp}(i\cdot \omega^-(e) \theta)\in S^1_{e^-}$.
  \item $f_{e^+}:\partial_{e^+} A=0\times S^1_e \rightarrow S^1_{e^+}$ is $f(1,e^{i\cdot \theta})=\textrm{exp}(i\cdot \omega^+(e) \theta)\in S^1_{e^+}$.
\end{itemize}

\noindent The resulting space is denoted $K(\Gamma,\omega)$.  A picture of a simple GBS complex is in Figure \ref{fig:GBSgraphAndComplex}.

The above construction is the topological analog of the graph of groups approach which produces the \emph{generalized Baumslag-Solitar (GBS)} groups.  These groups are obtained as follows.  Let $(\Gamma,\omega)$ be a directed weighted graph.  Associate to each $v\in V(\Gamma)$ a copy $(\mathbb{Z}_v,+)$ of the integers.  Associate to each directed edge $e\in E(\Gamma)$ a copy $(\mathbb{Z}_e,+)$ of the integers and a pair of monomorphisms $f_{e^{\pm \ }}:\mathbb{Z}_e\rightarrow \mathbb{Z}_{e^{\pm}}$ determined by mapping $1\in \mathbb{Z}_e$ into, respectively, the vertex groups $\mathbb{Z}_{e^{\pm}}$ via $1\overset{f_{e^{\pm}}}{\longmapsto}\omega^{\pm}(1)\cdot 1$.  This produces the group of the weighted graph $(\Gamma,\omega)$.  This group is denoted $\pi(\Gamma,\omega)$.  There is also a description of $\pi(\Gamma,\omega)$ in terms of generators and relations which reflects the structure of $(\Gamma,\omega)$.  It is obtained as follows.
\begin{itemize}
    \item Choose a spanning subtree $T$ of $\Gamma$.
    \item For each vertex $v \in V(\Gamma)$ there is generator $g_v$ and for each non-$T$ edge in $E(\Gamma)$ there is a generator $t_e$ in the group of the graph $\pi(\Gamma,\omega)$.
    \item For each directed edge $e=[u,v]$ in $T$, there is a relation of the form $g^{\omega^-(e)}_u = g^{\omega^+(e)}_v$.
    \item For each directed non-$T$ edge $e=[u,v]$ in $\Gamma$ there is a relation of the form $t^{-1}_e g^{\omega^-(e)}_u t_e = g^{\omega^+(e)}_v$.  Note that it is possible that $u=v$ here, i.e., $[u,v]$ is the loop $[u,u]$.
\end{itemize}
Note that $\pi(\Gamma,\omega)$ does not depend on the choice of the spanning tree $T$ and also that it is not the topological fundamental group of the graph $\Gamma$.  The topological fundamental group of a graph is always a free group and, with the exceptions of the cases where $\Gamma$ is a single vertex or a single straight edge with at least one weight $\pm 1$, $\pi(\Gamma,\omega)$ is not a free group.

Given a maximal subtree $T$ of $\Gamma$, there is a canonical selection of generators for the fundamental group $\pi_1(K(\Gamma,\omega),x_0)$ and a straightforward induction on the number of edges in $\Gamma$, using Van Kampen's Theorem, which show that $\pi_1(K(\Gamma,\omega),x_0) \cong \pi(\Gamma,\omega)$.  See  \cite{BedDelTimm2006}, \cite{Campbell_Quick_Robertson_Roney-Dougal_Stewart_2024}.  In recognition of this connection between the graph $(\Gamma,\omega)$, the group $\pi(\Gamma,\omega)$, and the space $K(\Gamma,\omega)$, call the weighted graphs \emph{generalized Baumslag-Solitar (GBS)} graphs and the spaces \emph{generalized Baumslag-Solitar (GBS)} complexes.

The selection of the canonical system of generators for $\pi_1(K(\Gamma,\omega),x_0)$ depends on the observation that $K(\Gamma,\omega)$ is a Seifert fibred 2-complex with base space a copy of $\Gamma$ and that the Seifert fibration has a section
$$
\Gamma \overset{i}{\longrightarrow} K(\Gamma,\omega) \overset{\eta}{\longrightarrow}\Gamma.
$$

\noindent The fibres in the Seifert fibration are the vertex circles $S^1_v$, $v\in V(\Gamma)$, and the circles $t\times S^1$, $t\in (0,1)_e$ in the interior of the annulus $A_e = [0,1]\times S^1$, associated to each edge $e\in E(\Gamma)$.

\bigskip
Let $(\Gamma, T)$ denote the embedded copy of the graph $\Gamma$ and a spanning tree $T$ in the GBS complex $K(\Gamma, \omega)$.  Fix the base point $x_0=v_0$ in some vertex circle $S^1_{v_0}\cap T$.  Then for every vertex $v$ of $\Gamma$, the associated vertex generator is the homotopy class of a loop based at $x_0=v_0$ which runs along the embedded copy of $T \subset \Gamma \subset K(\Gamma, \omega)$, then around the vertex circle $S^1_v$, and then returns to $x_0$ along $T$.  For each oriented non-$T$ edge $e=[u,v]$, the associated edge generator $t_e$ is the homotopy class of the directed loop based at $x_0$ which starts at $x_0$, runs along $T$ to $u$, crosses $[u,v]$, then follows $T$ from $v$ back to $x_0$.  Simplify the notation by suppressing base points in the notation and using the vertices and the non-$T$ edges themselves to represent the homotopy classes of their associated loops.

For example, for the two GBS graphs in Figure \ref{fig:SingleEdge}, choose $T=u$ and $T=[u,v]$, respectively.  This produces

\begin{itemize}
  \item $\pi_1(K([u^m,u^n]))\cong \pi([u^m,u^n]) \cong \langle u,t : t^{-1}u^mt=u^n \rangle$ and
  \item $\pi_1(K([u^m,v^n]))\cong \pi([u^m,u^n]) \cong \langle u,v : u^m=u^n \rangle.$
\end{itemize}

\noindent For the graph and maximal subtree selected in Figure \ref{fig:GBSgraphAndComplex}, the fundamental group is
$$
\pi_1(K(\Gamma,\omega))\cong \pi(\Gamma,\omega)\cong \langle u,v,x,s:u^5=x^5,u^{10}=v^{14},s^{-1}x^{-10}s=v^7\rangle.
$$

\noindent The groups for the graphs which have a single loop are the classical Baumslag-Solitar groups.  In the case of the graphs with a single edge, when $m$ and $n$ are relatively prime, these are the graphs of the $(m,n)$-torus knot groups, otherwise, they are fundamental groups of other simple Seifert fibred 3-manifolds.

\begin{figure}
  \centering
\begin{tikzpicture}
    \draw (0,0) -- (7,0);
    \draw [dashed] (7,0) -- (9,0);
    \draw[black,fill=black] (0,0) circle (.15cm) (2,0) circle (.15cm) (4,0) circle (.15) (6,0) circle (.15);
    \node (m1) at (.4,.3) {$m_1$};
    \node (n1) at (1.6,.3) {1};
    \node (m2) at (2.4,.3) {$m_2$};
    \node (n2) at (3.6,.3) {1};
    \node (m3) at (4.4,.3) {$m_3$};
    \node (n3) at (5.6,.3) {1};
    \node (m4) at (6.4,.3) {$m_4$};
    \node (a) at (-1,0) {\textrm{(a)}};
    \node (Km) at (9.75,0) { \ :$K(\Gamma_{\textrm{\textbf{m}}})$};
    \node (spacer) at (0,-1) {};
\end{tikzpicture}
\begin{tikzpicture} 
    \node (C1) at (0,-2) {$S^1$};
    \node (C2) at (2,-2) {$S^1$};
    \node (C3) at (4,-2) {$S^1$};
    \node (C4) at (6,-2) {$S^1$};
    \begin{scope}[every path/.style={-triangle 60}]
       \draw (1.6,-2) -- (.4,-2);
       \draw (3.6,-2) -- (2.4,-2);
       \draw (5.6,-2) -- (4.4,-2);
       \draw (7.2,-2) -- (6.4,-2);
    \end{scope}
       \draw [dashed] (7.2,-2) -- (9,-2);
    \node (m1) at (1.1,-1.75) {$\times m_1$};
    \node (m2) at (3.1, -1.75) {$\times m_2$};
    \node (m3) at (5.1, -1.75) {$\times m_3$};
    \node (m4) at (7.1,-1.75) {$\times m_4$};
    \node (b) at (-1,-2) {\textrm{(b)}};
    \node (Sigmam) at (9.5,-2) {:$\Sigma_{\textrm{\textbf{m}}}$};
\end{tikzpicture}
\caption{Let $\textrm{\textbf{m}}=(m_1,m_2,m_3,\ldots)$.  (a) The graph for $K(\Gamma_{\textrm{\textbf{m}}})$.  All edges in the picture are oriented from right to left. (b) A picture of the solenoid $\Sigma_{\textrm{\textbf{m}}}$.} \label{fig:GBS+Inv}
\end{figure}
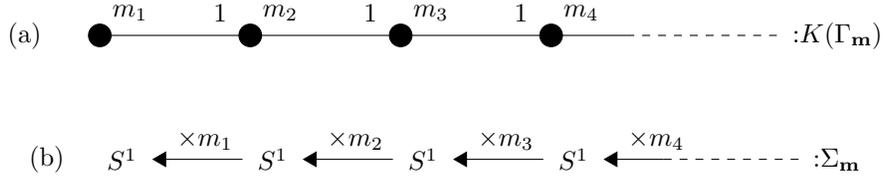

The standard approach to developing the group theory of $\pi(\Gamma,\omega)$ is to exploit structural aspects of the weighted graphs and their combinatorics, e.g., \cite{DelRobTimm2011}, Levitt \cite{Levitt2007}, Robinson \cite{Robinson2011}.  However, the connections between the combinatorics of $(\Gamma, \omega)$ and the topology of $K(\Gamma,\omega)$ also allow for topological methods to provide insight into the group theory \cite{Timm2022Interact}, \cite{Campbell_Quick_Robertson_Roney-Dougal_Stewart_2024}.  In particular, in the context of this paper, by \cite{BedDelTimm2006}, the GBS complex $K(\Gamma,\omega)$ has a nontrivial degree $d$ self cover for every $d$ which is relatively prime to all the weights on $(\Gamma,\omega)$.  Consequently, when $\Gamma$ is a finite graph, $K(\Gamma,\omega)$ has non-trivial self covers.  Therefore, it follows that every GBS group is non-cohopfian.  Also note that the presentation complex for $\pi(\Gamma,\omega)$ is obtained from the GBS complex $K(\Gamma,\omega)$, by choosing a spanning tree $T \subset \Gamma$, embedding the pair $(\Gamma, T)$ in $K(\Gamma,\omega)$ via the section $i: \Gamma \rightarrow K(\Gamma,\omega)$ for the Seifert fibration, then contracting the embedded copy $i(T) \subset \Gamma$ to a point.

\bigskip
It is interesting to compare the diagrams we use to denote the solenoids to certain infinite GBS graphs. Fix a sequence $\textrm{\textbf{m}}=(m_1,m_2,m_3,\ldots)$.  Form the corresponding GBS graph $\Gamma_{\textrm{\textbf{m}}}=\cup\{[v^{m_i}_i,v^1_{i+1}]:i\in \mathbb{N}\}$ and GBS complex $K(\Gamma_{\textrm{\textbf{m}}})$.  Also form the solenoid $\Sigma_{\textrm{\textbf{m}}}$.  The visual similarities of the pictures we draw to represent these spaces are striking, see Figure \ref{fig:GBS+Inv}, but there are also more fundamental topological connections as well.  In each picture, edges represent instructions for wrapping the $n+1^{\textrm{st}}$ copy of $S^1$ around the $n^{\textrm{th}}$.  Consequently, the GBS complex $K_{\textrm{\textbf{m}}}$ can be thought of as 2-dimensional representation of the inverse system which generates the solenoid $\Sigma_{\textrm{\textbf{m}}}$.

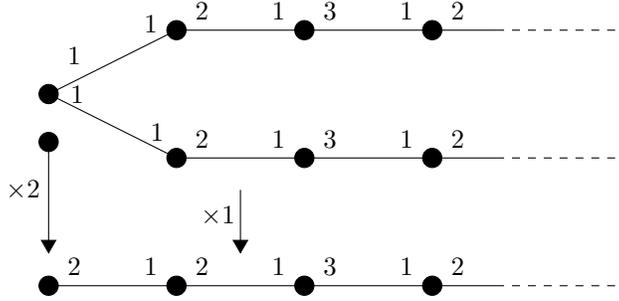
\begin{figure}[t]
  \centering
\begin{tikzpicture}[scale=.85]
    \draw (0,0) -- (7,0);
    \draw [dashed] (7,0) -- (9,0);
    \draw[black,fill=black] (0,0) circle (.15cm) (2,0) circle (.15cm) (4,0) circle (.15) (6,0) circle (.15);
    \node (m1) at (.4,.3) {$2$};
    \node (n1) at (1.6,.3) {1};
    \node (m2) at (2.4,.3) {$2$};
    \node (n2) at (3.6,.3) {1};
    \node (m3) at (4.4,.3) {$3$};
    \node (n3) at (5.6,.3) {1};
    \node (m4) at (6.4,.3) {$2$};
    \node (a) at (4.5, -.75) {  };
    \draw (0,3) -- (2,4); \draw (0,3) -- (2,2);
    \draw (2,4) -- (7,4); \draw[dashed] (7,4) -- (9,4);
    \draw (2,2) -- (7,2); \draw[dashed] (7,2) -- (9,2);
    \draw[black,fill=black] (0,3) circle (.15cm) (2,4) circle (.15cm) (4,4) circle (.15) (6,4) circle (.15);
    \draw[black,fill=black] (2,2) circle (.15cm) (4,2) circle (.15) (6,2) circle (.15);
    \node (m11) at (.45,3) {$1$}; \node (n11) at (1.7,2.4) {1};
    \node (m12) at (2.4,2.3) {$2$}; \node (n12) at (3.6,2.3) {1};
    \node (m13) at (4.4,2.3) {$3$}; \node (n13) at (5.6,2.3) {1};
    \node (m14) at (6.4,2.3) {$2$};

    \node (m11) at (.4,3.6) {$1$}; \node (n11) at (1.6,4.1) {1};
    \node (m12) at (2.4,4.3) {$2$}; \node (n12) at (3.6,4.3) {1};
    \node (m13) at (4.4,4.3) {$3$}; \node (n13) at (5.6,4.3) {1};
    \node (m14) at (6.4,4.3) {$2$};
    \draw[black,fill=black] (0,2.25) circle (.15cm);
    \begin{scope}[every path/.style={-triangle 60}]
        \draw (0,2.25) -- (0,.5); \draw (3,1.5) -- (3,.5);
    \end{scope}
    \node (deg) at (-.4,1.5) {$\times 2$}; \node (deg2) at (2.65,1.1) {$\times 1$};
\end{tikzpicture}
\caption{The degree 2 cover of $K_{\textrm{\textbf{p}}}$. All edges are again oriented from right to left.}\label{Deg2CovBlockPrime}
\end{figure}

The defining sequence $\textrm{\textbf{m}}$ determines the covering space properties of both $K_{\textrm{\textbf{m}}}$ and $\Sigma_{\textrm{\textbf{m}}}$.  For example, setting $m_i = 2$ for all $i\in \mathbb{N}$, both the dyadic GBS complex $K_{\textrm{\textbf{2}}}$ and the dyadic solenoid $\Sigma_{\textrm{\textbf{2}}}$ have odd degree self covers for every odd natural number.  At the other extreme, setting
$$
\textrm{\textbf{p}}=(2,2,3,2,3,5,2,3,5,7,2,3,5,7,11,\ldots)
$$
\noindent where one adds blocks of consecutive primes $2,3,5,\ldots,p_j$ of length $j\in \mathbb{N}$, neither $K_{\textrm{\textbf{p}}}$ nor $\Sigma_{\textrm{\textbf{p}}}$ have any nontrivial self-covers at all.  However, both $K_{\textrm{\textbf{2}}}$ and $K_{\textrm{\textbf{p}}}$ have many covering spaces which are not self covers. For example, $K_{\textrm{\textbf{2}}}$ has $\mathbb{Z}_{2^n}$-regular non-self covers for all $n\in \mathbb{N}$ and $K_{\textrm{\textbf{p}}}$ has a $\mathbb{Z}_{p_i}$-regular non-self cover for each prime $p_i$, corresponding to the first time $p_i$ appears in the sequence $\textrm{\textbf{p}}$. A picture for the degree 2 cover of $K_{\textrm{\textbf{p}}}$ is in Figure \ref{Deg2CovBlockPrime}.  The GBS graphs for these covering spaces of $K_{\textrm{\textbf{2}}}$ and $K_{\textrm{\textbf{p}}}$ are all directed rooted trees similar to that in Figure \ref{Deg2CovBlockPrime}.  More details about the connections between the combinatorics of the GBS graphs and the topology of the associated GBS complex are in \cite{Campbell_Quick_Robertson_Roney-Dougal_Stewart_2024}. Observe that the GBS graph $\Gamma_{\textbf{\textrm{m}}}$ produces the presentation  $\pi(\Gamma_{\textrm{\textbf{m}}}) \cong \langle x_i:x_{i+1}=x^{m_i}_i, \ i\in\mathbb{N}\rangle$, an infinite presentation for $\mathbb{Z}$.

\bigskip
The graph theoretic based approach to the GBS complexes $K(\Gamma,\omega)$ exploits the facts that each edge space is the product of the vertex spaces $S^1$ with an interval and, of relevance to the self covering phenomena of interest in this paper, that $S^1$ has nontrivial self covers.  Thus, the construction can easily be generalized.

Fix a space $Z$ which has nontrivial self covers. Let $D = \{d: Z \textrm{ has a degree } d \textrm{ self cover}\}$.  Let $(\Gamma,\omega)$ be a graph with directed edge set.  Let $\omega: E(\Gamma) \rightarrow D\times D$ be a weight function with $\omega(e)=(\omega^-(e),\omega^+(e))$.  For each $v \in V(\Gamma)$ let $Z_v$ be a copy of $Z$ and for each directed edge $e \in E(\Gamma)$ let $Z_e = Z\times [0,1,]$. For each edge $e$, choose a self cover cover $f_{e^-}:Z\times 0 \rightarrow Z_{e^-}$ of degree $\omega^-(e)\in D$ and a self cover $f_{e^+}:Z\times 1 \rightarrow Z_{e^+}$ of degree $\omega^+(e)\in D$ as gluing maps.  Let $\mathcal{F}=\{f_{e^{\pm}}: e\in E(\Gamma)\}$.  The resulting space is denoted $K(\Gamma,\omega,Z,\mathcal{F})$ and is called a \emph{graph of} $Z$'s.

\begin{question}
Which $K(\Gamma,\omega,Z,\mathcal{F})$ have non-trivial self covers?
\end{question}

The obvious version of this question to consider is the generalization of the GBS complex construction to 3-complexes which arise as graphs of $Z$'s in which $Z=S^1\times S^1$.  The associated group theoretic questions is then whether a graph of groups in which all vertex groups and edge groups are copies of $\mathbb{Z}\times \mathbb{Z}$ and all edge maps are monomorphisms, is non-cohopfian.  For example, consider the case where the graph is a single edge $e=[u,v]$ and the gluing maps are of the form $f_{e^-}(0\times (e^{i\cdot \alpha},e^{i\cdot \beta}))=(e^{i\cdot m\alpha},e^{i\cdot \beta})$ and $f_{e^+}(1 \times (e^{i\cdot \alpha},e^{i\cdot \beta}))=(e^{i\cdot n\alpha},e^{i\cdot \beta})$.  The resulting space is the product $K([u^m,v^n])\times S^1$.  It has non-trivial self covers of degree $d$ for each $d$ relatively prime to both $m$ and $n$ coming from those of the $K([u^m,v^n])$ factor and others coming from the $S^1$ factor, but what happens when the gluing maps are more complicated self covers of $S^1\times S^1$ is an open question.  In addition, the results in \cite{EdaMatijevic2006} and the discussion of the similarities and differences between the solenoids $\Sigma_{\textrm{\textbf{m}}}$ and the infinite GBS complexes $K(\Gamma_{\textrm{\textbf{m}}})$ suggest that consideration of graphs of $S^1\times S^1$'s using infinite graphs is warranted.

\bigskip
We close this section by noting that the graph-based methods of this section can be combined with the Cantor set controlled constructions of Section \ref{sec:CantorSetControl} to produce more low (and high) dimensional spaces which non-trivially self cover.

\begin{example}\label{ex:GBS+CantorControlled} Start with a GBS complex $K(\Gamma,\omega)$ and choose a fiber $t_0\times S^1$ in the Seifert fibration in the interior of one of the annuli in $K(\Gamma,\omega)$.  On this fibre, choose an appropriately embedded middle thirds Cantor set $C$ as control set, then between each pair $\frac{i}{3^j},\frac{i+1}{3^j}$ of consecutive rational points used to form $C$, delete the interiors of a collection of disks $\mathcal{D}$ which have diameters that decrease to 0 in levels.  The resulting continuum $K(\Gamma,\omega) \setminus \left(\underset{D\in \mathcal{D}}{\cup}D\right)$ has nontrivial self covers of degree $d$ for each $d$ relatively prime to all the weights on $(\Gamma,\omega)$.  Since an annulus can be represented as the GBS complex $K([x^1,y^1])$, Cantor's Annulus in Figure \ref{fig:CountPuncAnnTangent} is an examples of this sort of construction where the Seifert fibre chosen as the home of the Cantor set is the middle circle $0 \times S^1$ of the annulus.  This sort of construction can be done for arbitrary collections $\{S^1_{\alpha}:\alpha \in I\}$ of Seifert fibres which can, in fact, contain vertex circles of the GBS complex and still produce spaces which non-trivially self cover.
\end{example}

\section{Complex Dynamical Systems}\label{sec:DynamSystems}
The connections between self covering phenomena and complex dynamics discussed in this section were brought to my attention by Jim Belk.  The spaces of interest in this context are the complementary subsets of the Riemann sphere $\widehat{\mathbb{C}}=\mathbb{C}\cup \infty$ known as the Fatou and Julia sets of self maps of the Rieman sphere $f:\widehat{\mathbb{C}} \rightarrow \widehat{\mathbb{C}}$.  For the specific examples mentioned here, $f$ is either a complex polynomial or rational function, but more general situations are considered in the literature.  The terminology used is that of Milnor \cite{Milnor2006}.  Note that what follows depends heavily on the Open Mapping Theorem for holomorphic functions on Riemann surfaces.

Suppose $f:\widehat{\mathbb{C}} \rightarrow \widehat{\mathbb{C}}$ is a function.  For $z \in \mathbb{C}$, let $f^1(z) = f(z)$ and, in general, define $f^n:\widehat{\mathbb{C}} \rightarrow \widehat{\mathbb{C}}$ by $f^{n+1}(z) = f(f^n(z))$, the $n+1^{\textrm{st}}$ \emph{forward iterate of $f$}.  Similarily, let $f^{-1}(z) = \{x\in \widehat{\mathbb{C}} : f(x) = z\}$ and define the $n+1^{\textrm{st}}$ \emph{backward iterate $f^{-n-1}(z)$} by $f^{-1}\left(f^{-n}(z)\right)$.  Given a sequence of maps $f_n: S \rightarrow \widehat{\mathcal{C}}$  on a
Riemannn surface $S$, it \emph{converges locally uniformly} to the limit $g : S \rightarrow \widehat{\mathcal{C}}$, if for every compact subset $K\subset S$,  the sequence $fn|K$ of $f_n$ restricted to $K$ converges uniformly to $g|K$. A collection $\mathcal{F}$ of homomorphic maps $f_{\alpha}:S \rightarrow \widehat{\mathcal{C}}$ is \emph{normal} if every infinite sequence of maps from $\mathcal{F}$ converges locally uniformly to a limit.  For a given map $f:\widehat{\mathcal{C}} \rightarrow \widehat{\mathcal{C}}$, the collection $\mathcal{F} = \{f^n:n\in \mathbb{N}\}$ of forward iterates of $f$ (and their backward iterates when of interest) is the \emph{iterated function system} for $f$.

\begin{definition}
Let $f:\widehat{\mathcal{C}} \rightarrow \widehat{\mathcal{C}}$ be a non-constant holomorphic map.  Let $\mathcal{F} = \{f^n:n\in \mathbb{N}\}$, the iterated function system for $f$.  The \emph{Fatou set of $f$} is the set of all points $z_0 \in \widehat{\mathcal{C}}$ such that there is a neighborhood $N(z_0)$ of $z_0$ on which $\mathcal{F}$ is a normal family.  The \emph{Julia set} $\mathcal{J}(f)$ is the set of points in $\widehat{\mathcal{C}}$ for which no such neighborhood exists.
\end{definition}

\noindent
That is, the Fatou set of a holomorphic map $f$ is the set of points $z_0$ for which there is a neighborhood $N(z_0)$ on which the dynamics of $f$ are ``tame'' or ``well-behaved.''  The Julia set, on the other hand is the set of points where the dynamics in neighborhoods containing it are ``wild,'' ``chaotic,'' or ``exhibit sensitive dependence on initial conditions."  By definition, it is clear that the Fatou and Julia sets are a pair of disjoint sets whose union is all of $\widehat{\mathbb{C}}$.  The Julia set is a closed subset of $\widehat{\mathbb{C}}$ and the Fatou set is the complementary open set.

\bigskip
Familiar examples illustrate the basic dynamics of iterated function systems, including their connections to spaces which non-trivially self cover.  Consider, for example, $f:\widehat{\mathcal{C}} \rightarrow \widehat{\mathcal{C}}$ given by $f(z)=z^2$ where $f(\infty)=\infty$ and its forward (and backward) iterates.  See Milnor \cite{Milnor2006}, page 41.

In this case the Fatou set is $\{z: |z|<1\} \cup \{z:|z|>1\}$ and we see that on this set set the dynamics of $f$ are easy to understand.  On $D = \{z:|z| <1\}$, the forward iterates $f^n(z), z \in D$ clearly converge pointwise to 0.  Moreover it is easy to see that the forward iterates $f^n$ of $f$ converge uniformly to the constant function $g(z)=0$ on $D$.  On $\widehat{\mathcal{C}} \setminus cl(D)$ the iterates $f^n$ converge to the constant function $g(z)=\infty$.

However, for $z \in S^1 = \partial(D)$, the Julia set for $f$, the dynamics on neighborhoods of $z$ are more complicated.  Each neighborhood of $z \in S^1$ contains points whose forward iterates under $f$ converging to 0, points with forward iterates under $f$ converging to $\infty$, and points for which all forward iterates of $f$ remain in $S^1$.  Thus, there is a discontinuity in the behavior of the dynamics along $S^1$ which causes the dynamics for $f$ to exhibit sensitive dependence on initial conditions on neighborhoods of points in $S^1$.

Observe that $f$ has three fixed points, $z=0,z=\infty,$ and $z=1$.  That is, $f$ has three periodic points of period 1. Two of these periodic points, $z=0$ and $z=\infty$, are \emph{critical points} of $f$, i.e., there is no neighborhood of either 0 or $\infty$ on which $f$ is a homeomorphism. In particular, note that this occurs at $0$ because $f^{\prime}(0) = 0$.  Both 0 and $\infty$ are in the Fatou set.  Both are attracting periodic points.  On the other hand, $z=1$, the periodic point contained in $J(f)$ has $f^{\prime}(1) = 2 \neq 0$ and is a repelling periodic point of $f$.  Of greatest relevance to our current discussion, observe that the restriction $f|J(f) = f|S^1$ is a degree 2 self-cover.  Additionally, note that for each $z\in J(f) = S^1$, the derivative of $f$ evaluated as $z$ is $f^{\prime}(z) =2z \neq 0$ and that there are neighborhoods $N(z)$ and $U(f(Z))$ in $\mathbb{C}$ such that $f$ maps $J(f)\cap N(z)$ homeomorphically onto $J(f) \cap U(f(z))$.

Two of the more general results illustrated by the above observations about the dynamics of $f(z)=z^2$ are the following.

\begin{proposition}
\emph{\textbf{Invariance Lemma}} If $f:S\rightarrow S$ is a holomorphic map of the Riemann surface S and $J(f)$ is its Julia set, then $z \in J(f)$ if and only if $f(z) \in J(f)$.  Consequently, the full set of forward and backward iterates of a point $z \in J(f)$ is contained in $J(f)$.  \emph{\cite{Milnor2006}, Lemma 4.3}.
\end{proposition}

\begin{proposition}
Suppose $f:\widehat{\mathbb{C}} \rightarrow \widehat{\mathbb{C}}$ is a nonconstant holomorphic function with Julia set $J(f)$.  If $z\in J(f)$ and $f^{\prime}(z) \neq 0$, then are neighborhoods $N(z)$ and $U(f(z))$ such that the restriction of $f$ to $N(z)\cap J(f)$ is a conformal isomorphism of $N(z)\cap J(f)$ onto $U(f(z))\cap J(f)$. \emph{\cite{Milnor2006}, page 44.}
\end{proposition}

For our purposes it is worth taking a moment to compare $f(z)=z^2$ to two other functions $g(z) = z^2 - 2$ and $h(z)=z^2-6$.  By \cite{Milnor2006}, Lemma 7.1, the Julia set for $g$ is the closed interval $J(g)=[-2,2]$ on the $x$-axis in $\mathbb{C}$.  Critical points of $g$ are $0,\infty$ and its fixed points are $-1,2,\infty$.  In this case, even though $g(J(g)) = J(g)$, it is not the case that the restriction $g|J(g)$ to the Julia set is a self cover.  This occurs because the critical point $z=0$ for $g$ is in $J(g)$ and as a consequence, $g$ maps small open disks $N$ around $z=0$ to small open sets $U$ about $g(0)$, via a branched cover.  For $h$, critical points are again $0$ and $\infty$.  The fixed points are $z=3$ and $z=-2$.  By \cite{Milnor2006}, Problem 4-e, the Julia set is a Cantor set contained in $[-3,-\sqrt{3}]\cup[\sqrt{3},3]$.  Note the critical point $0\notin J(h)$ and that $g|J(h)$ is a double cover.

These three simple examples suggest some of the questions addressed in the literature about the Julia sets of holomorphic functions:  (1) For which holomorphic self maps $f$ on the Riemann sphere is the Julia set $J(f)$ smooth?  For which is $J(f)$ connected?  For which is $J(f)$ locally connected? For which is $J(f)$ an interval? For which is $f$ a self cover of $J(f)$?

\begin{figure}
\centering
    \includegraphics[width=.3\textwidth]{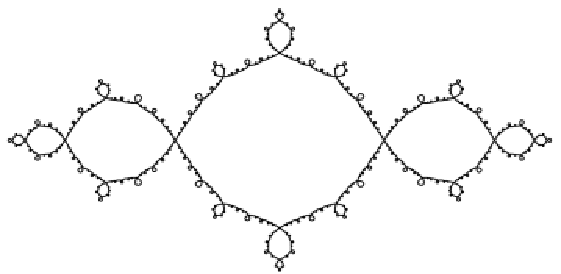}
    \includegraphics[width=.3\textwidth]{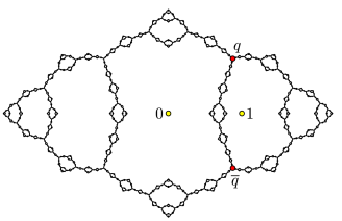}
    \includegraphics[width=.3\textwidth]{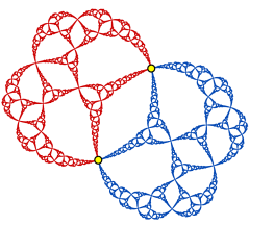}
\begin{tikzpicture}
\node (a) at (-5.05,0) {\textrm{(a)}};\node (b) at (0,0) {\textrm{(b)}};\node (c) at (5.05,0) {\textrm{(c)}};
\end{tikzpicture}
\caption{(a) The bascilica $J(z^2 - 1)$. (b) The keep $J(z^{-2} -1)$. (c) Orchids $J(\frac{w \cdot z^2 - 1}{z^2 - 1})$, $w=\exp(\frac{2\pi \dot i}{3})$.  Thanks to Jim Belk \cite{BelkPC} for these examples and the pictures.}\label{fig:juliasets}
\end{figure}

Examples of more complicated Julia sets which double cover are the basilica, the Julia set of $f(z) = z^2-1$; the ``keep,'' that of $g(z) = z^{-2} - 1$, and ``orchids,'' the Julia set of $h(z) = \frac{w\cdot z^2 - 1}{z^2 - 1}$, $w=\exp(\frac{2\pi \dot i}{3})$.  See Figure \ref{fig:juliasets}.

The dynamics of maps on the Riemann sphere and their Julia sets are closely related to the theory of self-similar groups, that is, groups generated by finite automata.  Like GBS groups, these groups are easy-to-describe groups with interesting properties.  The Grigorchuk group \cite{Grigorchuk1980}, an infinite finitely generated torsion group (with other interesting properties), and the iterated monodromy group of the basilica $f(z) = z^2 -1$, see  Bartoldi and Nekrashevych \cite{BartholdiNekrashevych2008} are examples of self-similar groups.  For a the type of of rational functions $f$ on the Riemann sphere called \emph{post critically finite rational functions}, the Julia set $J(f)$ can be reconstructed from the iterated monodromy group associated to $f$.  The basilica $f(z)=z^2 - 1$ is one such example.

We very briefly describe a bit of the group theory and topology for the special case of the iterated monodromy groups.  The preface to Nekrashevych \cite{Nekrashevych2005} contains longer expository introduction to the group theory and topology of the full class of self similar groups.

Suppose that $f:M_1 \rightarrow M$ is covering space of the manifold by $M$ by an open subset $M_1 \subset M$.  Let $M_1 = f^{-1}(M)$ and in general let $M_n=f^{-1}(M_{n-1})$  This generates an inverse system
$$
M \overset{f}{\longleftarrow} M_1 \overset{f}{\longleftarrow} M_2 \overset{f}{\longleftarrow} \cdots \overset{f}{\longleftarrow} M_n \overset{f}{\longleftarrow} M_{n+1} \overset{f}{\longleftarrow} \cdots \cdots $$

Choose a $t \in M$ and consider the set $T = \overset{\infty}{\underset{n=1}{\cup}} f^{-n}(t)$, the full set of inverse iterates of $t$ under $f$.  $T$ is naturally an $n$-ary rooted tree with root $t$.  The vertices of $T$ are $t$ and all its inverse iterates $f^{-n}(t)$, $n\in \mathbb{N}$.  The edges of $T$ are obtained by connecting a vertex $y$ of $T$ to each of its $n$ inverse images $f^{-1}(y) = \{x_1,x_2,\ldots x_n\}$.  There is an natural action of the pointed fundamental group $\pi_1(M,t)$ on $T$ acting by automorphisms of $T$.  Let $K_n$ be the kernel of this action on the $n^{\textrm{th}}$ level $f^{-n}(t)$.  The \emph{iterated monodromy group of $f$} is the quotient
$$
IMG(f) = \pi_1(M,t)/ \overset{\infty}{\underset{n=1}{\cap}}K_n.
$$
\noindent  It is from the action of $IMG(f)$ on $T$ that one can recover the Julia set of a post critically finite rational map $f$. See \cite{Nekrashevych2005}, Chapter 5.

David Bellamy \cite{Bellamy2005} uses properties of holomorphic maps to produce his self covers of the pseudo circle.  One wonders if there are deeper connections between his result and the dynamical processes discussed in this section.  In particular, is the pseudo circle the Julia set of a holomorphic self map $f$ of the Riemann sphere?

\section{Geometric Representability}\label{sec:GeoRep}
This section contains brief discussions of several topics which are peripherally related to self covering.  Their common theme is the geometric representability of algebraic properties via the fundamental group or other groups associated to topological spaces of particular dimensions or types.

\bigskip
For relatively prime pairs $p,q \geq 2$, the fundamental groups of the exteriors of the $(p,q)$-torus knots are the GBS groups $\pi([x^p,y^q])$.  This is a collection of finitely generated non-abelian groups with finite index clones and the self covering of the knot exteriors show that the non-cohopficity of these groups is realizable via 3-dimensional spaces.  The fundamental groups of the exteriors of certain knotted $S^2$'s in $\mathbb{R}^4$, see Yoshiawa \cite{Yoshiawa1988}, are the GBS groups $G=\langle a,b,t:a^p=b^q, t a^{\alpha}t^{-1}=b^{\beta} \rangle$, $p\beta-q\alpha = \pm 1$, and consequently these fundamental groups have finite index clones.  Given that the groups have finite index clones, it is reasonable to ask if the associated 4-manifolds have nontrivial self covers.  The fundamental groups of these knot exteriors are of independent interest -- they are finitely related deficiency 1 groups which have no free base \cite{Yoshiawa1988}.

\bigskip
Given a finite GBS graph $(\Gamma, \omega)$, the self covering of the associated GBS complex $K(\Gamma,\omega)$ provides for the geometric representability of the non-cohopficity of the GBS group $\pi(\Gamma,\omega)$ in dimension 2 within the class of aspherical 2-complexes.  Thomas \cite{Thomas2024}, investigates a second type of dimension 2 geometric representability for GBS groups, that of the geometric representability of algebraic $(G,2)$-complexes for an aspherical group $G$.  Note that the standard presentations for GBS groups associated to finite graphs are deficiency 1 presentations.  Consequently, the presentation complexes for these groups are aspherical and so are all such GBS groups.

\begin{definition} Given a group $G$, an \emph{algebraic (G,2)-complex} is an exact sequence
$$
A=(A_2 \overset{\partial_2}{\longrightarrow} A_1 \overset{\partial_1}{\longrightarrow} A_0 \overset{\epsilon}{\longrightarrow} \mathbb{Z} \rightarrow 0)
$$
\noindent where the $A_i$ are $\mathbb{Z}G$-modules.
\end{definition}

\noindent The modules of interest are non-free modules $M$, which are \emph{stably free}, that is, modules $M$ for which there is free module $F$ such that $M\oplus F$ is free.  Then, according to Thomas,

\begin{definition}
Given an aspherical group $G$, a stably free $\mathbb{Z}G$-module $M$ is \emph{geometrically realizable} if and only if there is a finite 2-complex $K$ such that $\pi_1(K)\cong G$ and $\pi_2(K)\cong M$.
\end{definition}

The GBS groups of interest in this paper are $T(10,15) = \pi([x^{10},y^{15}]) = \langle x,y : x^{10}=y^{15}\rangle$ and $BS(m,m\pm 1) = \pi([x^m,x^{m\pm1}]) = \langle t,x : t^{-1} x^m t^{-1} = x^{m \pm 1}\rangle$.  However, the GBS complexes associated to the GBS graphs which define $T(10,15)$ and $BS(m,m\pm1)$ and, as noted above, their presentation complexes for these groups are aspherical.  They therefore have trivial second homotopy group. Consequently, 2-complexes of other sorts must be used to satisfy Thomas' geometric representability criteria.

For $T(10,15)$, he begins with the following deficiency 0 presentations

$$
\begin{array}{cc}
  T(10,15) \cong \langle x, y : x^{10} = y^{15}, 1 \rangle \\ [7pt]
   T(10,15) \cong \langle x^{10}=y^{15}, a^2=b^3, x^{15}=a^3, y^{20}=b^4 \rangle.
\end{array}
$$

\noindent The non-aspherical 2-complex $S^2\lor J$ associated to the first presentation is obtained by identifying a point in $S^2$ with a point in the presentation complex $J$ for $\pi([x^{10},x^{15}])$.\footnote{For the \underline{GBS} complex $K=K([x^{10},y^{15}]$, there is a degree $d$ covering space of $K\lor S^2$ induced by a degree $d$ self cover of $K$.  The total space of this cover is obtained by attaching $d$ copies of $S^2$ to the GBS $K=K([x^{10},x^{15}])$.  It is, consequently, not a self cover because $d\geq 2$.} Then, by \cite{Thomas2024}, Theorem D, there exist non-free stably free $\mathbb{Z}BS(m,m\pm1)$ modules.  Question 1, \cite{Thomas2024} asks if these modules are geometrically realizable.

\bigskip
The Tinsley construction, as described above, provides high dimensional examples whose fundamental groups are examples of groups with proper finite index clones.  Daverman introduces the construction in the context of his interest in Coram and Duval's notion of \emph{approximate fibration} \cite{CoramDuval1977}.  It turns out that when $M$ is a manifold, the existence of a $\mathbb{Z}_d$-regular self cover, $2 \leq d < \infty$, of $M$ is an obstruction to it being the (approximate) fibre in an approximate fibration.  In particular, by  Corollary 5.9 in \cite{Daverman1993HyperhopfianGA}, because $\mathbb{RP}^3 \# \mathbb{RP}^3$ has a $\mathbb{Z}_2$-regular self cover, it \underline{cannot} be a \emph{codimension 2} fibrator. This result shows that $\mathbb{Z}_2\star \mathbb{Z}_2$ is the anomaly among free products of finite cyclic groups which are fundamental groups of closed 3-manifolds.

Daverman's work on approximate fibrations generated two questions, still of interest at this writing, which he asked at the 2004 June Workshop in Geometric Topology \cite{DavermanJuneWork2004}.  (1) Must a degree 1 continuous map from a manifold to itself induce an isomorphism of fundamental groups? (Recall that any such map must be an epimorphism at the $\pi_1$ level.)  (2) If so, must the map be a homotopy equivalence?  He is also interested in the relationship between manifolds which have $G$-regular self covers for finite non-cyclic and non-abelian groups $G$ and accordingly non-cohopfian groups $G$ with normal finite index clones $N$ such that $G/N$ is a finite non-abelian group.  A group worth consideration in this context is suggested in Example 2.1 of van Limbeck \cite{vanLimbeck2021}.

\bigskip
Finally, we return again briefly to the non-cohopfian condition and self covering properties for 3-manifolds and their fundamental group.  A more detailed look at this situation will be forth coming.

Because of the product structure, when $F$ is a 2-dimensional surface, the 3-manifold $S^1\times F$ has nontrivial self covers coming from those of the $S^1$ factor and, as noted in the introduction, if $h:F \rightarrow F$ is a periodic homeomorphism, the locally trivial $F$-bundle over the circle $M = [0,1] \times F/h$ has a $\mathbb{Z}_{k+1}$-regular self cover.  When $F$ is a 2-torus, Klein bottle, annulus, or Mobius band, $F\times I$ and $F\times S^1$ also have non-trivial self covers induced by those of $F$.  The twisted I-bundle $B$ over the Klein bottle and its double, the closed 3-manifold obtained by gluing two copies of $B$ together along their boundries, also have non-trival self covers.  Further generalizing the results for the torus knot complements and products $S^1\times F$, $F$ a 2-dimensional surface, many Seifert fibred 3-manifolds have non-trivial self covers.

It is more surprising that the connect sum $\mathbb{RP}^3\#\mathbb{RP}^3$ has nontrivial self covers.  In fact, it has degree $d$ self covers for all $d \geq 2$, but the only one of these which is regular is the $\mathbb{Z}_2$-regular one.  The proof of this is not hard, but depends on the fact that $S^3$ is the identity with respect to connect sum of 3-manifolds.  We illustrate by constructing the $\mathbb{Z}_2$-regular self cover.

Given two 3-manifolds $M_1$, $M_2$ (perhaps with boundary), the \emph{connect sum} $M_1\#M_2$ is formed as follows.  Choose small closed 3-balls $B^3_i \subset int(M_i)$ with boundaries $\partial B^3_i = S^2_i$ which are bi-collared in $M_i$.  Delete the interiors of the $B_i$ and glue the resulting punctured 3-manifolds $M_i\setminus int(B^3_i)$ together via a homeomorphism $h:S^2_1\rightarrow S^2$.  Under this definition of the sum $\#$, the 3-sphere $S^3$ is the identity.  That is, given any 3-manifold $M^3$, $M^3\#S^3$ is homeomorphic to $M^3$.

To see that $\mathbb{RP}^3\#\mathbb{RP}^3$ has a $\mathbb{Z}_2$-regular self cover $f$, write the domain of $f$ as $\mathbb{RP}^3 \#S^3 \#\mathbb{RP}^3$ and its range as $\mathbb{RP}^3\#\mathbb{RP}^3$.  Then $f$ maps the lefthand copy of $\mathbb{RP}^3$ in the domain to the left copy of it in the range via the identity, wraps the $S^3$ summand in the domain to the righthand copy of $\mathbb{RP}^3$ in the range via the antipodal map, and maps the righthand copy of $\mathbb{RP}^3$ in the domain to the left copy of $\mathbb{RP}^3$ in the range via the identity.  As a corollary, the free product of finite groups $\pi_1(\mathbb{RP^3}\# \mathbb{RP}^3)\cong \mathbb{Z}_2\star \mathbb{Z}_2$, is non-cohopfian and this fact is geometrically representable in dimension 3.  In general, the degree $d$ self covers of $\mathbb{RP}^3\#\mathbb{RP}^3$ are given by writing  $\mathbb{RP}^3\#\mathbb{RP}^3=\mathbb{RP}^3\# \left(\overset{d-1}{\underset{i=1}{\#}}S^2\right)\# \mathbb{RP}^3$.   That these self covers are not regular is a parity argument.  Tollefson \cite{Tollefson1969} proves that $\mathbb{RP}^3\#\mathbb{RP}^3$ is the only closed, connected, non-prime 3-manifold which non-trivially covers itself.

At the $\pi_1$ level, the preceding shows the geometric representability in dimension 3 of the fact that the free product $\mathbb{Z}_2 \star \mathbb{Z}_2$ has clones of every finite index.  Note that this can also be shown by considering the one-point joint $\mathbb{RP}^2 \lor \mathbb{RP}^2$ and mimicking the above construction in dimension 2.

\section{Concluding Remarks}\label{sec:conclussions}
The processes for constructing spaces which non-trivially self cover fall into roughly 6 broad overlapping categories.  They are

\begin{enumerate}
   \item Products and their generalizations, including
    \begin{enumerate}
        \item locally trivial bundles and
        \item Seifert fibred spaces.
    \end{enumerate}
   \item Categorical processes, including
    \begin{enumerate}
        \item that of \cite{DelgadoTimm2017} and
        \item The Tinsley Construction.
    \end{enumerate}
   \item Cantor Set controlled constructions.
   \item Graphs of spaces.
   \item Inverse limit constructions.
   \item Complex dynamical systems base constructions.
\end{enumerate}

\noindent  These processes can be combined to produce additional examples of metric spaces for which it is reasonable to ask whether or not the resulting spaces will non-trivially self cover.  The work of Eda et al \cite{EdaEtAl2005} on inverse limits of inverse sequences of 2-tori and Example \ref{ex:GBS+CantorControlled} above illustrate that sometimes combinations of these methods will provide spaces which non-trivially self cover, while other times they will not.  Further investigation is warranted.

\bigskip
In closing, if your favorite example of a space which non-trivially self covers has been omitted in this discussion, or you know of a another method for constructing spaces, metric or not, or even non-Hausdorff, which non-trivially self cover, not mentioned above, let me know.  I would like to continue to collect these spaces and process, and to included them in a sequel to this paper.

\bigskip

\noindent
\large{\textbf{Acknowledgements.}}
Thanks to to participants in the Dynamics and Continuum Theory Session of the 2025 Summer Conference on Topology and Its Applications for their discussion of this material. Thanks to Bob Daverman for our many discussion of self covering phenomena at many conferences over many years.  All of these conversations are much appreciated and helped make this a better paper.

\printbibliography

\end{document}